\documentclass[a4paper,12pt]{amsart}
\date{D\'ecembre 2001}
\usepackage[francais]{babel}
\usepackage{amssymb}
\input xypic
\xyoption{all}
\makeindex

\def\K0{{\mathrm{K}}}
\def\isom{{\mapright{\sim}}}

\def\Bb{{\mathbf{B}}}
\def\Cb{{\mathbf{C}}}
\def\Db{{\mathbf{D}}}
\def\Fb{{\mathbf{F}}}
\def\Gb{{\mathbf{G}}}
\def\Hb{{\mathbf{H}}}
\def\Kb{{\mathbf{K}}}
\def\Lb{{\mathbf{L}}}
\def\Pb{{\mathbf{P}}}
\def\Qb{{\mathbf{Q}}}
\def\Sb{{\mathbf{S}}}
\def\Tb{{\mathbf{T}}}
\def\Ub{{\mathbf{U}}}
\def\Vb{{\mathbf{V}}}
\def\Xb{{\mathbf{X}}}
\def\Yb{{\mathbf{Y}}}
\def\Zb{{\mathbf{Z}}}
\def\ab{{\mathbf{a}}}
\def\bb{{\mathbf{b}}}
\def\mb{{\mathbf{m}}}
\def\nb{{\mathbf{n}}}
\def\sb{{\mathbf{s}}}
\def\vb{{\mathbf{v}}}
\def\wb{{\mathbf{w}}}
\def\xb{{\mathbf{x}}}
\def\yb{{\mathbf{y}}}

\def\pib{{\mathbf{\pi}}}
\def\FM{{\mathbf{F}}}
\def\NM{{\mathbf{N}}}
\def\QM{{\mathbf{Q}}}
\def\ZM{{\mathbf{Z}}}
\def\mG{{\mathfrak{m}}}
\def\lG{{\mathfrak{l}}}
\def\AC{{\mathcal{A}}}
\def\CC{{\mathcal{C}}}
\def\EC{{\mathcal{E}}}
\def\FC{{\mathcal{F}}}
\def\GC{{\mathcal{G}}}
\def\HC{{\mathcal{H}}}
\def\IC{{\mathcal{I}}}
\def\OC{{\mathcal{O}}}
\def\RC{{\mathcal{R}}}
\def\SC{{\mathcal{S}}}
\def\UC{{\mathcal{U}}}
\def\XC{{\mathcal{X}}}

\newtheorem{theo}{Th\'eor\`eme}[section]
\newtheorem{prop}[theo]{Proposition}
\newtheorem{lem}[theo]{Lemme}
\newtheorem{coro}[theo]{Corollaire}

\theoremstyle{definition}
\newtheorem{rem}[theo]{Remarque}

\def\Hom{\mathop{\mathrm{Hom}}\nolimits}
\def\Id{\mathop{\mathrm{Id}}\nolimits}
\def\Im{\mathop{\mathrm{Im}}\nolimits}
\def\Ind{\mathop{\mathrm{Ind}}\nolimits}
\def\Irr{\mathop{\mathrm{Irr}}\nolimits}

\def\Res{\mathop{\mathrm{Res}}\nolimits}
\def\Spec{\mathop{\mathrm{Spec}}\nolimits}

\def\tete#1{\par\leavevmode\makebox[0.7cm]{$(\mathrm{#1})$}}

\def\longto{\longrightarrow}
\def\injto{\hookrightarrow}

\def\mapright#1{\hspace{0.2em}\smash{
     \mathop{\rightarrow}\limits^{\SS#1}}\hspace{0.2em}}

\def\longmapright#1{\hspace{0.3em}\smash{
     \mathop{\longrightarrow}\limits^{#1}}\hspace{0.3em}}
\def\longtrait#1{\hspace{0.3em}{\frac{~ \SS{#1} ~}{~}}\hspace{0.3em}}

\def\fonction#1#2#3#4#5{\begin{array}{rccc}
{#1} : & {#2} & \longto & {#3} \\
& {#4} & \longmapsto & {#5} 
\end{array}}

\def\fonctio#1#2#3#4{\begin{array}{ccc}
{#1} & \longto & {#2} \\
{#3} & \longmapsto & {#4} 
\end{array}}

\def\incl{\hspace{0.05cm}{\subset}\hspace{0.05cm}}
\def\vide{\emptyset}

\def\fq{\FM_q}
\def\fp{\FM_p}

\def\DS{\displaystyle}
\def\SS{\scriptstyle}

\def\lexp#1#2{\kern\scriptspace\vphantom{#2}^{#1}\kern-\scriptspace#2}

\def\micro{\hspace{0.05cm}}
\def\cad{c'est-\`a-dire }
\def\car{caract\`ere }
\def\cars{caract\`eres }
\def\ele{\'el\'ement }
\def\eles{\'el\'ements }

\def\irrs{irr\'eductibles }
\def\mor{morphisme }

\def\endo{endomorphisme }

\def\iso{isomorphisme }
\def\isos{isomorphismes }
\def\para{sous-groupe parabolique }

\def\resp{respectivement }
\def\ssi{si et seulement si }
\def\tor{tore maximal }

\def\can{{\mathrm{can}}}

\def\equat{\refstepcounter{theo}$$~}
\def\endequat{\leqno{\boldsymbol{(\arabic{section}.\arabic{theo})}}~$$}

\def\sem{{\mathrm{sem}}}
\def\nablab{{\boldsymbol{\nabla}}}

\def\parf{\operatorname{\!-parf}\nolimits}
\def\proj{\operatorname{\!-proj}\nolimits}
\def\RRR{{\mathrm{R}}}

\def\LLL{{\mathrm{L}}}
\def\Mod{\operatorname{\!-mod}\nolimits}

\def\equivb#1{\hspace{0.1cm}{\boldsymbol{\equiv}_{#1}}\hspace{0.1cm}}

\def\opp{{\mathrm{op}}}

\begin{document}

\title{Cat\'egories d\'eriv\'ees et vari\'et\'es de Deligne-Lusztig}
\author{C\'edric Bonnaf\'e et Rapha\"el Rouquier}
\address{C\'edric Bonnaf\'e~:
        Universit\'e de Franche-Comt\'e, D\'epartement de Math\'ematiques 
	(CNRS UMR 6623), 16 Route de Gray, 25000 Besan\c{c}on, FRANCE.\newline
         {\it E-mail~:~}{\tt bonnafe@math.univ-fcomte.fr}}
\address{Rapha\"el Rouquier~:
        UFR de Math\'ematiques et Institut de Math\'ematiques de Jussieu
        (CNRS UMR 7586), Universit\'e Paris 7, 2 place Jussieu,
        75251 Paris Cedex 05, FRANCE.\newline
        {\it E-mail~:~}{\tt rouquier@math.jussieu.fr}}

\maketitle
\tableofcontents

\section{Introduction}

Les classes de conjugaison d'un groupe $G$ r\'eductif connexe sur un corps
fini se d\'ecrivent \`a partir d'une classe d'\'el\'ements semi-simples et
une classe d'\'el\'ements unipotents du centralisateur d'un \'el\'ement
semi-simple (d\'ecomposition de Jordan). Deligne et Lusztig ont
construit une d\'ecomposition analogue des
caract\`eres irr\'eductibles de $G$~: ils sont param\'etr\'es par
une classe semi-simple du groupe dual et un caract\`ere unipotent du
centralisateur d'un \'el\'ement de cette classe.
Nous d\'emontrons ici un r\'esultat du m\^eme type pour des repr\'esentations
modulaires de $G$.

\smallskip
Choisissons un nombre premier $\ell$ diff\'erent de la
caract\'eristique du corps de d\'efinition du groupe.
Brou\'e et Michel ont montr\'e que la d\'ecomposition en $\ell$-blocs
des caract\`eres raffine la d\'ecomposition suivant la
$\ell'$-partie de la classe semi-simple associ\'ee \`a un caract\`ere.
Brou\'e \cite{B} a conjectur\'e que les foncteurs d'induction de
Deligne-Lusztig devaient produire une \'equivalence de Morita entre
la r\'eunion des blocs associ\'es \`a un $\ell'$-\'el\'ement semi-simple
$s$ du groupe dual et la r\'eunion des blocs unipotents du groupe
dual du centralisateur de $s$ (lorsque ce dernier est un sous-groupe de Levi).

La motivation de ce travail est la preuve de cette conjecture
sur la d\'ecomposition de Jordan des blocs. Nous d\'emontrons
qu'un bloc d'un groupe r\'eductif connexe sur un corps fini (en
caract\'eristique transverse) est
Morita-\'equivalent \`a un bloc quasi-isol\'e d'un 
sous-groupe de Levi.

Comme l'avait montr\'e Brou\'e, le probl\`eme est g\'eom\'etrique.
On dispose de foncteurs d'induction de Deligne-Lusztig, qui envoient
repr\'esentations sur complexes de repr\'esentations. Le probl\`eme
est de d\'emontrer que, sous certaines hypoth\`eses, les complexes
qui interviennent sont concentr\'es en un degr\'e. Brou\'e a ramen\'e
le probl\`eme \`a \'etablir la concentration en degr\'e moiti\'e de
la cohomologie de certains faisceaux localement
constants sur des vari\'et\'es de Deligne-Lusztig.
Ce probl\`eme a \'et\'e r\'esolu, dans le cas de faisceaux mod\'er\'ement
ramifi\'es associ\'es \`a des caract\`eres de tores, par Deligne et Lusztig
\cite{delu} et notre travail consiste \`a r\'esoudre ce probl\`eme en
g\'en\'eral (th\'eor\`eme \ref{i 0}).

\smallskip
Le probl\`eme qui se pose est double. D'abord, les vari\'et\'es
concern\'ees $X$ sont associ\'ees \`a des vari\'et\'es de drapeaux paraboliques
et on manque d'une compactification lisse avec diviseurs \`a croisements
normaux. Ensuite, les syst\`emes locaux $\FC$ en jeu sont sauvagement
ramifi\'es. Du coup, il ne suffit plus
de montrer que le faisceau ne s'\'etend pas \`a l'infini dans une bonne
compactification pour d\'eduire que sa cohomologie est concentr\'ee
en degr\'e moiti\'e.

\smallskip
Nous parvenons \`a contourner ces difficult\'es en deux \'etapes.

\smallskip
Tout d'abord, nous introduisons des vari\'et\'es $f_i:X_i\to X$ au-dessus de
$X$ et des syst\`emes locaux mod\'er\'ement ramifi\'es
$\FC_i$ sur $X_i$ tels que $\FC$ est dans la sous-cat\'egorie triangul\'ee
de $D^b(X)$ engendr\'ee par les $\RRR(f_i)_*\FC_i$.
Nous d\'eduisons cela d'un r\'esultat pr\'esentant un int\'er\^et
ind\'ependant~: la cat\'egorie des complexes parfaits de modules
pour $G$ est engendr\'ee par les images des foncteurs de Deligne-Lusztig
(\S \ref{sec engendrement}).

\smallskip
Ensuite, une \'etude pr\'ecise de la ramification \`a l'infini des
faisceaux $\FC_i$ nous permet de conclure. Plus pr\'ecis\'ement,
on dispose d'une immersion ouverte $j_i:X_i\to Y_i$ et d'un
morphisme propre $Y_i\to\bar{X}$, o\`u $j:X\to \bar{X}$ est une
compactification de $X$ (singuli\`ere en g\'en\'eral). L'\'etude de la
ramification montre que $(\RRR(j_i)_*\FC_i)_{|f_i^{-1}(\bar{X}-X)})=0$,
ce qui suffit pour \'etablir que le morphisme canonique
$j_!\FC\to \RRR j_*\FC$ est un quasi-isomorphisme.

Pour comprendre pr\'ecis\'ement comment
les syst\`emes locaux $\FC_i$ associ\'es \`a des caract\`eres de tores se
ramifient \`a l'infini (\S \ref{section monodromie}), nous construisons
de nouvelles vari\'et\'es de type
Deligne-Lusztig (\S \ref{sub recollement}). Le recollement de
rev\^etements repose sur la construction d'isomorphismes canoniques entre
quotients de tores. C'est l'objet de la partie \S \ref{sub sous}, o\`u
nous obtenons aussi une nouvelle approche de la conjugaison rationnelle
de caract\`eres de tores.

\smallskip
Les r\'esultats principaux sont contenus dans les parties \S 10 et \S 11.
Dans les parties \S \ref{section monodromie} et \S \ref{sub f},
nous supposons que l'anneau de coefficients $\Lambda$ est un corps.
L'appendice \S \ref{appendice} regroupe quelques r\'esultats g\'eom\'etriques.


\smallskip
Le lecteur verra sans peine l'importance consid\'erable
de notre dette \`a l'\'egard de
l'article fondateur de Deligne et Lusztig \cite{delu}. Nous utilisons
aussi certaines techniques d\'evelopp\'ees par Digne et Michel \cite{dimirou}
(vari\'et\'es associ\'ees \`a des suites d'\'el\'ements du groupe de Weyl
et compactifications partielles).

\smallskip
Dans un second travail, nous aborderons les applications \`a la
th\'eorie locale des blocs (morphisme de Brauer, \'equivalences splendides,
finitude des alg\`ebres de source, engendrement de la cat\'egorie
d\'eriv\'ee, variation du sous-groupe parabolique, groupes non connexes).

\smallskip
Nous remercions Marc Cabanes et Michel Enguehard pour leurs nombreux
commentaires et suggestions.

\section{Notations g\'en\'erales}
\subsection{Groupes, mono\"{\i}des}
\label{Groupes}
Soit $G$ un groupe. Nous noterons $G^\opp$\index{g@$G^\opp$}
le groupe oppos\'e \`a $G$.
Pour $M$ un $G$-ensemble, on note $M^\opp$\index{m@$M^\opp$} le $G^\opp$-ensemble,
restriction de $M$ \`a travers l'isomorphisme
$G^\opp\isom G,\ g\mapsto g^{-1}$.

Pour $g \in G$, nous noterons $(g)$ (ou $(g)_G$ s'il y a confusion possible)
la classe de conjugaison de $g$ dans $G$. Si $g$ et $g'$ sont deux \eles de
$G$ nous \'ecrirons
$g \sim g'$ (ou $g \sim_G g'$) pour dire que $g$ et $g'$ sont conjugu\'es
dans $G$.

 Si $X$ est une partie de $G$ stable par conjugaison, $X/\sim$
(ou $X/\sim_G$) d\'esignera l'ensemble des classes de conjugaison de $G$
contenues dans $X$. Si $\pi$ est un ensemble de nombres premiers (ou un
nombre premier) nous noterons $X_\pi$ (\resp $X_{\pi'}$) l'ensemble des
$\pi$-\eles (\resp $\pi'$-\'el\'ements) de $X$. Si $g \in G$ est d'ordre
fini, nous noterons $g_\pi$ et $g_{\pi'}$ les uniques \eles de $G_\pi$ et
$G_{\pi'}$ respectivement tels que $g=g_\pi g_{\pi'}=g_{\pi'}g_\pi$.

Si $E$ est un ensemble, nous noterons $\Sigma(E)$\index{sig@$\Sigma(E)$}
le mono\"{\i}de libre sur
l'ensemble $E$, \cad l'ensemble des suites finies d'\eles de $E$ muni de
la loi de concat\'enation. Son \'el\'ement neutre sera
not\'e $()$. Si $\sigma$ est un endomorphisme de $E$, nous noterons
parfois ${^\sigma e}$ l'\'el\'ement $\sigma(e)$.

\subsection{Anneaux, corps\label{sub anneaux}}
Dans cet article, nous fixons une fois pour toutes deux nombres premiers
distincts $p$\index{p@$p$} et $\ell$\index{l@$\ell$}.
Nous noterons $\FM$\index{f@$\FM$} une cl\^oture alg\'ebrique du corps
fini \`a $p$ \eles $\fp$. Nous fixons une puissance $q$\index{q@$q$}
de $p$ et notons $\fq$
le sous-corps de $\FM$ \`a $q$ \'el\'ements. Nous appellerons vari\'et\'e (\resp
groupe alg\'ebrique) une vari\'et\'e alg\'ebrique quasi-projective d\'efinie
(\resp un groupe alg\'ebrique d\'efini) sur le corps $\FM$.

Nous fixons aussi une extension alg\'ebrique $K$\index{k@$K$} du corps
$\ell$-adique $\QM_\ell$.
Soit $R$\index{r@$R$} son anneau d'entiers sur $\ZM_\ell$ et soit
$k$\index{k@$k$} le corps r\'esiduel de
$R$~: c'est une extension alg\'ebrique du corps fini \`a $\ell$ \eles $\FM_\ell$.
Nous notons $\lG$\index{l@$\lG$} l'id\'eal maximal de $R$.
Nous supposerons que le corps $K$ est ``assez gros'', \cad que, pour tout groupe
fini $H$ rencontr\'e dans cet article, $K$ contient les racines $e$-i\`emes de
l'unit\'e, o\`u $e$ est l'exposant de $H$.

Nous fixons une fois pour toutes un \iso
        $$\imath : (\QM/\ZM)_{p'} \longmapright{\sim} \FM^\times
	\index{im@$\imath$}$$
ainsi qu'un \mor injectif
        $$\jmath : (\QM/\ZM)_{p'} \injto \overline{\QM}_\ell^\times\index{jm@$\jmath$},$$
o\`u $\overline{\QM}_\ell$ est une cl\^oture alg\'ebrique de $K$. 
Soit
        $$\kappa = \jmath \circ \imath^{-1} : \FM^\times \injto 
	\overline{\QM}_\ell^\times
	\index{ka@$\kappa$}.$$

Dans tout cet article, $\Lambda$ sera l'un des anneaux $K$ ou
$R/\lG^n$ ($n\ge 0$).

\subsection{Cat\'egories d\'eriv\'ees}
Si $\AC$ est une cat\'egorie additive, nous noterons $\K0(\AC)$ son groupe 
de Grothendieck et $K^b(\AC)$ la cat\'egorie homotopique des complexes 
born\'es d'objets de $\AC$.

Si  $\CC$  est  une  cat\'egorie ab\'elienne, nous noterons 
$\K0(\CC)$ son groupe de Grothendieck,
$D^b(\CC)$ (respectivement $D^-(\CC)$) la cat\'egorie d\'eriv\'ee 
des complexes \`a cohomologie born\'ee (respectivement born\'ee
sup\'erieurement) et
$\CC\parf$\index{parf@$\CC\parf$} la cat\'egorie de ses complexes parfaits~;
$\K0(\CC)$\index{k0@$\K0(\CC)$} est
aussi le groupe de Grothendieck  de  la  cat\'egorie triangul\'ee  $D^b(\CC)$.
Si $E$ est un ensemble d'objets d'une cat\'egorie triangul\'ee $D$,
nous appellerons {\it sous-cat\'egorie
de  $D$ engendr\'ee par  $E$}  la plus petite sous-cat\'egorie
triangul\'ee pleine de $D$ stable par facteurs directs et contenant $E$.
Nous identifierons $\CC$  avec  la  sous-cat\'egorie pleine de  $D^b(\CC)$
des complexes concentr\'es en degr\'e $0$ via le foncteur canonique.

Lorsque $\CC$ a assez de projectifs, \'etant donn\'es $C$ et $C'$ deux objets
de $D^b(\CC)$, nous noterons 
$\RRR\Hom^\bullet(C,C')$ le complexe total associ\'e au complexe
double des homomorphismes d'une r\'esolution projective de $C$ vers $C'$.

\subsection{Alg\`ebres}
Soient $A$ et $B$ deux $\Lambda$-alg\`ebres.     Nous noterons $A^\opp$
l'alg\`ebre oppos\'ee \`a $A$.   Par un $A$-module, nous entendrons un module
\`a gauche pour $A$. Tous les modules consid\'er\'es seront de type fini. Nous
noterons $A\Mod$ la cat\'egorie des $A$-modules.  On identifie la
cat\'egorie des $(A,B)$-bimodules \`a celle
des $(A\otimes_\Lambda B^\opp)$-modules.
Nous noterons $\K0(A)$, $D^b(A)$, $D^-(A)$ et
$A\parf$ les constructions obtenues pour $\CC=A\Mod$.

Pour $M$ un complexe born\'e de $A$-modules, nous d\'efinissons le complexe de
$A$-modules \`a droite
$M^*=\RRR\Hom_\Lambda^\bullet(M,\Lambda)$.

Nous noterons $A\proj$ la sous-cat\'egorie additive pleine de $A\Mod$
form\'ee par les $A$-modules projectifs.

\subsection{Groupes finis}
Si  $G$ est un groupe fini,  nous  d\'esignerons  par 
$\Irr_K G$\index{irr@$\Irr_K G$} l'ensemble de ses
\cars \irrs sur  $K$.  Nous  identifierons  alors $\K0(KG)$ avec $\ZM \Irr_K G$,
le  $\ZM$-module libre  de base  $\Irr_K G$.  Le produit scalaire sur  $\K0(KG)$
faisant de  $\Irr_K G$  une  base  orthonormale sera not\'e
$\langle \cdot,\cdot \rangle_G$.
Si  $\theta : G \to \Lambda^\times$ est un \car lin\'eaire de $G$,  nous noterons
$\Lambda_\theta$\index{lam@$\Lambda_\theta$}
le $\Lambda G$-module sur lequel un \ele $g \in G$ agit par multiplication par $\theta(g)$,
et $e_\theta$\index{eth@$e_\theta$}
(ou $e_{\Lambda,\theta}$, ou $e_{\Lambda,\theta}^G$ s'il est n\'ecessaire de pr\'eciser)
l'idempotent primitif central de $\Lambda G$ qui agit comme l'identit\'e sur
$\Lambda_\theta$.

\subsection{Groupes alg\'ebriques}
Soit $\Hb$ un groupe alg\'ebrique. On d\'esignera par
$\Hb^\circ$\index{hci@$\Hb^\circ$} sa composante neutre
et $\Hb_\sem$\index{hsem@$\Hb_\sem$}
l'ensemble de ses \eles semi-simples. Si $h \in \Hb$, nous
noterons
$C_\Hb^\circ(h)$ la composante neutre de son  centralisateur dans $\Hb$.
Si $\Hb$ est ab\'elien et si $F$ est un
endomorphisme de $\Hb$, alors nous notons
$N_{F^n/F} : \Hb \to \Hb$\index{no@$N_{F^n/F}$},
$h \mapsto h~\lexp{F}{h}\dots \lexp{F^{n-1}}{h}$ le \mor {\it norme de $F^n$
\`a $F$} (o\`u $n \in \NM^*$).

\subsection{Faisceaux}
Soit $\Xb$ une vari\'et\'e et soit $A$ une $\Lambda$-alg\`ebre.
Nous noterons $D^?_A(\Xb)$\index{d?@$D^?_A(\Xb)$} la sous-cat\'egorie
pleine de la cat\'egorie d\'eriv\'ee de la cat\'egorie des  faisceaux
de $A$-modules sur $X$
form\'ee des complexes \`a cohomologie constructible, born\'ee si
$?=b$, born\'ee sup\'erieurement si $?=-$.

Notons $\pi_\Xb : \Xb \to \Spec\FM$\index{pi@$\pi_\Xb$} le \mor canonique.
Un $A$-module $M$ sera identifi\'e avec le  $A$-faisceau constant sur
$\Spec\FM$
associ\'e \`a  $M$.  Nous  posons  alors
$M_\Xb = \pi_\Xb^* M$\index{mx@$M_\Xb$}~:  c'est le
$A$-faisceau constant sur $\Xb$ associ\'e au $A$-module $M$.  De m\^eme, un
complexe de $A$-modules $C$ sera identifi\'e avec  le complexe de faisceaux
constants sur $\Spec\FM$ correspondant et nous noterons $C_\Xb=\pi_\Xb^* C$.

Nous noterons
$\RRR\Gamma(\Xb,C)=\RRR(\pi_{\Xb})_*(C)$
\index{rg@$\RRR\Gamma(\Xb,C)$} (\resp
$\RRR\Gamma_c(\Xb,C)=\RRR(\pi_{\Xb})_!(C)$\index{rgc@$\RRR\Gamma_c(\Xb,C)$}) le
complexe de cohomologie (\resp de cohomologie \`a support compact)
d'un complexe $C\in D^b_A(\Xb)$. Nous \'ecrirons $\RRR\Gamma(\Xb)$ 
(\resp $\RRR\Gamma_c(\Xb)$) pour $\RRR\Gamma(\Xb,\Lambda_\Xb)$
(\resp $\RRR\Gamma_c(\Xb,\Lambda_\Xb)$).
Pour $G$ un groupe fini agissant sur $\Xb$ de telle sorte que les stabilisateurs
des points de $\Xb$ sont des sous-groupes de $G$ d'ordre inversible dans
$\Lambda$, alors $\RRR\Gamma_c(\Xb)$ est un complexe parfait de
$\Lambda G$-modules
\cite[\S 3.8]{delu}.

\section{Vari\'et\'es et groupes finis\label{GvarietesH}}~

Les vari\'et\'es de Deligne-Lusztig sont munies d'actions \`a gauche et \`a
droite de groupes r\'eductifs finis. Leur cohomologie permet alors de d\'efinir
des foncteurs entre cat\'egories de repr\'esentations. Nous rappelons ici
dans un cadre g\'en\'eral comment d\'efinir ces foncteurs ainsi que quelques-unes
de leurs propri\'et\'es \'el\'ementaires.

\subsection{D\'efinitions\label{sousousous}}
Nous fixons trois groupes finis $G$, $H$ et $K$. Une {\it $G$-vari\'et\'e}
(\resp une {\it vari\'et\'e-$G$}) est une vari\'et\'e sur laquelle $G$ agit \`a
gauche (\resp \`a droite). Une {\it $G$-vari\'et\'e-$H$} est \`a la fois une
$G$-vari\'et\'e et une vari\'et\'e-$H$, les actions de $G$ et $H$ commutant.
Si $\Xb$ est une vari\'et\'e-$G$ et si $\Yb$ est une $G$-vari\'et\'e, nous
noterons $\Xb \times_G \Yb$ la vari\'et\'e quotient de $\Xb \times \Yb$ par
$G$, le groupe $G$ agissant \`a gauche sur $\Xb \times \Yb$ par
        $$\fonctio{G \times (\Xb \times \Yb)}{\Xb \times
                                     \Yb}{(g,(x,y))}{(xg^{-1},gy).}$$
Une $G$-vari\'et\'e, ou une vari\'et\'e-$G$, est dite
{\it r\'eguli\`ere} si le groupe $G$ agit librement.

\subsection{Foncteurs\label{sub foncteur}}
Soit $\Xb$ une $G$-vari\'et\'e-$H$.  Les complexes $\RRR\Gamma_c(\Xb)$ et
$\RRR\Gamma(\Xb)$ sont des  complexes de $(\Lambda G,\Lambda H)$-bimodules.
Ils induisent donc deux foncteurs
        $$\fonction{\RC_H^G(\Xb)}{D^-(\Lambda H)}{D^-(\Lambda G)}{C}{
        \RRR\Gamma_c(\Xb) \otimes_{\Lambda H}^\LLL C}\index{rh@$\RC_H^G(\Xb)$}$$
        $$\fonction{\SC_H^G(\Xb)}{D^-(\Lambda H)}{D^-(\Lambda G)}{C}{
        \RRR\Gamma(\Xb) \otimes_{\Lambda H}^\LLL C.}\leqno{\mathrm{et}}
        \index{sh@$\SC_H^G(\Xb)$}$$
On consid\'erera aussi le foncteur
        $$\fonction{\FC_H^\Xb}{D^-(\Lambda H)}{D^-_\Lambda(\Xb/H)}{M}{\pi_* 
	\Lambda_{\Xb}
        \otimes_{\Lambda H}^\LLL M_{\Xb/H}.}\index{fh@$\FC_H^\Xb$}$$
o\`u $\pi:\Xb\to \Xb/H$ est le morphisme quotient.

\bigskip

\begin{rem}
\label{systeme}
Supposons que $\Xb$ est une vari\'et\'e-$H$
r\'eguli\`ere. Alors,
le foncteur $\FC_H^\Xb$ envoie un $\Lambda H$-module sur un $\Lambda$-syst\`eme local
sur la vari\'et\'e $\Xb/H$. En outre, les foncteurs
$\RC_H^G(\Xb)$, $\SC_H^G(\Xb)$ et $\FC_H^\Xb$ se restreignent en
des foncteurs entre cat\'egories d\'eriv\'ees born\'ees.
\end{rem}

\begin{lem}\label{lemme foncteur}
On a des isomorphismes canoniques de foncteurs
de $D^b(\Lambda H)$ vers $D^-(\Lambda G)$,
$$\RRR(\pi_{\Xb/H})_!\circ \FC_H^\Xb\simeq \RC_H^G(\Xb)$$
$$\RRR(\pi_{\Xb/H})_*\circ \FC_H^\Xb\simeq \SC_H^G(\Xb).$$

Si les stabilisateurs de points de $\Xb$ sous
l'action de $H$ sont d'ordre inversible dans $\Lambda$, 
alors, pour $M$ un $\Lambda H$-module, on a un isomorphisme
canonique
$$\FC_H^\Xb M\simeq \pi_* \Lambda_{\Xb}\otimes_{\Lambda H} M_{\Xb/H}$$
\end{lem}

\bigskip

\begin{proof} Puisque $\RRR(\pi_{\Xb/H})_! \circ
\pi_*\simeq \RRR(\pi_{\Xb/H})_! \circ \RRR\pi_!\simeq \RRR(\pi_\Xb)_!$,
(noter que $\pi_*=\pi_!$ est un foncteur exact car $\pi$ est un \mor fini),
la commutativit\'e du
premier diagramme est \'equivalente \`a la propri\'et\'e suivante~:
        $$\RRR(\pi_{\Xb/H})_!(\pi_* \Lambda_\Xb \otimes_{\Lambda H}^\LLL \pi_{\Xb/H}^* -)
        \simeq (\RRR(\pi_{\Xb/H})_!\pi_* \Lambda_\Xb )\otimes_{\Lambda H}^\LLL -.$$
Cela r\'esulte du lemme \ref{!}.

La commutativit\'e du deuxi\`eme diagramme se d\'emontre de la m\^eme
mani\`ere mais en utilisant cette fois le lemme \ref{*}
au lieu du lemme \ref{!}.

Pour la derni\`ere partie du lemme, on utilise que
les fibres de $\pi_* \Lambda_{\Xb}$ sont des $\Lambda H$-modules projectifs.
\end{proof}

\subsection{Composition}\label{Composition}
Soit maintenant $\Yb$ une $H$-vari\'et\'e-$K$. La vari\'et\'e $\Xb \times_H \Yb$
est une $G$-vari\'et\'e-$K$ et la formule de K\"unneth montre que
$$\RRR\Gamma_c(\Xb) \otimes^\LLL_\Lambda \RRR\Gamma_c(\Yb) \simeq
\RRR\Gamma_c(\Xb \times \Yb).$$

Supposons jusqu'\`a la fin du \S \ref{Composition}
que les stabilisateurs des points de $\Xb\times\Yb$
sous l'action diagonale de $H$ sont d'ordre inversible dans $\Lambda$.
Alors,
$\RRR\Gamma_c(\Xb\times\Yb)\otimes^\LLL_{\Lambda H}\Lambda
\simeq \RRR\Gamma_c(\Xb \times_H \Yb)$ (cf lemme \ref{lemme foncteur}).
On en d\'eduit, apr\`es application du foncteur
$-\otimes^\LLL_{\Lambda H}\Lambda$ \`a la formule de K\"unneth, que
\equat\label{encore}
\RRR\Gamma_c(\Xb) \otimes_{\Lambda H}^\LLL \RRR\Gamma_c(\Yb) \simeq
\RRR\Gamma_c(\Xb \times_H \Yb).
\endequat
dans $D^b(\Lambda G\otimes \Lambda K)$. En particulier, on obtient
\equat\label{composition}
\RC_H^G(\Xb) \circ \RC_K^H(\Yb) \simeq \RC_K^G(\Xb \times_H \Yb).
\endequat

\section{Pr\'eliminaires sur les groupes r\'eductifs finis\label{GF}}~
\subsection{Groupes r\'eductifs}
Dor\'enavant, et ce jusqu'\`a la fin de cet article, nous fixons un groupe
r\'eductif connexe $\Gb$\index{g@$\Gb$}
muni d'une isog\'enie $F : \Gb \to \Gb$\index{f@$F$}
dont une puissance $F^\delta$\index{de@$\delta$}
est l'\endo de Frobenius sur $\Gb$ associ\'e
\`a une structure rationnelle sur $\fq$. Il est
\`a noter que les entiers $q$ et $\delta$ ne sont pas uniquement d\'etermin\'es
par la donn\'ee de $\Gb$ et $F$. Cependant, le r\'eel positif $q^{1/\delta}$ l'est.
Soit $\Bb$ un sous-groupe de Borel $F$-stable de $\Gb$ et soit $\Tb$ un \tor 
$F$-stable de $\Bb$. Le radical unipotent de $\Bb$ sera not\'e $\Ub$. On note
$(\Gb^*,\Tb^*,F^*)$\index{g@$\Gb^*$}\index{t@$\Tb^*$}\index{f@$F^*$}
un triplet dual de $(\Gb,\Tb,F)$
\cite[d\'efinition 5.21]{delu}.

\medskip

Le sujet g\'en\'eral de cet article est l'\'etude des repr\'esentations
modulaires du groupe fini $\Gb^F$. Pour s'y pr\'eparer, nous fixons dans
ce \S \ref{GF} quelques notations (racines, groupe de Weyl...),
avant d'\'etablir quelques r\'esultats combinatoires sur les tores
$F$-stables de $\Gb$ (cf \S\ref{sub sous}). Nous
introduisons tout d'abord
deux groupes r\'eductifs connexes, du m\^eme type que $\Gb$,
qui joueront le r\^ole d'auxiliaires techniques. Ils seront utiles dans
certaines constructions (s\'eries rationnelles, recollement de vari\'et\'es
de Deligne-Lusztig...) et leur usage facilitera certaines d\'emonstrations.

\bigskip

\begin{rem}\label{hf}
Si $\Hb$ est un sous-groupe ferm\'e connexe distingu\'e de $\Gb$ (mais 
non n\'ecessairement 
$F$-stable), et si on pose $\Kb=\{g \in \Gb~|~g^{-1}F(g) \in \Hb\}$, alors 
il est imm\'ediat que $\Kb$ est un sous-groupe ferm\'e de $\Gb$ contenant $\Gb^F$. 
De plus, 
$$\Kb=\Gb^F\cdot\Kb^\circ.$$
En effet, l'application de Lang $\Gb \to \Gb$, $g \mapsto g^{-1}F(g)$ est 
un morphisme fini surjectif, donc l'image de $\Kb^\circ$ par cette application 
est une sous-vari\'et\'e irr\'eductible ferm\'ee de dimension $\dim \Kb^\circ
=\dim \Kb=\dim \Hb$ de $\Hb$. Or, $\Hb$ est ici suppos\'e connexe, 
donc l'image de $\Kb^\circ$ est $\Hb$. Le r\'esultat d\'ecoule 
alors du fait que les fibres de l'application de Lang sont des 
$\Gb^F$-orbites (pour l'action de $\Gb^F$ par multiplication \`a gauche).
\end{rem}

\bigskip

\subsubsection{}
On sait \cite[preuve du corollaire 5.18]{delu} qu'il existe un groupe
r\'eductif connexe $\tilde{\Gb}$\index{gt@$\tilde{\Gb}$}
muni d'une isog\'enie $F : \tilde{\Gb} \to \tilde{\Gb}$ telle
que $F^\delta$ soit l'\endo de Frobenius sur $\tilde{\Gb}$ associ\'e \`a une 
$\fq$-structure et v\'erifiant les conditions suivantes~:

\smallskip

\tete{1} {\it $\Gb$ est un sous-groupe ferm\'e $F$-stable de $\tilde{\Gb}$, contenant
le groupe d\'eriv\'e de $\tilde{\Gb}$~;}

\smallskip

\tete{2} {\it le centre $Z(\tilde{\Gb})$ de $\tilde{\Gb}$ est connexe.}

\smallskip

\noindent On a alors $\tilde{\Gb}=\Gb\cdot Z(\tilde{\Gb})$.
Nous noterons $\tilde{\Tb}$ (\resp $\tilde{\Bb}$) l'unique \tor (\resp
sous-groupe de Borel) de $\tilde{\Gb}$ contenant $\Tb$ (\resp $\Bb$). D'autre part,
nous noterons $(\tilde{\Gb}^*,\tilde{\Tb}^*,F^*)$ un triplet dual de 
$(\tilde{\Gb},\tilde{\Tb},F)$.
L'inclusion $\Gb \injto \tilde{\Gb}$ induit un \mor surjectif de groupes alg\'ebriques
$i^* : \tilde{\Gb}^* \to \Gb^*$ commutant avec $F$, dont le noyau est un
tore central de $\tilde{\Gb}^*$ et tel que $i^*(\tilde{\Tb}^*)=\Tb^*$. Ce
\mor $i^*$ est bien d\'efini \`a conjugaison pr\`es par un \ele de
$\Tb^{*F^*}$.

\subsubsection{}
De mani\`ere duale, il existe un groupe r\'eductif connexe $\hat{\Gb}$\index{gh@$\hat{\Gb}$}
muni d'une
isog\'enie $F : \hat{\Gb} \to \hat{\Gb}$ telle que $F^\delta$ soit l'\endo de Frobenius
sur $\hat{\Gb}$ associ\'e \`a une $\fq$-structure et un tore central
$F$-stable $\Cb$ de $\hat{\Gb}$ v\'erifiant les deux propri\'et\'es suivantes~:

\smallskip

\tete{1} $\Gb=\hat{\Gb}/\Cb$~;

\smallskip

\tete{2} {\it Toute coracine $\alpha^\vee : \FM^\times \to \hat{\Gb}$ relative \`a un
tore maximal de $\hat{\Gb}$ est injective.}

\smallskip

Notons que si le groupe d\'eriv\'e de $\Gb$ est simplement connexe, alors
toute coracine est injective.

\smallskip

\noindent Nous noterons $\rho : \hat{\Gb} \to \Gb$\index{rho@$\rho$} la projection canonique.
Nous noterons $\hat{\Tb}$ (\resp $\hat{\Bb}$) le \tor (\resp le sous-groupe de Borel)
de $\hat{\Gb}$ \'egal \`a $\rho^{-1}(\Tb)$ (\resp $\rho^{-1}(\Bb)$).

\smallskip

Plus g\'en\'eralement, nous noterons
$\tilde{?}$ (\resp $\hat{?}$) l'objet associ\'e \`a $\tilde{\Gb}$ (\resp $\hat{\Gb}$)
d\'efini de la m\^eme fa\c{c}on que l'objet $?$ associ\'e \`a $\Gb$.

\subsection{Racines, groupe de Weyl, groupe de tresses}
Soit $W$\index{w@$W$} le groupe de Weyl de $\Gb$ relatif
\`a $\Tb$, soit $n\mapsto \bar{n}$\index{nb@$\bar{n}$} le morphisme canonique
$N_\Gb(\Tb)\to W$, soit $\Phi$\index{phi@$\Phi$}
le syst\`eme de racines de $\Gb$ relatif \`a
$\Tb$ et soit $\Phi^\vee$\index{phiv@$\Phi^\vee$}
le syst\`eme de coracines. Soit $\Phi^+$\index{phiplus@$\Phi^+$}
(\resp $\Delta$\index{delta@$\Delta$}) le syst\`eme de racines positives (\resp la base)
de $\Phi$ associ\'e (\resp associ\'ee) \`a $\Bb$. Si $\alpha\in \Phi$, nous noterons
$\alpha^\vee$\index{alphavee@$\alpha^\vee$} sa coracine, 
$\Tb_{\alpha^\vee}$\index{talpha@$\Tb_{\alpha^\vee}$}le sous-tore de $\Tb$
image de $\alpha^\vee$, $\Ub_\alpha$\index{ualpha@$\Ub_\alpha$} le sous-groupe unipotent
\`a un param\`etre normalis\'e par $\Tb$ associ\'e \`a $\alpha$,
$s_\alpha$\index{salpha@$s_\alpha$}
la r\'eflexion de $W$ par rapport \`a $\alpha$ et $\Gb_\alpha$\index{galpha@$\Gb_\alpha$}
le sous-groupe
de $\Gb$ engendr\'e par $\Ub_\alpha$ et $\Ub_{-\alpha}$. Nous noterons $\phi : \Phi
\to \Phi$\index{phi@$\phi$}
la bijection d\'efinie par $\lexp{F}{\Ub_\alpha}=\Ub_{\phi(\alpha)}$
pour tout $\alpha\in \Phi$. Nous noterons
$Y(\Tb)$\index{yt@$Y(\Tb)$}
le r\'eseau des sous-groupes \`a un param\`etre
de $\Tb$.

Soit $S$\index{s@$S$}
l'ensemble des r\'eflexions simples de $W$ et soit
$\bar{S}=S \cup \{1\}$\index{sbar@$\bar{S}$}.
Si $w \in W$, nous noterons $l(w)$\index{l@$l$}
la longueur de $w$ relativement \`a $S$ (c'est aussi
le cardinal de l'ensemble $\{\alpha\in \Phi^+~|~w(\alpha) \in -\Phi^+\}$).
On notera encore $l$ la fonction sur $N_\Gb(\Tb)$ donn\'ee par
$l(\bar{n})$\index{l@$l$} pour $n\in N_\Gb(\Tb)$.
L'ordre de Bruhat sur $W$ associ\'e \`a $S$ sera not\'e $\le$\index{$\le$}.
Si $v$ et $w$
sont deux \eles de $W$, nous \'ecrirons $v<w$ pour dire que $v \le w$
et $v \not= w$.

Nous noterons $B$\index{b@$B$} le groupe de tresses associ\'e \`a $(W,S)$,
de g\'en\'erateurs $\{\sb_\alpha\}_{\alpha\in\Delta}$\index{salpha@$\sb_\alpha$}.
Soit $f : B \to W$ \index{f@$f$}
le morphisme canonique (il v\'erifie
$f(\sb_\alpha)=s_\alpha$ pour tout $\in \Delta$). Il existe alors
une et une seule application (qui n'est pas un \mor de groupes)
$\sigma : W \to B$\index{sigma@$\sigma$} v\'erifiant
\begin{itemize}
\item $\sigma(vw)=\sigma(v)\sigma(w)$ si $v$ et $w$ sont deux \eles de
$W$ tels que $l(vw)=l(v)+l(w)$ ;
\item $\sigma(s_\alpha)=\sb_\alpha$ pour tout $\alpha\in \Delta$. 
\end{itemize}

Cette application v\'erifie $f \circ \sigma = \Id_W$.

\bigskip

Nous fixons un \mor de groupes $\hat{\varphi} : B \to N_{\hat{\Gb}}(\hat{\Tb})$
v\'erifiant les propri\'et\'es suivantes~:

\smallskip

\tete{1} {\it Le diagramme
$$\diagram
B \rrto^{\DS{\hat{\varphi}}} \ddrrto_{\DS{f}} && N_{\hat{\Gb}}(\hat{\Tb}) \ddto \\
&& \\
&& W \\
\enddiagram$$
est commutatif. (Dans ce diagramme, la fl\`eche verticale est la surjection
canonique.)

\smallskip

\tete{2} Pour tout $\alpha\in \Delta$, $\hat{\varphi}(\sb_\alpha)$ appartient au 
sous-groupe de $\hat{\Gb}$ engendr\'e par $\hat{\Ub}_\alpha$ et $\hat{\Ub}_{-\alpha}$. }

\bigskip

L'existence d'un tel \mor de groupes est assur\'ee par \cite[\S 4.6]{tits}.
Gr\^ace \`a la propri\'et\'e (1), il existe un prolongement de $\hat{\varphi}$
en un \mor de groupes $\hat{\varphi} : B\hat{\Tb} \to N_{\hat{\Gb}}(\hat{\Tb})$ qui est
l'identit\'e sur $\hat{\Tb}$.
Ici, $B\hat{\Tb}$ d\'esigne le produit semi-direct $B \ltimes \hat{\Tb}$, l'action
de $B$ sur $\hat{\Tb}$ s'effectuant via le \mor $f : B \to W$ et l'action
naturelle de $W$ sur $\hat{\Tb}$. Nous noterons $\varphi : B\Tb \to N_\Gb(\Tb)$
le \mor induit par $\hat{\varphi}$.

Si $w \in W$, posons $\dot{w}=\varphi(\sigma(w))$\index{wdot@$\dot{w}$}.
Alors, $\dot{w}$ est un repr\'esentant
de $w$ dans $N_\Gb(\Tb)$. De plus, si $v$ et $w$ sont deux \eles de $W$ tels
que $l(vw)=l(v)+l(w)$, alors $\dot{x}=\dot{v}\dot{w}$, o\`u $x=vw$. Remarquons
aussi que $\dot{1}=1$. Si $n \in N_\Gb(\Tb)$ et si $w$
d\'esigne sa classe dans $W$, alors il existe un unique \ele $t \in \Tb$ tel que
$n=\dot{w} t$.
On pose alors $\dot{\sigma}(n)=\sigma(\dot{w})t \in B\Tb$\index{sigmadot@$\dot{\sigma}$}.

\subsection{Suites d'\eles de $W$}
Pour $\wb=(w_1,\dots,w_r) \in \Sigma(W)$ (cf \S \ref{Groupes}),
nous posons
$\dot{\wb}=(\dot{w}_1,\dots,\dot{w}_r)\in\Sigma(N_\Gb(\Tb))$\index{wdot@$\dot{\wb}$}.
Pour $\nb=(n_1,\dots,n_r)\in \Sigma(N_\Gb(\Tb))$, on pose
$\bar{\nb}=(\bar{n}_1,\dots,\bar{n}_r)$.

Dans la suite, nous noterons $\sigma : \Sigma(W) \to B$\index{sigma@$\sigma$} et
$l : \Sigma(W) \to (\NM,+)$\index{l@$l$}
les uniques morphismes de mono\"{\i}des prolongeant $\sigma$ et $l$ respectivement.
En d'autres termes, $\sigma(\wb)=\sigma(w_1)\dots\sigma(w_r)$ et
$l(\wb)=l(w_1)+\dots + l(w_r)$. Le nombre entier $l(\wb)$ est
appel\'e la {\it longueur} de $\wb$.
 De m\^eme, nous noterons
$\dot{\sigma} : \Sigma(N_\Gb(\Tb)) \to B\Tb$ l'unique morphisme de mono\"{\i}des
prolongeant $\dot{\sigma}$.

Si $\wb=(w_1,\dots,w_r)$ et $\vb=(v_1,\dots,v_{r'})$
sont deux \eles de $\Sigma(W)$, nous dirons que $\vb \le \wb$\index{$\le$}
si $r=r'$
et $v_i \le w_i$ pour tout $1 \le i \le r$.

Nous noterons encore $\wb$ (respectivement $\nb$) l'automorphisme de $\Tb$
induit par la conjugaison par $\varphi\sigma(\wb)$ (respectivement 
$\varphi\dot{\sigma}(\nb)$).
 Nous avons ainsi prolong\'e
les actions de $W$ et de $N_\Gb(\Tb)$ sur $\Tb$ aux mono\"{\i}des libres
correspondants.


\bigskip

\subsection{Tores et caract\`eres\label{sub sous}}
Pour \'etudier les repr\'esentations du groupe $\Gb^F$, nous aurons besoin de
quelques r\'esultats de nature combinatoire sur les groupes
$\Tb^{\wb F}$. Leur \'enonc\'e et leur preuve font l'objet des prochains
paragraphes.

\subsubsection{Conjugaison g\'eom\'etrique, s\'eries rationnelles\label{sec conju geo}}
Nous fixons maintenant un entier naturel non nul $d$\index{d@$d$}
multiple de $\delta$ tel que tout \tor $F$-stable
de $\Gb$ soit d\'eploy\'e sur le corps $\FM_{q^{d/\delta}}$. En d'autres termes,
on a $(wF)^d(t)=F^d(t)=t^{q^{d/\delta}}$ pour tout $w \in W$ et pour tout $t \in \Tb$.
Il est \`a noter que $N_{F^d/wF}$ induit une surjection $\Tb^{F^d} \longto \Tb^{wF}$.
De plus, si $t \in \Tb$ est tel que $N_{F^d/wF}(t) \in \Tb^{wF}$,
alors $t \in \Tb^{F^d}$. Pour finir, remarquons que $F^d$ agit
trivialement sur $W$ et que, pour tout $w \in W$, on a
$wF(w)\dots F^{d-1}(w)=1$.

Soit $\zeta=\imath(1/(q^{d/\delta}-1))$\index{zeta@$\zeta$}, o\`u 
$\imath : (\QM/\ZM)_{p'} \longmapright{\sim}
\FM^\times$ est l'\iso choisi dans le \S\ref{sub anneaux} ($\zeta$ est un
g\'en\'erateur du groupe cyclique $\FM_{q^{d/\delta}}^\times$). Si $\wb$ est une
suite d'\'el\'ements de $W$, soit
$$\fonction{N_\wb}{Y(\Tb)}{\Tb^{\wb F}}{\lambda}{N_{F^d/\wb F}(\lambda)(\zeta).}\index{nw@$N_\wb$}$$
Alors, $N_\wb$ est un \mor surjectif de groupes et la suite
\equat\label{norme exacte}
\diagram
0 \rrto && Y(\Tb) \rrto^{\DS{\wb F-1}} && Y(\Tb) \rrto^{\DS{N_\wb}} && \Tb^{\wb F}
\rrto && 0
\enddiagram
\endequat
est exacte \cite[proposition 13.7 (ii)]{DM}.

\medskip

Nous noterons $\nablab(\Tb,W,F)$\index{nabla@$\nablab(\Tb,W,F)$} (ou
$\nablab_\Lambda(\Tb,W,F)$)
l'ensemble des couples $(\wb,\theta)$
o\`u $\wb$ est une suite finie d'\eles de $W$
et $\theta : \Tb^{\wb F} \to \Lambda^\times$ est un \car lin\'eaire
d'ordre inversible dans $\Lambda$. 

L'identification de $W$ au groupe de Weyl de $\tilde{\Gb}$ relatif
\`a $\tilde{\Tb}$ (\cad $W \isom N_{\tilde{\Gb}}(\tilde{\Tb})/\tilde{\Tb}$), permet de d\'efinir
l'application
$$\fonction{\Res}{\nablab(\tilde{\Tb},W,F)}{\nablab(\Tb,W,F)}{(\wb,\tilde{\theta})}{(\wb,
\Res_{\Tb^{\wb F}}^{\tilde{\Tb}^{\wb F}} \tilde{\theta}).}$$
Nous dirons alors que $\Res(\wb,\tilde{\theta})$ est la {\it restriction} de 
$(\wb,\tilde{\theta})$ \`a $\Gb$, ou que $(\wb,\tilde{\theta})$ est une 
{\it extension} de $\Res(\wb,\hat{\theta})$ \`a $\tilde{\Gb}$.

\smallskip
Deux \eles $(\wb,\theta)$ et $(\wb',\theta')$ de $\nablab(\Tb,W,F)$ sont
{\it g\'eom\'etriquement conjugu\'es}  si les \cars lin\'eaires $\theta \circ
N_{F^d/\wb F}$ et $\theta' \circ N_{F^d/\wb' F}$ de
$\Tb^{F^d}$ sont conjugu\'es
sous $W$. D'autre part, les couples $(\wb,\theta)$ et $(\wb',\theta')$ 
appartiennent \`a la m\^eme {\it s\'erie rationnelle} (ce que nous noterons par
$(\wb,\theta) \equivb{W} (\wb',\theta')$) s'il existe des extensions de $(\wb,\theta)$ et
$(\wb',\theta')$ \`a $\tilde{\Gb}$ qui sont g\'eom\'etriquement conjugu\'ees.

\begin{rem}
\label{remarque conju}
Si deux \eles de $\nablab(\Tb_,W,F)$ sont
dans la m\^eme s\'erie rationnelle, alors ils sont g\'eom\'etriquement
conjugu\'es. La r\'eciproque n'est pas vraie en g\'en\'eral, comme le montre
l'exemple du groupe $\Sb\Lb_2(\FM)$ pour $\FM$ de caract\'eristique
diff\'erente de $2$.

La r\'eciproque est n\'eanmoins correcte pour $\tilde{\Gb}$~:
deux \eles de $\nablab(\tilde{\Tb},W,F)$ sont g\'eom\'etriquement conjugu\'es
\ssi ils appartiennent \`a la m\^eme s\'erie rationnelle.
En particulier, la relation d'\'equivalence $\equivb{W}$ sur $\nablab(\Tb,W,F)$
ne d\'epend pas du choix du groupe r\'eductif $\tilde{\Gb}$ \`a centre connexe
choisi (voir le th\'eor\`eme \ref{serie et equivalence} pour une approche ne
faisant pas intervenir $\tilde{\Gb}$).
\end{rem}

\subsubsection{Sous-groupes de $\Tb^{\wb F}$}\label{sousgroupestores}
Nous fixons ici, et ce jusqu'\`a la fin du paragraphe \ref{sub sous}, quatre
suites finies d'\eles de $\bar{S}$,
$$\vb=(s_1^\prime,\dots,s_r^\prime)\ \le\
\xb=(t_1^\prime,\dots,t_r^\prime)\ \le\
\yb=(t_1,\dots,t_r)\ \le\
\wb=(s_1,\dots,s_r).$$

Soit $1 \le i \le r$. Si $s_i\not=1$, nous notons
$\alpha_{\wb,i}$\index{alphawvi@$\alpha_{\wb,i}$}
l'unique racine simple telle que $s_i= s_{\alpha_{\wb,i}}$.
Lorsque $s_i=1$, nous posons $\alpha_{\wb,i}=0$.

Soit $I_{\wb,\vb}=\{1 \le i \le r~|~s_i^\prime\not= s_i\}$\index{iwv@$I_{\wb,\vb}$}.
On pose
$$\alpha_{\wb,\vb,i}=
\begin{cases}
\alpha_{\wb,i}&\text{ si }i\in I_{\wb,\vb} \\
0 & \text{ sinon.}
\end{cases}
$$
\index{alphawvi@$\alpha_{\wb,\vb,i}$}
On a $(s_i-s'_i) Y(\Tb) \incl \ZM \alpha_{\wb,i}^\vee$ pour tout $i$.
Par ailleurs, nous posons
$$\beta_{\wb,\vb,i}^\vee=s_1 \dots s_{i-1}(\alpha_{\wb,\vb,i}^\vee)
\index{betawi@$\beta_{\wb,\vb,i}^\vee$} \text{ et }
\beta_{\wb,i}^\vee=s_1 \dots s_{i-1}(\alpha_{\wb,i}^\vee).
\index{betawi@$\beta_{\wb,i}^\vee$}$$

Nous allons \'etudier ici quelques propri\'et\'es du
sous-r\'eseau $Y_{\wb,\vb}$\index{ywv@$Y_{\wb,\vb}$} de $Y(\Tb)$ engendr\'e
par les $(\beta_{\wb,\vb,i}^\vee)_i$.

\begin{prop}\label{reseau}
Avec les hypoth\`eses et notations ci-dessus, alors

\tete{1} on a $(\wb-\yb)Y(\Tb) \incl Y_{\wb,\vb}$ ;

\tete{2} le groupe $Y_{\wb,\vb}$ est engendr\'e par la famille
$(t_1\dots t_{i-1}(\alpha_{\wb,\vb,i}^\vee))_i$ ;

\tete{3} on a $(\wb F - 1) Y(\Tb) + Y_{\wb,\vb}=(\vb F-1)Y(\Tb) + Y_{\wb,\vb}$~;

\tete{4} les morphismes $N_\wb$ et $N_\yb$ induisent un isomorphisme
$$\Tb^{\wb F}/N_\wb(Y_{\wb,\vb}) \simeq \Tb^{\yb F}/N_\yb(Y_{\wb,\vb}).$$
\end{prop}

\begin{proof}
D\'emontrons tout d'abord (1).
Puisque
$Y_{\wb,\yb} \incl Y_{\wb,\vb}$, il suffit de d\'emontrer que 
$(\wb-\yb)Y(\Tb) \incl Y_{\wb,\yb}$, ce que nous faisons
par r\'ecurrence sur $m=|I_{\wb,\yb}|$.

Si $m=0$, alors $\wb=\yb$ et le r\'esultat en d\'ecoule. Supposons
maintenant que $m=1$. On a alors $I_{\wb,\yb}=\{i\}$.
On pose $\ab=(s_1, \dots ,s_{i -1})$, $\sb = (s_{i})$ et $\bb=(s_{i+1}, \dots ,
s_r)$. Alors,
$$(\wb-\yb)Y(\Tb)=(\ab \sb \bb-\ab  \bb)Y(\Tb)=\ab((\sb-1) Y(\Tb)) \incl
\ab(\ZM\alpha_{\wb,i}^\vee)=\ZM\beta_{\wb,i}^\vee=Y_{\wb,\yb}.$$

Supposons maintenant que $m \ge 2$ et que le r\'esultat est vrai
pour tout couple $(\wb',\yb')$ tel que $\yb' \le \wb'$ et $m_{\wb',\yb'} < m$.
On \'ecrit
$$I_{\wb,\yb}=\{i_1,\dots,i_m\}$$
o\`u $i_1 < \dots < i_m$. On note $\yb'$ l'unique suite
telle que $\yb'\le \wb$ et
$$I_{\wb,\yb'}=\{i_2,\dots, i_m\}.$$
On a
$$(\wb-\yb)Y(\Tb) \incl (\wb-\yb')Y(\Tb) + (\yb'-\yb) Y(\Tb).$$
Or, par hypoth\`ese de r\'ecurrence, $(\wb-\yb')Y(\Tb) \incl Y_{\wb,\yb'}$ et
$(\yb'-\yb)Y(\Tb) \incl Y_{\yb',\yb}$. Le choix que nous avons fait pour $\yb'$
implique que $\beta_{\yb',i_1}^\vee=\beta_{\wb,i_1}^\vee$, donc
$Y_{\wb,\yb'}+Y_{\yb',\yb}=Y_{\wb,\yb}$. Cela termine la d\'emonstration
de (1).

\smallskip
D\'emontrons maintenant (2).
Il suffit de d\'emontrer que, pour tout $i \in I_{\wb,\vb}$,
il existe une famille $(\mu_j)_{j \in I_{\wb,\vb}, j < i}$ d'entiers tels que
$$t_1\dots t_{i-1}(\alpha_{\wb,i}^\vee)=
\beta_{\wb,i}^\vee+\sum_{\substack{j \in I_{\wb,\vb}\\ j<i}}
\mu_j \beta_{\wb,j}^\vee$$
\cad tels que
$$(s_1 \dots s_{i-1} -t_1\dots t_{i-1})(\alpha_{\wb,i}^\vee)=
-\sum_{\substack{j \in I_{\wb,\vb}\\ j<i}}
\mu_j \beta_{\wb,j}^\vee$$
Cela r\'esulte de (1) appliqu\'e dans le cas o\`u $\wb=(s_1,\dots,s_{i-1})$
et $\yb=\vb=(t_1,\dots,t_{i-1})$.

\smallskip
Pour finir la d\'emonstration de la proposition \ref{reseau},
il suffit de remarquer que (3) d\'ecoule imm\'ediatement de (1)
tandis que (4) r\'esulte de (3) et de la suite exacte \ref{norme exacte}.
\end{proof}

\begin{coro}\label{reseau vxyw}
On a $Y_{\wb,\vb}=Y_{\wb,\yb}+Y_{\yb,\vb}$.
\end{coro}

\begin{proof}
On a
$I_{\wb,\vb}=I_{\wb,\yb}\coprod I_{\yb,\vb}$.
Tout d'abord,
$Y_{\yb,\vb}=\sum_{i\in I_{\yb,\vb}} \Zb t_1\cdots
t_{i-1}(\alpha_{\wb,i}^\vee)$ car
$\beta_{\yb,i}^\vee=t_1\cdots t_{i-1}(\alpha_{\wb,i}^\vee)$
pour $i\in  I_{\yb,\vb}$.
La proposition \ref{reseau} (2) appliqu\'ee au cas $\vb=\yb$
montre que
$Y_{\wb,\yb}=\sum_{i\in I_{\wb,\yb}} \Zb t_1\cdots
t_{i-1}(\alpha_{\wb,i}^\vee)$.
Enfin,
$Y_{\wb,\vb}=\sum_{i\in I_{\wb,\vb}} \Zb t_1\cdots
t_{i-1}(\alpha_{\wb,i}^\vee)$
d'apr\`es la proposition \ref{reseau} (2).
\end{proof}

\begin{coro}\label{iso quotient}
Soit $\theta$ un \car lin\'eaire de
$\Tb^{\wb F}$ qui est trivial sur $N_\wb(Y_{\wb,\vb})$.
Notons $\theta_\yb$ le \car lin\'eaire de $\Tb^{\yb F}$ trivial sur
$N_\yb(Y_{\wb,\vb})$ d\'efini via
l'\iso de la proposition $\ref{reseau}$ $(4)$ $($\cad tel que
$\theta_\yb \circ N_\yb=\theta \circ N_\wb)$. Alors, les paires
$(\wb,\theta)$ et $(\yb,\theta_\yb)$ sont dans la m\^eme s\'erie rationnelle.
\end{coro}

\begin{proof}
Pour d\'emontrer le corollaire \ref{iso quotient}, on peut se ramener au cas
o\`u $\Gb=\tilde{\Gb}$~: c'est alors imm\'ediat.
\end{proof}

\medskip

Fixons maintenant un \car lin\'eaire $\theta : \Tb^{\wb F} \to \Lambda^\times$ et
notons $\IC(\wb,\theta)$\index{iwtheta@$\IC(\wb,\theta)$} l'ensemble
des $\wb' \le \wb$ tels que $\theta$ est trivial sur $N_\wb(Y_{\wb,\wb'})$.
L'ensemble $\IC(\wb,\theta)$ admet un plus petit \'el\'ement, not\'e
$\wb_\theta=
(s_{1,\theta},\dots,s_{r,\theta})$\index{wtheta@$\wb_\theta$}, qui est d\'efini par
$$s_{i,\theta}=\left\{\begin{array}{ll}
1& {\mathrm{si}}~\theta(N_\wb(\beta_{\wb,i}^\vee))=1, \\
s_i & {\mathrm{sinon}}.
\end{array}\right.$$
On a donc
$$\IC(\wb,\theta)=\{\xb \in \Sigma(\bar{S})~|~\wb_\theta \le \xb \le \wb\}.$$
D'apr\`es la proposition \ref{reseau} (4), si $\yb \in \IC(\wb,\theta)$,
alors $\theta$ d\'efinit un \car lin\'eaire
$\theta_\yb$ de $\Tb^{\yb F}$ trivial sur $N_\yb(Y_{\wb,\wb_\theta})$. Notons que
$\theta_\wb=\theta$.

\begin{coro}\label{intervalle}
Supposons que $\yb \in \IC(\wb,\theta)$. Alors
$$\IC(\yb,\theta_\yb)=\{\wb' \in \IC(\wb,\theta)~|~\wb_\theta \le \wb' \le \yb\}.$$
En d'autres termes, $\yb_{\theta_\yb}=\wb_\theta$.
\end{coro}

\begin{proof}
Pour tout $\xb\le\yb$, on a
$$\theta\circ N_\wb(Y_{\wb,\xb})=
\theta\circ N_\wb(Y_{\wb,\yb}+Y_{\yb,\xb})$$
d'apr\`es le corollaire \ref{reseau vxyw}. Or, 
$\theta\circ N_\wb(Y_{\wb,\yb})=1$ car $\yb \in \IC(\wb,\theta)$.
Enfin, $\theta\circ N_\wb=\theta_\yb\circ N_\yb$, ce qui montre
que 
$$\theta\circ N_\wb(Y_{\wb,\xb})=\theta_\yb\circ N_\yb(Y_{\yb,\xb}),$$
donc que
$\xb\in \IC(\wb,\theta)$ si et seulement si  $\xb\in \IC(\yb,\theta_\yb)$.
\end{proof}

On dit que $(\wb,\theta),(\wb',\theta')\in\nablab(\Tb,W,F)$ sont
\'el\'ementairement \'equivalents si l'une des deux conditions
suivantes est r\'ealis\'ee~:
\begin{enumerate}
\item
il existe $v\in W$ tel que $f(\sigma(\wb'))=vf(\sigma(\wb))F(v)^{-1}$
et $\theta=\theta'\circ v$ ({\em $W$-conjugaison})
\item
$\wb \in \Sigma(S)$, $\wb'=\wb_\theta$ et $\theta'=\theta_{\wb_\theta}$.
\end{enumerate}
On d\'efinit la relation d'\'equivalence $\sim_W$ sur $\nablab(\Tb,W,F)$ comme
la cl\^oture transitive (et sym\'etrique)
de la relation d'\'equivalence \'el\'ementaire.

\smallskip
Le th\'eor\`eme suivant fournit une nouvelle caract\'erisation de la
conjugaison rationnelle de caract\`eres de tores, qui ne fait pas
intervenir de groupe $\tilde{\Gb}$.

\begin{theo}\label{serie et equivalence}
Deux couples $(\wb,\theta)$ et $(\wb',\theta')$ sont $\sim_W$-\'equivalents
si et seulement si ils sont dans la m\^eme s\'erie rationnelle.
\end{theo}

\begin{proof}
Gr\^ace au corollaire \ref{iso quotient}, il est clair que la relation $\sim_W$ est plus fine 
que la relation $\equiv_W$. Par cons\'equent, le th\'eor\`eme 
\ref{serie et equivalence} est cons\'equence du lemme suivant.
\end{proof}

\begin{lem}\label{serie et minimalite}
Soient $w,w'\in W$ et $((w),\theta)$ et $((w'),\theta')$ deux \'el\'ements de 
$\nablab(\Tb,W,F)$ v\'erifiant les propri\'et\'es suivantes~: 

\tete{1} $((w),\theta)$ et $((w'),\theta')$ sont dans la m\^eme s\'erie
rationnelle ;

\tete{2} $((w),\theta)$ et $((w'),\theta')$ sont minimaux (pour les longueurs
de $(w)$ et $(w')$) 
dans leurs classes respectives pour la relation $\sim_W$.

Alors, $((w),\theta)$ et $((w'),\theta')$ sont conjugu\'es sous $W$.
\end{lem}

\begin{proof} Tout comme dans le corollaire \ref{iso quotient}, on peut se ramener au cas 
o\`u le centre de $\Gb$ est connexe. Posons 
$$\Phi_{w,\theta}^\vee =\{\alpha^\vee \in \Phi^\vee~|~\theta \circ N_w(\alpha^\vee)=1\}\quad{\mathrm{et}}\quad
\Phi_{w,\theta}^{\vee +}=\Phi_{w,\theta}^\vee \cap 
\{\alpha^\vee~|~\alpha\in\Phi^+\}.$$
Le sous-ensemble $\Phi_{w,\theta}^\vee$ de $\Phi^\vee$ est clos et
sym\'etrique. C'est un donc un sous-syst\`eme de $\Phi^\vee$ dont
$\Phi_{w,\theta}^{\vee +}$ est un syst\`eme de coracines positives. On note
$W_{w,\theta}$ le groupe engendr\'e par les $s_\alpha$ avec 
$\alpha^\vee \in \Phi_{w,\theta}^\vee$. C'est aussi le stabilisateur de $\theta \circ N_w$ 
(cela r\'esulte de la connexit\'e du centre de $\Gb$ et de
\cite[th\'eor\`eme 5.13]{delu}).

Il est clair que $w\phi$ stabilise $\Phi_{w,\theta}^\vee$. 
Nous allons montrer que l'hypoth\`ese (2) implique que $w\phi$ stabilise 
$\Phi_{w,\theta}^{\vee +}$. Pour cela, choisissons une d\'ecomposition r\'eduite 
$(s_1,\dots,s_r)$ de $w$ o\`u les $s_i$ sont dans $S$ et notons 
$\alpha_i$ la racine simple associ\'ee \`a $s_i$. Alors les $s_1 \dots s_{i-1}(\alpha_i^\vee)$ 
sont les coracines positives rendues n\'egatives par l'action de $w^{-1}$. Par minimalit\'e 
de $((w),\theta)$, on a $(s_1,\dots,s_r)_\theta=(s_1,\dots,s_r)$. Donc 
$s_1 \dots s_{i-1}(\alpha_i^\vee)$ n'appartient pas \`a $\Phi_{w,\theta}^{\vee}$. 
Cela montre que $\phi^{-1}w^{-1}(\Phi_{w,\theta}^{\vee +}) \incl \Phi^+$.
Finalement, on a bien que $w\phi$ stabilise $\Phi_{w,\theta}^{\vee +}$. 

D'autre part, (1) montre qu'il existe $x \in W$ tel que $\theta \circ N_w = \theta' \circ N_{w'}\circ x$. 
Par cons\'equent $x(\Phi_{w,\theta}^\vee)=\Phi_{w',\theta'}^\vee$, donc 
il existe un \'el\'ement $a \in W_{w',\theta'}$ 
tel que $ax(\Phi_{w,\theta}^{\vee +})=\Phi_{w',\theta'}^{\vee +}$. 
Or, $\Phi_{w,\theta}^{\vee +}$ est $w\phi$-stable, et $\Phi_{w',\theta'}^{\vee +}$ est $w'\phi$-stable. 
Donc $axwF(ax)^{-1}\phi\phi^{-1}w^{\prime -1}=axwF(ax)^{-1}w^{\prime -1}$ stabilise 
$\Phi_{w',\theta'}^{\vee +}$. 

Or, $\theta' \circ N_{w'} \circ axwF(ax)^{-1}w^{\prime -1}=\theta \circ N_w \circ 
wF(ax)^{-1}w^{\prime -1}$. Soit $\lambda \in Y(\Tb)$. Alors
\begin{eqnarray*}
\theta \circ N_{F^d/wF} (\lexp{wF(ax)^{-1}w^{\prime -1}}{\lambda(\zeta)}) &=& 
\theta \circ N_{F^d/wF} (\lexp{wF}{(}\lexp{(ax)^{-1}F^{-1}w^{\prime -1}}{\lambda(\zeta)}))\\
&=& \theta \circ N_{F^d/wF} (\lexp{(ax)^{-1}F^{-1}w^{\prime -1}}{\lambda(\zeta)}) \\
&=&\theta' \circ N_{w'}(\lexp{F^{-1}w^{\prime -1}}{\lambda(\zeta)}) \\
&=&\theta' \circ N_{w'}({\lambda(\zeta)}). 
\end{eqnarray*}
Ici, $F^{-1}$ est vu comme un automorphisme du groupe $\Tb$.
Ceci montre que $axwF(ax)^{-1}w^{\prime -1}$ est dans $W_{w',\theta'}$. Puisqu'il stabilise 
$\Phi_{w',\theta'}^{\vee +}$, il est trivial. Par cons\'equent, 
on a $w'=axwF(ax)^{-1}$. En outre, $\theta'=\theta \circ ax$ car la norme est surjective.
\end{proof}

\begin{rem}
La preuve montre aussi qu'un couple minimal correspond \`a un couple
maximalement d\'eploy\'e au sens de Deligne et Lusztig
\cite[d\'efinition 5.25]{delu}. Lorsque le centre de $\Gb$ est connexe,
il est d\'emontr\'e dans \cite[proposition 5.26]{delu} que deux couples
$((w),\theta)$ et $((w'),\theta')$ appartenant \`a une m\^eme s\'erie
rationnelle et maximalement d\'eploy\'es sont conjugu\'es sous $W$, cas
particulier du lemme \ref{serie et minimalite}.
\end{rem}

\subsubsection{Une application\label{application}}
L'objet de ce paragraphe est la d\'efinition et l'\'etude d'un groupe
diagonalisable qui interviendra dans la construction de recollements
de vari\'et\'es de Deligne-Lusztig au \S \ref{sub recollement}.

%

Fixons tout d'abord quelques notations. D\'efinissons 
$$
\Sb_{\wb,\vb}=\bigl\{(a_1,\dots,a_r) \in \Tb^r~|~
a_i^{-1}~\lexp{s_i}{a_{i+1}} \in \Im \alpha_{\wb,\vb,i}^\vee\ \mathrm{ pour }\ 
1 \le i \le r-1\
\mathrm{et}\ a_r^{-1}~\lexp{s_r F}{a_1} \in
 \Im \alpha_{\wb,\vb,r}^\vee ~\bigr\}\index{swv@$\Sb{\wb,\vb}$}
$$
et
$$\fonction{\mu_\yb}{\Tb}{\Tb^r}{a}{(a,\lexp{t_1}{a},
\dots,\lexp{t_{r-1}\dots t_1}{a}).}\index{muy@$\mu_\yb$}$$

\begin{prop}\label{Sw}
Avec les notations ci-dessus, on a~:

\tete{1} $\mu_\yb(\Tb^{\yb F}) \incl \Sb_{\wb,\vb}$.

\tete{2} $\Sb_{\wb,\vb}=\mu_\yb(\Tb^{\yb F})\cdot \Sb_{\wb,\vb}^\circ$.

\tete{3} $\mu_\yb(N_\yb(Y_{\wb,\vb}))\incl \Sb_{\wb,\vb}^\circ$
et le diagramme suivant est commutatif
$$\xymatrix{
 & \Tb^{\yb F}/N_\yb(Y_{\wb,\vb}) \ar[dr]^{\mu_\yb} \\
Y(\Tb)\ar[ur]^{N_\yb} \ar[dr]_{N_\wb} & & \Sb_{\wb,\vb}/\Sb_{\wb,\vb}^\circ \\
 & \Tb^{\wb F}/N_\wb(Y_{\wb,\vb}) \ar[ur]_{\mu_\wb} \\
}$$ 

\tete{4} Si les coracines de $\Gb$ sont injectives, alors
$\mu_\yb(\Tb^{\yb F}) \cap \Sb_{\wb,\vb}^\circ = \mu_\yb(N_\yb(Y_{\wb,\vb}))$.
\end{prop}

\begin{proof}
Posons $\Gb'=\Gb^r$, $\Tb'=\Tb^r$ et $W'=W^r$. Dans cette preuve, 
les suites $\vb$, $\yb$ et $\wb$ seront vues comme des \'el\'ements 
de $W'$. Posons d'autre part
$$\Tb_{\wb,\vb}^\prime=\Im \alpha_{\wb,\vb,1}^\vee \times \dots \times 
\Im \alpha_{\wb,\vb,r}^\vee.$$
On a alors
$$\lexp{\wb}{t}~\lexp{\yb}{t}^{-1} \in \Tb_{\wb,\vb}^\prime\leqno{(*)}$$ 
pour tout $t \in \Tb'$. Pour finir, notons
$$\fonction{F'}{\Gb'}{\Gb'}{(g_1,\dots,g_r)}{(g_2,\dots,g_r,F(g_1)).}$$
Alors $F^{\prime r\delta}$ est un endomorphisme de Frobenius sur 
$\Gb'$ qui munit ce dernier d'une structure rationnelle 
sur $\fq$. Avec ces notations, le groupe diagonalisable $\Sb_{\wb,\vb}$
peut \^etre d\'efini ainsi~:
$$\Sb_{\wb,\vb}=\{t \in \Tb'~|~t^{-1}~\lexp{\wb F'}{t} \in \Tb_{\wb,\vb}^\prime\}.$$
D'apr\`es $(*)$, cette caract\'erisation est \'equivalente 
\`a la suivante~:
$$\Sb_{\wb,\vb}=\{t \in \Tb'~|~t^{-1}~\lexp{\yb F'}{t} \in \Tb_{\wb,\vb}^\prime\}.$$
Compte tenu de cette derni\`ere d\'efinition, on a 
$\Tb^{\prime \yb F'} \incl \Sb_{\wb,\vb}$ et 
$\Sb_{\wb,\vb}=\Tb^{\prime \yb F'}\cdot\Sb_{\wb,\vb}^\circ$ (voir remarque 
\ref{hf}). Or, un calcul facile montre que $\Tb^{\prime \yb F'}=\mu_\yb(\Tb^{\yb F})$, 
ce qui termine la preuve de (1) et (2). 

D\'emontrons maintenant (3). 
On a, pour tout $w' \in W'$, $(w'F')^{rd}=F^{\prime rd}$, et 
$F^{\prime rd}$ est un endomorphisme de Frobenius d\'eploy\'e sur $\Tb'$. 
On d\'efinit alors
$$\fonction{N_{w'}^\prime}{Y(\Tb')}{\Tb^{\prime w'F'}}{\lambda}{
N_{F^{\prime rd}/w'F'}(\lambda(\zeta)).}$$

Remarquons tout d'abord que $N_{F^{\prime rd}/\yb F'}(\Tb_{\wb,\vb}^\prime)$ 
est contenu dans $\Sb_{\wb,\vb}$ car
$\yb F'\circ N_{F^{\prime rd}/\yb F'}=F^{\prime rd}$ est
d\'eploy\'e donc stabilise $\Tb_{\wb,\vb}^\prime$.
Puisque $N_{F^{\prime rd}/\yb F'}(\Tb_{\wb,\vb}^\prime)$
est connexe, il est en fait contenu dans $\Sb_{\wb,\vb}^\circ$. D'autre part, 
le noyau de $N_{F^{\prime rd}/\yb F'}$ \'etant fini (il est contenu dans 
$\Tb^{\prime F^{\prime rd}}$), on a 
$$\dim \Tb_{\wb,\vb}^\prime=\dim N_{F^{\prime rd}/\yb F'}(\Tb_{\wb,\vb}^\prime)
=\dim \Sb_{\wb,\vb}.$$
Par cons\'equent, 
$$\Sb_{\wb,\vb}^\circ=N_{F^{\prime rd}/\yb F'}(\Tb_{\wb,\vb}^\prime).\leqno{(a)}$$
Ce r\'esultat \'etant vrai pour tout $\yb$ compris entre $\vb$ et $\wb$, 
on a aussi 
$$\Sb_{\wb,\vb}^\circ=N_{F^{\prime rd}/\wb F'}(\Tb_{\wb,\vb}^\prime).\leqno{(a')}$$

Pour simplifier les notations, posons $Y_{\wb,\vb}^\prime=Y(\Tb_{\wb,\vb}^\prime)$. 
C'est un sous-r\'eseau facteur direct de $Y(\Tb')$. On a alors, d'apr\`es 
(a), 
$$Y(\Sb_{\wb,\vb}^\circ)=\Bigl(\Qb \otimes_\Zb 
N_{F^{\prime rd}/\yb F'}(Y_{\wb,\vb}^\prime) \Bigr) \cap Y(\Tb').\leqno{(b)}$$
De m\^eme, d'apr\`es (a'), on a 
$$Y(\Sb_{\wb,\vb}^\circ)=\Bigl(\Qb \otimes_\Zb 
N_{F^{\prime rd}/\wb F'}(Y_{\wb,\vb}^\prime) \Bigr) \cap Y(\Tb').\leqno{(b')}$$
L'\'egalit\'e (b) montre en particulier que 
$$N_\yb^\prime(Y_{\wb,\vb}^\prime) \incl \Sb_{\wb,\vb}^\circ.\leqno{(c)}$$
De m\^eme (b') montre que 
$$N_\wb^\prime(Y_{\wb,\vb}^\prime) \incl \Sb_{\wb,\vb}^\circ.\leqno{(c')}$$
D'autre part, si $\lambda' \in Y(\Tb')$, alors il existe 
$\mu' \in \Qb \otimes_\Zb Y(\Tb')$ tel que $\lambda'=(\wb F' -1)(\mu')$. 
Par suite, 
\begin{eqnarray*}
&&N_{F^{\prime rd}/\wb F'}(\lambda')-N_{F^{\prime rd}/\yb F'}(\lambda')\\
&=&N_{F^{\prime rd}/\wb F'} \circ (\wb F' -1) (\mu') 
- N_{F^{\prime rd}/\yb F'} \circ (\yb F' -1) (\mu')
- N_{F^{\prime rd}/\yb F'} \circ (\yb-\wb) \circ F'(\mu')\\
&=&(q^{d/\delta}-1)\mu'-(q^{d/\delta}-1)\mu'-
N_{F^{\prime rd}/\yb F'} \circ (\yb-\wb) \circ F'(\mu')\\
&=&-N_{F^{\prime rd}/\yb F'} \circ (\yb-\wb) \circ F'(\mu').
\end{eqnarray*}
Or, $(\yb - \wb)\circ F'(\mu') \in \Qb \otimes_\Zb Y_{\wb,\vb}^\prime$, donc 
$$N_{F^{\prime rd}/\wb F'}(\lambda')-N_{F^{\prime rd}/\yb F'}(\lambda')
\in \Bigl(\Qb \otimes_\Zb N_{F^{\prime rd}/\yb F'}(Y_{\wb,\vb}^\prime)\Bigr)
\cap Y(\Tb')=Y(\Sb_{\wb,\vb}^\circ).$$
Cela montre que le diagramme 
$$\xymatrix{
 & \Tb^{\prime \yb F'}/N_\yb^\prime(Y_{\wb,\vb}^\prime) \ar[dr] \\
Y(\Tb')\ar[ur]^{N_\yb^\prime} \ar[dr]_{N_\wb^\prime} & & 
\Sb_{\wb,\vb}/\Sb_{\wb,\vb}^\circ \\
 & \Tb^{\prime \wb F'}/N_\wb^\prime(Y_{\wb,\vb}^\prime) \ar[ur] \\
}$$
est commutatif.

Pour d\'eduire (3), il faut encore comparer $\mu_\yb\circ N_\yb$ et 
$N_\yb^\prime$.
Soit $p_1$ la premi\`ere projection $\Tb'\to \Tb$.
On a
$$p_1\circ N_\yb^\prime(0,\ldots,0,\lambda,0,\ldots,0)=
N_{F^d/\yb F}(t_1\cdots t_{i-1}(\lambda)(\zeta))=
N_\yb(t_1\cdots t_{i-1}(\lambda)),$$
o\`u le terme $\lambda$ a \'et\'e plac\'e en $i$-\`eme 
position.
Puisque $p_1\circ\mu_\yb$ est l'identit\'e, on a
$$p_1\circ\mu_\yb \circ N_\yb(t_1\dots t_{i-1}(\lambda))=
p_1\circ N_\yb^\prime(0,\dots,0,\lambda,0,\dots,0).$$
Or, $N_\yb^\prime(0,\dots,0,\lambda,0,\dots,0)\in \Tb^{\prime \yb F'}=
\mu_\yb(\Tb^{\yb F})$, donc
$$\mu_\yb \circ N_\yb(t_1\dots t_{i-1}(\lambda))=
N_\yb^\prime(0,\dots,0,\lambda,0,\dots,0).\leqno{(d)}$$

L'inclusion
$\mu_\yb(N_\yb(Y_{\wb,\vb}))\incl \Sb_{\wb,\vb}^\circ$
d\'ecoule de (c), (d) et de la proposition \ref{reseau}, (2).
La commutativit\'e du diagramme dans (3) r\'esulte de la commutativit\'e du
diagramme pr\'ec\'edent et de (d).

Passons enfin \`a (4). Compte tenu de (a), (d) et de la proposition
\ref{reseau}, (2), il suffit de montrer que
$$N_{F^{\prime rd}/\yb F'}(\Tb_{\wb,\vb}^\prime) \cap \Tb^{\prime \yb F'} = 
N_\yb^\prime(Y_{\wb,\vb}^\prime).$$
Or, si $t \in \Tb'$, $N_{F^{\prime rd}/\yb F'}(t) \in \Tb^{\prime \yb F'}$ 
si et seulement si $t \in \Tb^{\prime F^{\prime rd}}$. 
Par cons\'equent, 
$$N_{F^{\prime rd}/\yb F'}(\Tb_{\wb,\vb}^\prime) \cap \Tb^{\prime \yb F'} = 
N_{F^{\prime rd}/\yb F'}(\Tb_{\wb,\vb}^{\prime F^{\prime rd}}).$$
Il suffit donc de prouver que 
$$\Tb_{\wb,\vb}^{\prime F^{\prime rd}}=\{\lambda(\zeta)~|~\lambda\in
Y_{\wb,\vb}^\prime\}.$$
Mais 
$$\Tb_{\wb,\vb}^{\prime F^{\prime rd}}=(\Im \alpha_{\wb,\vb,1}^\vee)^{F^d}
\times \dots \times (\Im \alpha_{\wb,\vb,r}^\vee)^{F^d}.$$
Puisque les coracines de $\Gb$ sont injectives et que l'endomorphisme 
$F^d$ de $\Tb$ est un endomorphisme de Frobenius d\'eploy\'e sur
$\Fb_{q^{d/\delta}}$, 
on a 
$$(\Im \alpha_{\wb,\vb,i}^\vee)^{F^d}=\alpha_{\wb,\vb,i}^\vee(\Fb_{q^{d/\delta}}^\times)$$
pour $1 \le i \le r$. 
Cela montre que
$$\Tb_{\wb,\vb}^{\prime F^{\prime rd}}=
\alpha_{\wb,\vb,1}^\vee(\Fb_{q^{d/\delta}}^\times)
\times \dots \times \alpha_{\wb,\vb,r}^\vee(\Fb_{q^{d/\delta}}^\times).$$
La preuve du (4) est maintenant compl\`ete.
\end{proof}

\begin{rem}
\label{SW}
La proposition pr\'ec\'edente donne une
autre construction de l'isomorphisme de la proposition \ref{reseau} (4),
montrant au passage qu'elle est ind\'ependante des choix des
isomorphismes de \S \ref{sub anneaux}.
\end{rem}

\section{Vari\'et\'es de Deligne-Lusztig g\'en\'eralis\'ees}~

\subsection{Comparaison de mod\`eles}
Nous commen\c{c}ons cette partie par un r\'esultat g\'eom\'e\-trique 
concernant une classe de vari\'et\'es qui contient les vari\'et\'es de 
Deligne-Lusztig. Les hypoth\`eses ont \'et\'e choisies pour les
applications des parties \`a venir.

Fixons un \ele $n$ de $\Gb$, un sous-groupe parabolique
$\Pb$ de $\Gb$ de radical unipotent $\Vb$ et soit $\Lb$ un compl\'ement
de Levi $nF$-stable ($(n F)^\delta$ est alors un 
endomorphisme de Frobenius relatif \`a une $\fq$-structure sur $\Lb$).

Fixons aussi un sous-groupe ferm\'e connexe $nF$-stable 
$\Hb$ de $\Lb$, un sous-groupe distingu\'e $\Hb'$ de $\Hb$ (non
n\'ecessairement $nF$-stable) et posons 
$$\Kb=\{h \in \Hb~|~h^{-1}~\lexp{n F}{h} \in \Hb'\}.$$
Puisque $\Hb'$ est distingu\'e dans $\Hb$, $\Kb$ est un 
sous-groupe ferm\'e de $\Hb$. D'autre part, $\Hb^{n F} \incl \Kb$.

Posons maintenant
$$\Yb=\{g\Vb \in \Gb/\Vb~|~g^{-1}F(g) \in \Hb'\cdot(\Vb n
~\lexp{F}{\Vb})\}$$
$$\Xb=\{g\Hb\cdot\Vb \in \Gb/(\Hb\cdot\Vb)~|~g^{-1}F(g) \in
\Hb\cdot(\Vb n ~\lexp{F}{\Vb})\}.
\leqno{\mathrm{et}}$$
Alors, le groupe $\Gb^F$ agit par translation \`a gauche sur $\Yb$ et
$\Xb$, 
le groupe $\Kb$ agit librement par translation \`a droite sur $\Yb$ et
le 
morphisme 
$$\fonction{\pi}{\Yb}{\Xb}{g\Vb}{g\Hb\cdot\Vb}$$
est $\Gb^F$-\'equivariant. De plus, ses fibres sont des $\Kb$-orbites. 

\bigskip

\begin{prop}\label{existence quotient}
Le morphisme $\pi$ induit un isomorphisme $\Yb/\Kb \isom \Xb$. 
\end{prop}

\bigskip

\begin{proof} Tout d'abord, 
$$\Yb/\Kb=\{g\Kb\cdot \Vb\in \Gb/(\Kb\cdot\Vb)~|~g^{-1}F(g) \in \Hb'\cdot
(\Vb n ~\lexp{F}{\Vb})\}.$$
Notons 
$$\fonction{\pi'}{\Yb/\Kb}{\Xb}{g\Kb\cdot\Vb}{g\Hb\cdot\Vb.}$$
Il s'agit de montrer que $\pi'$ est un isomorphisme de vari\'et\'es. 
Nous allons pour cela construire un inverse explicite.

L'application
$$\fonctio{\Hb \times \Vb n ~\lexp{F}{\Vb}}{\Hb\cdot(\Vb n 
~\lexp{F}{\Vb})}{(h,x)}{hx}$$
est un isomorphisme de vari\'et\'es. Nous noterons 
$\lambda : \Hb\cdot(\Vb n ~\lexp{F}{\Vb})\to \Hb$ le compos\'e de 
l'inverse de cet isomorphisme avec la premi\`ere projection. 

D'autre part, l'application de Lang $\Hb \to \Hb$, 
$g \mapsto g~\lexp{n F}{g}^{-1}$ est un morphisme \'etale, donc elle
induit 
un isomorphisme $\gamma_0:\Hb/\Hb^{n F} \isom \Hb$. Nous noterons 
$\gamma$ le morphisme compos\'e
$$\diagram
\Hb \rto^{\gamma_0^{-1}\hspace{0.3cm}} & \Hb/\Hb^{n F} \rto
&\Gb/(\Kb\cdot\Vb),
\enddiagram$$
o\`u la deuxi\`eme fl\`eche est l'application canonique. 

Pour finir, posons 
$$\Xb'=\{g \in \Gb~|~g^{-1}F(g) \in \Hb\cdot(\Vb n
~\lexp{F}{\Vb})\},$$
et notons 
$$\fonction{\mu}{\Xb'}{\Yb/\Kb}{g}{g \gamma(\lambda(g^{-1}F(g)))
\Kb\cdot\Vb.}$$
Montrons que $\mu$ est bien d\'efini. Prenons
$g\in\Xb'$ et posons $h=\lambda(g^{-1}F(g))$. On a

$$(g\gamma_0^{-1}(h))^{-1}\lexp{F}{(g\gamma_0^{-1}(h))}\in
\gamma_0^{-1}(h)^{-1}h\lexp{nF}{\gamma_0^{-1}(h)} \Vb n
\lexp{F}{\Vb}=\Vb n \lexp{F}{\Vb}
$$
ce qui montre que $\mu$ est bien d\'efini --- c'est de plus
un morphisme de vari\'et\'es. Pour $g \in \Xb'$, $h\in\Hb$ et
$v\in\Vb$, on a
$$\gamma_0^{-1}(\lambda((ghv)^{-1}F(ghv)))=
h^{-1}\gamma_0^{-1}(\lambda(g^{-1}F(g))).$$
On en d\'eduit que $\mu(ghv)=\mu(g)$.
Par passage au quotient, $\mu$ induit un morphisme de vari\'et\'es 
$\mu' : \Xb'/\Hb\cdot\Vb=\Xb \to \Yb/\Kb$. C'est l'isomorphisme
inverse de $\pi'$.
Cela compl\`ete la preuve de la proposition.
\end{proof}

\subsection{D\'efinition} Pour $w \in W$ (\resp $n \in N_\Gb(\Tb)$) et
pour $g$ et $h$
deux \eles de $\Gb$, nous
\'ecrirons $g\Bb \longtrait{w} h\Bb$ (\resp $g\Ub \longtrait{n}
h\Ub$) pour dire que $g^{-1}h \in \Bb w\Bb$ (\resp $g^{-1}h \in \Ub n \Ub$).
Pour $\wb=(w_1,\dots,w_r) \in \Sigma(W)$ et
pour $\nb=(n_1,\dots,n_r) \in \Sigma(N_\Gb(\Tb))$, nous posons
$$\Xb(\wb)=\bigl\{(g_1\Bb,\dots,g_r\Bb) \in (\Gb/\Bb)^r~|~
g_1\Bb \longtrait{w_1} g_2\Bb \longtrait{w_2} \cdots \longtrait{w_{r-1}} g_r \Bb
\longtrait{w_r} F(g_1) \Bb\}\index{xw@$\Xb(\wb)$}$$
et
$$\Yb(\nb)=\bigl\{(g_1\Ub,\dots,g_r\Ub) \in (\Gb/\Ub)^r~|~
g_1\Ub \longtrait{n_1} g_2\Ub \longtrait{n_2} \cdots \longtrait{n_{r-1}} g_r \Ub
\longtrait{n_r} F(g_1) \Ub\}.\index{yn@$\Yb(\nb)$}$$
Supposons $\wb=\bar{\nb}$.
Alors, le groupe $\Tb^{\wb F}$ agit \`a droite
sur $\Yb(\nb)$
de la mani\`ere suivante~: si $(g_1\Ub,\dots,g_r\Ub) \in
\Yb(\nb)$ et si $t \in \Tb^{\wb F}$, alors
$$(g_1\Ub,\dots,g_r\Ub).t=(g_1t\Ub,g_2 \lexp{n_1^{-1}}{t} \Ub,
\dots,g_r\lexp{(n_1\dots n_{r-1})^{-1}}{t}\Ub).$$
Notons que cette action est libre.
D'autre part, $\Gb^F$ agit diagonalement sur $\Yb(\nb)$ par multiplication
\`a gauche et cette action
commute avec celle de $\Tb^{\wb F}$. Cela fait de $\Yb(\nb)$ une
$\Gb^F$-vari\'et\'e-$\Tb^{\wb F}$. 

\bigskip

Nous noterons $\Yb^\Gb(\nb)$ ou $\Yb^{\Gb,F}(\nb)$ 
(respectivement $\Xb^\Gb(\wb)$ ou $\Xb^{\Gb,F}(\wb)$) la vari\'et\'e 
$\Yb(\nb)$ (respectivement $\Xb(\wb)$) lorsque nous aurons besoin de pr\'eciser 
le groupe ambiant, voire l'isog\'enie consid\'er\'ee. 

\begin{rem}
\label{B=W}
Posons
$\Gb'=\Gb^r$, $\Bb'=\Bb^r$, $\Ub'=\Ub^r$, $\Tb'=\Tb^r$,
$W'=W^r$ et
$$\fonction{F'}{\Gb'}{\Gb'}{(g_1,\dots,g_r)}{(g_2,\dots,g_r,F(g_1)).}
\leqno{\mathrm{et}}$$
Alors, $F'$ est une isog\'enie sur $\Gb'$, $F^{\prime r\delta}$ d\'efinit
une $\fq$-structure sur le groupe $\Gb'$ et $\nb$ s'identifie
\`a un \ele du groupe $N_{\Gb'}(\Tb')$. On v\'erifie
que les applications
\begin{align*}
\Tb^{\nb F}&\isom \Tb^{\prime \nb F'}\\
t&\mapsto (t,\lexp{n_1^{-1}}{t},\dots, \lexp{(n_1\dots n_{r-1})^{-1}}{t})\\
\Gb^F&\isom \Gb^{\prime F'}\\
g&\mapsto (g,\dots,g)
\end{align*}
sont des \isos de groupes et que l'application
$$\fonctio{\Yb^{\Gb,F}(\nb)}{\Yb^{\Gb',F'}(\nb)}{(g_1\Ub,
\dots,g_r\Ub)}{(g_1,\dots,g_r)\Ub'}$$
est un \iso de $\Gb^F$-vari\'et\'es-$\Tb^{\nb F}$, o\`u $\Yb^{\Gb',F'}(\nb)$
est vue comme une $\Gb^F$-vari\'et\'e-$\Tb^{\nb F}$ via les \isos
pr\'ec\'edents.

Une interpr\'etation analogue s'en d\'eduit pour les vari\'et\'es
$\Xb^{\Gb,F}(\wb)$. Cela montre que les vari\'et\'es de Deligne-Lusztig
g\'en\'eralis\'ees ne sont en fait que des vari\'et\'es de Deligne-Lusztig
ordinaires pour un groupe diff\'erent.
\end{rem}

Le morphisme canonique
$$\fonction{\pi_\nb}{\Yb(\nb)}{\Xb(\wb)}{
(g_1\Ub,\dots,g_r\Ub)}{(g_1\Bb,\dots,g_r\Bb)}\index{pin@$\pi_\nb$}$$
induit un isomorphisme de $\Gb^F$-vari\'et\'es 
\equat\label{YversX}
\Yb(\nb)/\Tb^{\wb F} \isom \Xb(\wb).
\endequat
Pour voir cela, 
il faut se ramener au cas o\`u $r=1$ en utilisant la remarque \ref{B=W}, 
puis appliquer la proposition \ref{existence quotient} dans le cas 
o\`u $\Pb=\Bb$, $\Lb=\Hb=\Tb$, $n=n_1$ et $\Hb'=\{1\}$. 

\bigskip

Puisque $\Yb(\nb)$ est une $\Gb^F$-vari\'et\'e-$\Tb^{\wb F}$,
elle induit deux foncteurs entre cat\'egories d\'eriv\'ees qui seront not\'es
$$\fonction{\RC_\nb}{D^b(\Lambda\Tb^{\wb F})}{D^b(\Lambda\Gb^F)}{C}{
\RRR\Gamma_c(\Yb(\nb)) \otimes_{\Lambda\Tb^{\wb F}}^\LLL C}\index{rn@$\RC_\nb$}$$
$$\fonction{\SC_\nb}{D^b(\Lambda\Tb^{\wb F})}{D^b(\Lambda\Gb^F)}{C}{
\RRR\Gamma(\Yb(\nb)) \otimes_{\Lambda\Tb^{\wb F}}^\LLL C.}
\leqno{\mathrm{et}}\index{sn@$\SC_\nb$}$$
Ce sont les foncteurs
$\RC_{\Tb^{\wb F}}^{\Gb^F}(\Yb(\nb))$ et $\SC_{\Tb^{\wb F}}^{\Gb^F}(\Yb(\nb))$
de \S\ref{sub foncteur}.

Le foncteur $\RC_\nb$
induit une application lin\'eaire entre les groupes de Grothendieck des
cat\'egories
$D^b(\Lambda\Tb^{\wb F})$ et $D^b(\Lambda\Gb^F)$.
Cette application sera not\'ee
$$\RRR_\nb : \K0(\Lambda\Tb^{\wb F}) \longto
\K0(\Lambda\Gb^F).\index{rn@$\RRR_\nb$}$$

\bigskip

Les vari\'et\'es $\Xb(\wb)$ et $\Yb(\nb)$ ne d\'ependent, \`a isomorphisme
pr\`es, que des images de $\wb$ et $\nb$ dans $B$ \cite{dimirou}~:

\begin{prop}\label{remark tresse}
Soient $\vb$ et $\wb$ deux suites d'\eles de $W$
telles que $\sigma(\vb)=\sigma(\wb)$. Alors,
les $\Gb^F$-vari\'et\'es $\Xb(\wb)$ et $\Xb(\vb)$ sont
isomorphes. En particulier, une
d\'ecomposition r\'eduite de $\sigma(\wb)$ dans $B$ fournit une suite 
$\vb \in \Sigma(S)$
telle que $\Xb(\wb)$ est isomorphe \`a $\Xb(\vb)$.

De m\^eme, soient $\mb$ et $\nb$  deux suites d'\eles de $N_\Gb(\Tb)$ telles
que $\dot{\sigma}(\mb)\Tb=\dot{\sigma}(\nb)\Tb$. Alors,
on a un isomorphisme de $\Gb^F$-vari\'et\'es-$\Tb^{\nb F}$ entre
$\Yb(\mb)$ et $\Yb(\nb)$.  Par cons\'equent, il existe
$\wb\in\Sigma(S)$ telle que
$\Yb(\nb)$ est isomorphe \`a $\Yb(\dot{\wb})$.
\end{prop}

\begin{proof}
Soit $i$ tel que $l(\bar{n}_i\bar{n}_{i+1})=l(\bar{n}_i)+l(\bar{n}_{i+1})$. Alors, on a
un isomorphisme de $\Gb^F$-vari\'et\'es-$\Tb^{\wb F}$

\refstepcounter{theo}
\begin{align*}
\label{concatenationY}
\Yb(\nb)&\isom
 \Yb(n_1,\ldots,n_{i-1},n_in_{i+1},n_{i+2},\ldots,n_r),
    \tag{$\boldsymbol{\arabic{section}.\arabic{theo}}$} \\
(g_1\Ub,\ldots,g_r\Ub)&\mapsto
 (g_1\Ub,\ldots,g_{i-1}\Ub,g_{i+1}\Ub,\ldots,g_r\Ub).
\end{align*}
Ceci montre que $\Yb(\nb)$ et
$\Yb(\mb)$ sont isomorphes lorsque $\dot{\sigma}(\mb)=\dot{\sigma}(\nb)$.

\smallskip
Supposons maintenant $n_r\in \Tb$. Soit $t\in \Tb$ tel que
$t^{-1}~\lexp{\nb F}{t}=\lexp{\nb}{n_r}$. Alors, on a un
isomorphisme de $\Gb^F$-vari\'et\'es-$\Tb^{\wb F}$
\begin{align*}
\Yb(\nb)&\isom \Yb(n_1,\ldots,n_{r-1}),\\
(g_1\Ub,\ldots,g_r\Ub)&\mapsto (g_1t\Ub,g_2\lexp{n_1^{-1}}{t}\Ub,\ldots,
g_{r-1}\lexp{n_{r-2}^{-1}\cdots n_1^{-1}}{t}\Ub).
\end{align*}
On en d\'eduit donc que $\Yb(\nb)$ et 
$\Yb(\mb)$ sont isomorphes lorsque 
$\dot{\sigma}(\mb)\Tb=\dot{\sigma}(\nb)\Tb$. Puisque 
$\Xb(\wb)\simeq \Yb(\nb)/\Tb^{\wb F}$,
lorsque $\wb$ est l'image de $\nb$, on en d\'eduit aussi la premi\`ere
partie de la proposition.
\end{proof}

\begin{rem}
D'apr\`es la proposition \ref{remark tresse},
l'\'etude des vari\'et\'es $\Xb(\wb)$
et $\Yb(\dot{\wb})$ o\`u $\wb\in\Sigma(W)$ se ram\`ene au cas o\`u 
$\wb\in\Sigma(S)$.
\end{rem}

\begin{prop}\label{goodguy}
La vari\'et\'e $\Yb(\nb)$ est quasi-affine, lisse et purement de dimension
$l(\nb)$. Les stabilisateurs des points de $\Yb(\nb)$ dans $\Gb^F\times
\Tb^{\nb F}$ sont des $p$-groupes.
\end{prop}

\begin{proof}
La remarque \ref{B=W} ram\`ene la preuve de la proposition au cas
de vari\'et\'es de Deligne-Lusztig ordinaires. Alors,
la quasi-affinit\'e est due \`a Haastert \cite[th\'eor\`eme 2.3]{Haastert}
et les autres propri\'et\'es sont faciles \cite[\S 1.3]{delu}.
\end{proof}

\begin{rem}
Il est connu que les vari\'et\'es $\Yb(\nb)$ sont affines pour $q^{1/\delta}$
assez grand \cite[th\'eor\`eme 9.7]{delu}, mais ce r\'esultat
demeure inconnu pour $q^{1/\delta}$ quelconque.
\end{rem}

\begin{rem}
\label{UNU}
Pour $n \in N_\Gb(\Tb)$, posons
$$\Yb^*(n)=\{g (\Ub \cap \lexp{n}{\Ub}) \in \Gb/(\Ub \cap \lexp{n}{\Ub})~|~g^{-1}F(g)
\in n\Ub\}.$$
Alors, l'application canonique $\Yb^*(n) \to \Yb(n)$, $g(\Ub \cap \lexp{n}{\Ub})
\mapsto g\Ub$ est un isomorphisme de $\Gb^F$-vari\'et\'es-$\Tb^{nF}$.
On trouve dans la litt\'erature indiff\'eremment l'un ou l'autre de
ces mod\`eles.
\end{rem}

\begin{rem}
Puisque la vari\'et\'e $\Xb(\wb)$ est lisse,
la dualit\'e de Poincar\'e fournit
un isomorphisme fonctoriel de $\Lambda\Gb^F$-modules
\equat\label{iso dual}
\left((\RC_\nb M)^*\right)^\opp \simeq \SC_\nb \left((M^*)^\opp\right)[-2 l(\wb)]
\endequat
pour tout tout $\Lambda\Tb^{\wb F}$-module $M$, libre sur $\Lambda$.
\end{rem}

\subsection{Liens entre les groupes $\Gb$, $\tilde{\Gb}$ et $\hat{\Gb}$\label{lien lien}}
En identifiant
les groupes de Weyl de $\Gb$, $\tilde{\Gb}$ et $\hat{\Gb}$ relatifs \`a $\Tb$, 
$\tilde{\Tb}$ et $\hat{\Tb}$ respectivement, nous noterons, lorsque $\wb$ est une 
suite d'\eles de $W$, $\tilde{\Xb}(\wb)$ et $\hat{\Xb}(\wb)$ les vari\'et\'es 
d\'efinies gr\^ace aux groupes $\tilde{\Gb}$ et $\hat{\Gb}$ respectivement. De 
m\^eme, si $\nb$ est une suite d'\eles de $N_{\tilde{\Gb}}(\tilde{\Tb})$ ou de 
$N_{\hat{\Gb}}(\hat{\Tb})$, on d\'efinit des vari\'et\'es
$\tilde{\Yb}(\nb)$\index{ynt@$\tilde{\Yb}(\nb)$}
et $\hat{\Yb}(\nb)$\index{ynh@$\hat{\Yb}(\nb)$} et nous noterons
$\tilde{\RC}_\nb$\index{rnt@$\tilde{\RC}_\nb$}
et $\hat{\RC}_\nb$\index{rnh@$\hat{\RC}_\nb$}
les foncteurs associ\'es aux vari\'et\'es $\tilde{\Yb}(\nb)$
et $\hat{\Yb}(\nb)$ respectivement.

Dans ce paragraphe, nous allons rappeler les liens
entre ces vari\'et\'es. Dans le cas des vari\'et\'es $\Xb(\wb)$,
$\tilde{\Xb}(\wb)$ et $\hat{\Xb}(\wb)$, la comparaison est facile. Les morphismes
canoniques $\Gb\to \tilde{\Gb}$ et $\hat{\Gb}\to\Gb$ induisent des isomorphismes
$\Gb/\Bb\isom \tilde{\Gb}/\tilde{\Bb}$ et $\hat{\Gb}/\hat{\Bb}\isom \Gb/\Bb$, d'o\`u 

\begin{prop}\label{iso xtilde}
Soit $\wb$ une suite d'\eles de $W$. Alors, les $\Gb^F$-vari\'et\'es
$\Xb(\wb)$, $\tilde{\Xb}(\wb)$ et $\hat{\Xb}(\wb)$ sont canoniquement isomorphes.
\end{prop}

Les relations entre les vari\'et\'es $\Yb(\nb)$, $\tilde{\Yb}(\nb)$ et $\hat{\Yb}(\nb)$,
sont un peu plus compliqu\'ees. Commen\c{c}ons par les vari\'et\'es
$\tilde{\Yb}(\nb)$.

Le produit fournit des isomorphismes canoniques
$\tilde{\Gb}^F\times_{\Gb^F}\Gb^F\isom \tilde{\Gb}^F$ et
$\Gb^F\times_{\Tb^{\wb F}}\tilde{\Tb}^{\wb F}\isom \tilde{\Gb}^F$. On en d\'eduit~:

\begin{prop}\label{ytilde}
Soit $\nb$ une suite d'\eles de $N_\Gb(\Tb)$. Alors
on a des isomorphismes canoniques
$$\tilde{\Yb}(\nb)\isom \tilde{\Gb}^F \times_{\Gb^F} \Yb(\nb)
\isom \Yb(\nb)\times_{\Tb^{\wb F}} \tilde{\Tb}^{\wb F}.$$
\end{prop}

\begin{proof}
Gr\^ace \`a la remarque \ref{B=W}, on peut supposer que $\nb=(n_1)$, o\`u
$n_1 \in N_\Gb(\Tb)$. Dans ce cas, la preuve est par exemple analogue \`a
\cite[lemme 2.1.2 (b)]{bonnafe}.
\end{proof}

Soit $\nb$ une suite d'\eles de $N_\Gb(\Tb)$.
La proposition \ref{ytilde} a des cons\'equences au niveau des foncteurs
$\RC_\nb$ et $\tilde{\RC}_\nb$. En effet, on en d\'eduit que
$$\RRR\Gamma_c(\tilde{\Yb}(\nb)) \simeq \Lambda\tilde{\Gb}^F \otimes_{\Lambda\Gb^F}
\RRR\Gamma_c(\Yb(\nb))\hspace{0.7cm}{\mathrm{dans}}~
D^b(\Lambda\tilde{\Gb}^F\otimes \Lambda\Tb^{\wb F})$$
et que
$$\RRR\Gamma_c(\tilde{\Yb}(\nb)) \simeq
\RRR\Gamma_c(\Yb(\nb)) \otimes_{\Lambda\Tb^{\wb F}} \Lambda\tilde{\Tb}^{\wb F}
\hspace{0.7cm}{\mathrm{dans}}~
D^b(\Lambda\Gb^F\otimes\Lambda\tilde{\Tb}^{\wb F}).$$
Par suite, on a
\equat\label{ind et lusztig}
\tilde{\RC}_\nb \circ \Ind_{\Tb^{\wb F}}^{\tilde{\Tb}^{\wb F}} \simeq 
\Ind_{\Gb^F}^{\tilde{\Gb}^F} \circ \micro
\RC_\nb
\endequat
et
\equat\label{res et lusztig}
\RC_\nb \circ \Res_{\Tb^{\wb F}}^{\tilde{\Tb}^{\wb F}} \simeq 
\Res_{\Gb^F}^{\tilde{\Gb}^F} \circ \micro
\tilde{\RC}_\nb.
\endequat

\begin{rem}
\label{res ind l}
Les isomorphismes
(\ref{ind et lusztig}) et (\ref{res et lusztig}) induisent les \'egalit\'es
\equat\label{indu et lusztig}
\tilde{\RRR}_\nb \circ \Ind_{\Tb^{\wb F}}^{\tilde{\Tb}^{\wb F}} = 
\Ind_{\Gb^F}^{\tilde{\Gb}^F} \circ \micro
\RRR_\nb
\endequat
et
\equat\label{rest et lusztig}
\RRR_\nb \circ \Res_{\Tb^{\wb F}}^{\tilde{\Tb}^{\wb F}} = 
\Res_{\Gb^F}^{\tilde{\Gb}^F} \circ \micro \tilde{\RRR}_\nb.
\endequat
\end{rem}

Concernant les rapports entre les vari\'et\'es $\hat{\Yb}(\hat{\nb})$ et $\Yb(\nb)$,
nous avons le r\'esultat suivant (rappelons que $\rho : \hat{\Gb} \to \Gb$
d\'esigne la projection canonique)~:

\begin{prop}\label{yhat}
Soit $\hat{\nb}=(\hat{n}_1,\dots,\hat{n}_r)$ une suite finie d'\eles de $N_{\hat{\Gb}}(\hat{\Tb})$ et
posons $\nb=(\rho(\hat{n}_1),\dots,\rho(\hat{n}_r))$. Alors, l'application $\rho$
induit des \isos de $\Gb^F$-vari\'et\'es-$\Tb^{\wb F}$
$$\hat{\Yb}(\hat{\nb})/{\Cb^F} \isom \Yb(\nb)\ \  \text{ et }\ \ \ 
{\Cb^F}\setminus \hat{\Yb}(\hat{\nb}) \isom \Yb(\nb).$$
\end{prop}

\begin{proof}
Encore une fois, comme pour la proposition \ref{ytilde}, on peut
supposer que $\nb$ est une suite form\'ee d'un seul \'el\'ement. Dans ce cas, 
il suffit d'appliquer la proposition 
\ref{existence quotient}, en rempla\c{c}ant le groupe $\Gb$ par $\hat{\Gb}$, 
$\Pb$ par $\hat{\Bb}$, $\Lb$ par $\hat{\Tb}$, $n$ par $\hat{n}_1$, 
le groupe $\Hb$ par $\Cb$, et le groupe $\Hb'$ par $\{1\}$. Dans ce cas, 
la vari\'et\'e qui \'etait not\'ee $\Yb$ devient la vari\'et\'e $\hat{\Yb}(\hat{\nb})$, 
tandis que $\Xb$ devient $\Yb(\nb)$ et $\Kb=\Cb^F$. Cela compl\`ete la preuve 
de la proposition \ref{yhat}. 
\end{proof}

Si $\Res_{\hat{\Gb}^F}^{\Gb^F} : D^b(\Lambda\Gb^F) \longto 
D^b(\Lambda\hat{\Gb}^F)$
d\'esigne le foncteur de restriction \`a travers $\rho$,
alors la proposition \ref{yhat} montre que
\equat\label{commute GHAT}
\hat{\RC}_{\hat{\nb}} \circ \micro \Res_{\hat{\Tb}^{\wb F}}^{\Tb^{\wb F}} \simeq
\Res_{\hat{\Gb}^F}^{\Gb^F} \circ \RC_\nb.
\endequat

\begin{rem}
\label{qwerty}
Tous les r\'esultats \'enonc\'es dans ce \S\ref{lien lien} restent vrais
m\^eme lorsque le centre de $\tilde{\Gb}$ n'est pas connexe ou lorsque les coracines
de $\hat{\Gb}$ ne sont pas injectives. D'autre part, si $\check{\Gb}$ est un groupe
alg\'ebrique muni d'une isog\'enie $F : \check{\Gb} \to \check{\Gb}$ telle que
$F^\delta$ soit l'\endo de Frobenius associ\'e \`a une $\fq$-structure sur $\check{\Gb}$,
et si $\gamma : \check{\Gb} \to \Gb$ est un morphisme de groupes alg\'ebriques commutant
avec $F$ et induisant un isomorphisme de groupes, alors les groupes finis
$\check{\Gb}^F$ et $\Gb^F$ sont isomorphes via $\gamma$ et $\gamma$ induit des 
\'equivalences \'equivariantes de sites \'etales entre les vari\'et\'es 
$\Yb^{\check{\Gb}}(\check{\nb})$ et $\Yb^\Gb(\nb)$ (avec des notations \'evidentes). 
En particulier, on a 
$\RC_{\check{\nb}}^{\check{\Gb}} \circ \Res_{\check{\Tb}^{\check{\nb} F}}^{\Tb^{\nb F}} 
= \Res_{\check{\Gb}^F}^{\Gb^F} \circ \RC_\nb^\Gb$.  
\end{rem}

\section{Recollements\label{sub recollement}}~

\subsection{Construction} \label{Construction}
{\it Dans ce paragraphe \ref{Construction}, et dans ce paragraphe seulement, 
nous supposerons que $\Gb=\hat{\Gb}$.}
Soient $\vb$ et $\wb$ deux suites
d'\eles de $\bar{S}$ telles que $\vb \le \wb$. Nous noterons $\wb=(s_1,\dots,s_r)$
et $\vb=(s_1^\prime,\dots,s_r^\prime)$. Posons
$$\Xb(\wb;\vb)=\coprod_{\vb \le \yb \le \wb} \Xb(\yb).\index{xwv@$\Xb(\wb;\vb)$}$$
Puisque $\Bb s \Bb \coprod \Bb$ est une sous-vari\'et\'e ferm\'ee lisse
de $\Gb$ pour tout $s \in S$ (c'est m\^eme un \para de $\Gb$), $\Xb(\wb;\vb)$ est une
sous-vari\'et\'e localement ferm\'ee lisse de $(\Gb/\Bb)^r$. Nous allons ici
construire un rev\^etement \'etale (en g\'en\'eral non trivial)
de la vari\'et\'e $\Xb(\wb;\vb)$ en recollant des quotients des vari\'et\'es
$\Yb(\dot{\yb})$ pour $\vb \le \yb \le \wb$.

Soit $\alpha\in \Delta$. Alors, $\Ub \dot{s}_\alpha \Ub =
\Ub \dot{s}_\alpha \Ub_\alpha (\Ub \cap \lexp{s_\alpha}{\Ub})$.
Cela montre que
$$\Ub \dot{s}_\alpha \Ub \Tb_{\alpha^\vee} \coprod \Ub \Tb_{\alpha^\vee}
= \Gb_\alpha(\Ub \cap \lexp{s_\alpha}{\Ub})=\Gb_\alpha\Ub.$$
En particulier, c'est une sous-vari\'et\'e lisse localement ferm\'ee
de $\Gb$. Pour $1 \le i \le r$, on pose
$$\mathbf{\UC}_{\wb,\vb}^{(i)}=\left\{\begin{array}{ll}
\Gb_{\alpha_{\wb,i}}\Ub & {\mathrm{si~}} i \in I_{\wb,\vb}, \\
\Ub \dot{s}_i \Ub & {\mathrm{sinon}}
\end{array}\right.$$
o\`u $I_{\wb,\vb}$ a \'et\'e d\'efini au \S \ref{sousgroupestores}.

On  d\'efinit alors la vari\'et\'e
\begin{multline*}
\Yb'(\wb;\vb)=\bigl\{(g_1\Ub,\dots,g_r\Ub) \in (\Gb/\Ub)^r~|~
g_i^{-1}g_{i+1} \in \mathbf{\UC}^{(i)}_{\wb,\vb}\ \textrm{pour}\
1 \le i \le r-1\\
\mathrm{et}\ \ g_r^{-1}F(g_1) \in \mathbf{\UC}_{\wb,\vb}^{(r)}\}.
\end{multline*}
\index{ywvp@$\Yb'(\wb;\vb)$}
C'est une sous-vari\'et\'e localement ferm\'ee lisse de $(\Gb/\Ub)^r$, sur laquelle
le groupe $\Gb^F$ agit diagonalement par multiplication \`a gauche.
De plus, nous avons un morphisme canonique de $\Gb^F$-vari\'et\'es
$\pi_{\wb,\vb}^\prime : \Yb'(\wb;\vb) \to \Xb(\wb;\vb)$ induit par
la projection $(\Gb/\Ub)^r\to (\Gb/\Bb)^r$.

D'autre part, pour $(g_1 \Ub,\dots,g_r\Ub) \in \Yb'(\wb;\vb)$ et pour
$(t_1,\dots,t_r)
\in \Tb^r$, alors\\
$(g_1 t_1\Ub,\dots,g_r t_r\Ub) \in \Yb'(\wb;\vb)$ \ssi $(t_1,\dots,t_r)
\in \Sb_{\wb,\vb}$.

\begin{lem}\label{separabilite}
L'application canonique $\Yb'(\wb;\vb)/\Sb_{\wb,\vb} \to \Xb(\wb;\vb)$
induite par $\pi'_{\wb,\vb}$
est un \iso de $\Gb^F$-vari\'et\'es.
\end{lem}

\begin{proof}
Puisque le morphisme est bijectif et que $\Xb(\wb;\vb)$ est lisse,
il suffit de v\'erifier que c'est un isomorphisme au-dessus
de l'ouvert $\Xb(\wb)$. Cela r\'esulte alors de la remarque \ref{B=W} et 
de la proposition \ref{existence quotient} appliqu\'ee au groupe $\Gb'=\Gb^r$ 
muni de l'isog\'enie d\'efinie dans la remarque \ref{B=W}, avec
$\Pb=\Bb^r$, $\Lb=\Hb=\Tb^r$, $n=\dot{\wb}$, 
et $\Hb'=\Tb_{\wb,\vb}^{(1)}\times\cdots\Tb_{\wb,\vb}^{(r)}$.
Il faut noter qu'alors $\Kb=\Sb_{\wb,\vb}$. 
\end{proof}

D\'efinissons
\equat\label{definition recollement 1}
\Yb(\wb;\vb)=\Yb'(\wb;\vb)/\Sb_{\wb,\vb}^\circ,\index{ywv@$\Yb(\wb;\vb)$}
\endequat
et notons $\pi_{\wb,\vb} : \Yb(\wb,\vb) \to
\Xb(\wb,\vb)$\index{piwv@$\pi_{\wb,\vb}$} le morphisme
de $\Gb^F$-vari\'et\'es induit par $\pi_{\wb,\vb}^\prime$.
Le lemme \ref{separabilite} montre que $\pi_{\wb,\vb}$ est un
$(\Sb_{\wb,\vb}/\Sb_{\wb,\vb}^\circ)$-torseur.

Consid\'erons maintenant une suite $\yb$ telle que $\vb \le \yb \le \wb$.
Alors, l'injection canonique $\Yb(\dot{\yb}) \injto \Yb'(\wb;\vb)$ est un morphisme
de $\Gb^F$-vari\'et\'es-$\Tb^{\yb F}$~: ici, $\Yb(\wb;\vb)$ est vue comme une
vari\'et\'e-$\Tb^{\yb F}$ via le morphisme $\mu_\yb : \Tb^{\yb F} \injto
\Sb_{\wb,\vb}$ d\'efini dans le \S\ref{application}. Puisque
$$\mu_\yb(\Tb^{\yb F}) \cap \Sb_{\wb,\vb}^\circ = \mu_\yb(N_\yb(Y_{\wb,\vb})),$$
(cf proposition \ref{Sw} (4)), cela induit un morphisme de
$\Gb^F$-vari\'et\'es-$(\Tb^{\yb F}/N_\yb(Y_{\wb,\vb}))$
$$\Yb(\dot{\yb})/N_\yb(Y_{\wb,\vb}) \longto \Yb(\wb;\vb).$$

\begin{lem}\label{clef}
L'image du morphisme $\Yb(\dot{\yb}) \to \Yb(\wb;\vb)$ est \'egale \`a
$\pi_{\wb,\vb}^{-1}(\Xb(\yb))$. Ce morphisme induit un \iso de
$\Gb^F$-vari\'et\'es-$(\Tb^{\yb F}/N_\yb(Y_{\wb,\vb}))$
$$\Yb(\dot{\yb})/N_\yb(Y_{\wb,\vb}) \isom
\pi_{\wb,\vb}^{-1}(\Xb(\yb)).$$
\end{lem}

\begin{proof}
L'application $\Yb(\dot{\yb}) \to \Xb(\yb)$ est surjective. Donc, d'apr\`es
le lemme \ref{separabilite}, l'application
$$\begin{array}{ccc}
\Yb(\dot{\yb}) \times_{\mu_\yb(\Tb^{\yb F})} \Sb_{\wb,\vb} &\isom&
 \pi_{\wb,\vb}^{\prime -1}(\Xb(\yb))\\
(\gamma,\tau)&\mapsto&\gamma.\tau
\end{array}$$
est un isomorphisme de vari\'et\'es. Le r\'esultat d\'ecoule alors de ce que
$$\Sb_{\wb,\vb}=\mu_\yb(\Tb^{\yb F})\cdot\Sb_{\wb,\vb}^\circ,$$
ce qui est d\'emontr\'e dans la proposition \ref{Sw} (2).
\end{proof}

\subsection{Le cas g\'en\'eral} Revenons maintenant au cas g\'en\'eral~:
nous ne
supposons plus que $\Gb=\hat{\Gb}$. Le lecteur aura remarqu\'e
que le lemme \ref{clef} n'est plus n\'ecessairement vrai
(nous avons utilis\'e
la proposition \ref{Sw} (4), qui n'est pas vraie en g\'en\'eral).
C'est pourquoi
nous allons devoir passer par le groupe $\hat{\Gb}$ pour d\'efinir un bon
rev\^etement de la vari\'et\'e $\Xb(\wb;\vb)$.

Notons $\hat{\Yb}'(\wb;\vb)$ (\resp $\hat{\Sb}_{\wb,\vb}$) la vari\'et\'e (\resp le groupe
diagonalisable) d\'efinie (\resp d\'efini) dans le groupe $\hat{\Gb}$
comme l'a \'et\'e $\Yb'(\wb;\vb)$ (\resp $\Sb_{\wb,\vb}$) dans le
paragraphe pr\'ec\'edent. Nous posons alors~:
\equat\label{definition recollement}
\Yb(\wb;\vb)=\hat{\Yb}'(\wb;\vb)/(\hat{\Sb}_{\wb,\vb}^\circ\cdot\Cb^F),
\index{ywv@$\Yb(\wb;\vb)$}
\endequat
le groupe $\Cb^F$ agissant diagonalement sur $(\hat{\Gb}/\hat{\Bb})^r$ par multiplication.
Notons $\pi_{\wb,\vb} : \Yb(\wb;\vb) \to \Xb(\wb;\vb)$\index{piwv@$\pi_{\wb,\vb}$}
le \mor canonique~:
c'est un $(\hat{\Sb}_{\wb,\vb}/(\hat{\Sb}_{\wb,\vb}^\circ\cdot\Cb^F))$-torseur.

Compte tenu des propositions \ref{iso xtilde} et \ref{yhat} et du lemme \ref{clef},
nous d\'eduisons le r\'esultat suivant.

\begin{theo}\label{theo recollement}
Le morphisme de $\Gb^F$-vari\'et\'es
$$\pi_{\wb,\vb} : \Yb(\wb;\vb) \to \Xb(\wb;\vb)$$
est un $\Tb^{\wb F}/N_\wb(Y_{\wb,\vb})$-torseur \'etale. De plus,
pour $\vb \le \yb \le \wb$,
l'application canonique $\Yb(\dot{\yb}) \to \Yb(\wb;\vb)$ induit un
\iso $\Yb(\dot{\yb})/N_\yb(Y_{\wb,\vb}) \isom \pi_{\wb,\vb}^{-1}(\Xb(\yb))$
de $\Gb^F$-vari\'et\'es-$(\Tb^{\yb F}/N_\yb(Y_{\wb,\vb}))$ rendant
commutatif le diagramme suivant
\equat\label{diagramme recollement}
\diagram
\Yb(\dot{\yb})/N_\yb(Y_{\wb,\vb}) \rrto^{\DS{\sim}} \ddrrto_{\DS{\pi_\yb}} &&
\pi_{\wb,\vb}^{-1}(\Xb(\yb)) \ddto \rrto|<\ahook&& \Yb(\wb;\vb)
\ddto^{\DS{\pi_{\wb,\vb}}} \\
&&&& \\
&& \Xb(\yb) \rrto|<\ahook && \Xb(\wb;\vb)\\
\enddiagram
\endequat
\end{theo}

\begin{rem}
Dans le diagramme \ref{diagramme recollement}, le rev\^etement
\'etale $\Yb(\dot{\yb})/N_\xb(Y_{\wb,\vb}) \to \Xb(\yb)$ est un
$\Tb^{\yb F}/N_\yb(Y_{\wb,\vb})$-torseur, alors que le rev\^etement \'etale
$\pi_{\wb,\vb} : \pi_{\wb,\vb}^{-1}(\Xb(\yb)) \to \Xb(\yb)$
est quant \`a lui
un $\Tb^{\wb F}/N_\wb(Y_{\wb,\vb})$-torseur. Ces descriptions sont
compatibles, via l'isomorphisme canonique $\Tb^{\yb F}/N_\yb(Y_{\wb,\vb})
\isom \Tb^{\wb F}/N_\wb(Y_{\wb,\vb})$ (cf proposition \ref{Sw}).
\end{rem}

\subsection{Recollement de recollements} Soient $\vb \le \xb \le \yb \le \wb$.
La vari\'et\'e $\Xb(\yb;\xb)$
est une sous-vari\'et\'e localement ferm\'ee de $\Xb(\wb;\vb)$ et
$\pi_{\wb,\vb}^{-1}(\Xb(\yb;\xb)) \to \Xb(\yb;\xb)$ est un
$\Tb^{\wb F}/N_\wb(Y_{\wb,\vb})$-torseur \'etale.

D'autre part, $\pi_{\yb,\xb} : \Yb(\yb;\xb) \to
\Xb(\yb;\xb)$ est un $\Tb^{\yb F}/N_\yb(Y_{\yb,\xb})$-torseur \'etale. D'apr\`es le
corollaire \ref{reseau vxyw},
on a $Y_{\yb,\xb} \incl Y_{\wb,\vb}$. Par cons\'equent,
$\Tb^{\yb F}/N_\yb(Y_{\wb,\vb})$ est un quotient de
$\Tb^{\yb F}/N_\yb(Y_{\yb,\xb})$ et
le \mor canonique $\Yb(\yb;\xb)/\bigl(N_\yb(Y_{\wb,\vb})/N_\yb(Y_{\yb,\xb})\bigr) \to
\Xb(\yb;\xb)$ est un $\Tb^{\yb F}/N_\yb(Y_{\wb,\vb})$-torseur \'etale.
Mais d'apr\`es la proposition \ref{reseau} (4), on a un \iso
canonique
$$\Tb^{\yb F}/N_\yb(Y_{\wb,\vb}) \simeq \Tb^{\wb F}/N_\wb(Y_{\wb,\vb}).$$

Par cons\'equent, nous avons construit deux $\Tb^{\wb F}/N_\wb(Y_{\wb,\vb})$-torseurs
\'etales au-dessus de $\Xb(\yb;\xb)$. Le corollaire
suivant montre qu'ils sont en fait isomorphes.

\begin{coro}\label{iso recollement}
Le \mor canonique $\Yb(\yb;\xb) \to \Yb(\wb;\vb)$ induit un \iso $\Gb^F$-\'equivariant
de $\Tb^{\wb F}/N_\wb(Y_{\wb,\vb})$-torseurs \'etales
$$\Yb(\yb;\xb)/\bigl(N_\yb(Y_{\wb,\vb})/N_\yb(Y_{\yb,\xb})\bigr) \isom
\pi_{\wb,\vb}^{-1}(\Xb(\yb;\xb))$$
au-dessus de $\Xb(\yb;\xb)$.
\end{coro}

\begin{proof}
Le \mor canonique provient de l'inclusion
$\Yb'(\yb;\xb)\hookrightarrow \Yb'(\wb;\vb)$. Il est
clair que son image est \'egale \`a la
vari\'et\'e lisse $\pi_{\wb,\vb}^{-1}(\Xb(\yb;\xb))$.
Pour d\'emontrer le corollaire \ref{iso recollement}, il suffit de d\'emontrer
que le \mor
$$\Yb(\yb;\xb)/\bigl(N_\yb(Y_{\wb,\vb})/N_\yb(Y_{\yb,\xb})\bigr) \longto
\pi_{\wb,\vb}^{-1}(\Xb(\yb;\xb))$$
ainsi d\'efini induit un \iso au-dessus d'un ouvert de $\Xb(\yb;\xb)$.
On peut par exemple prendre l'ouvert $\Xb(\yb)$ de $\Xb(\yb;\xb)$.
Mais ce fait r\'esulte de la commutativit\'e du diagramme
du th\'eor\`eme \ref{theo recollement}.
\end{proof}

\bigskip

Nous r\'esumons le corollaire \ref{iso recollement} dans le diagramme commutatif
suivant~:
$$\diagram
\Yb(\yb;\xb) \rrto \ddrrto \xto[4,2]_{\DS{\pi_{\yb,\xb}}} &&
\Yb(\yb;\xb)/\bigl(N_\yb(Y_{\wb,\vb})/N_\yb(Y_{\yb,\xb})\bigr)
\ddto^{\DS{\sim}} && \\
&&&&\\
&& \pi_{\wb,\vb}^{-1}(\Xb(\yb;\xb)) \rrto|<\ahook \ddto && \Yb(\wb;\vb)
\ddto^{\DS{\pi_{\wb,\vb}}}\\
&&&& \\
&& \Xb(\yb;\xb) \rrto|<\ahook && \Xb(\wb;\vb).
\enddiagram$$

\subsection{Remarques}
Pour construire le rev\^etement $\pi_{\wb,\vb} : \Yb(\wb;\vb) \to
\Xb(\wb;\vb)$, il a \'et\'e n\'ecessaire de passer par un groupe $\hat{\Gb}$
dont les coracines sont injectives et tel qu'il existe un tore central
$F$-stable $\Cb$ de $\hat{\Gb}$ tel que $\Gb\simeq \hat{\Gb}/\Cb$. Disons dans ce paragraphe
qu'une paire $(\hat{\Gb},\Cb)$ v\'erifiant ces propri\'et\'es est {\it sympathique}.

Le lecteur pourra v\'erifier que, dans le cas o\`u $\Gb=\Pb\Gb\Lb_2(\FM)$
et o\`u $p \not= 2$, alors on ne peut effectuer directement la construction
dans $\Gb$. En effet, on obtiendrait ainsi un rev\^etement connexe
de $\Yb(s;1)$ (ici, $s$ d\'esigne l'unique r\'eflexion simple de $W$), alors
que l'on attend un rev\^etement non connexe. Il est donc n\'ecessaire
de passer par $\hat{\Gb}=\Gb\Lb_2(\FM)$.

Soient $(\hat{\Gb}_1,\Cb_1)$ et $(\hat{\Gb}_2,\Cb_2)$ deux
paires sympathiques et $\hat{\Gb}_0$ le produit fibr\'e de $\hat{\Gb}_1$ et $\hat{\Gb}_2$
au-dessus de $\Gb$. On note $\hat{\Yb}_i(\wb;\vb)$ la vari\'et\'e
d\'efinie dans $\hat{\Gb}_i$ comme $\Yb(\wb;\vb)$ a \'et\'e d\'efinie dans
$\Gb$ lorsque les coracines de $\Gb$ sont injectives.
La paire $(\hat{\Gb}_0,\Cb_1 \times \Cb_2)$
est sympathique et
$\hat{\Yb}_i(\wb;\vb) \isom \hat{\Yb}_0(\wb;\vb)/\Cb_j^F$, o\`u $\{i,j\}=\{1,2\}$.
Cela fournit un isomorphisme
$\hat{\Yb}_1(\wb;\vb)/\Cb_1^F \simeq \hat{\Yb}_2(\wb;\vb)/\Cb_2^F$. On en d\'eduit que
la vari\'et\'e
$\Yb(\wb;\vb)$ ne d\'epend pas du choix de $\hat{\Gb}$, \`a isomorphisme
unique pr\`es.

\section{Calculs de monodromie\label{section monodromie}}~

Nous allons appliquer ici les r\'esultats du paragraphe pr\'ec\'edent
au calcul des images directes sup\'erieures (dans une compactification
lisse adapt\'ee) des faisceaux localement constants sur les vari\'et\'es
$\Xb(\wb)$ provenant du rev\^etement \'etale $\Yb(\dot{\wb}) \to \Xb(\wb)$.
Le r\'esultat (le th\'eor\`eme \ref{theo mono}),
g\'en\'eralise \cite[lemme 9.13]{delu}. Nous obtiendrons comme cons\'equences
des r\'esultats sur les foncteurs de Deligne-Lusztig (cf \S\ref{sub f}).

\medskip

\begin{centerline}{\it Dans toute cette section,
nous supposons que $\Lambda$ est un corps.}\end{centerline}

\medskip

On fixe dans cette section une suite $\wb$ d'\'el\'ements de $\bar{S}$.

Soient $\vb$, $\xb$
et $\yb$ des suites d'\eles de $\bar{S}$ telles que
$\vb \le \xb \le \yb \le \wb$.
Nous noterons $j_{\yb;\xb}^{\wb;\vb} : \Xb(\yb;\xb) \injto
\Xb(\wb;\vb)$\index{jyxwv@$j_{\yb;\xb}^{\wb;\vb}$}
 l'injection canonique. Pour simplifier, nous remplacerons,
lorsque $\yb=\xb$, l'indice
ou l'exposant $\yb;\xb$ par $\yb$. D'autre part, nous posons
\begin{eqnarray*}
\Xb(\overline{\xb})&=&\coprod_{\xb' \le \xb} \Xb(\xb') \\
(&=& \Xb(\xb;(1,\dots,1))\hspace{1cm}).
\index{xxover@$\overline{\Xb}(\xb)$}
\end{eqnarray*}
Si $\vb=(1,\dots,1)$, nous remplacerons l'indice ou l'exposant $\xb;\vb$ par
${\overline{\xb}}$. Par exemple, $j_\yb^{\overline{\xb}}$ d\'esigne l'injection canonique
$\Xb(\yb) \injto \Xb(\overline{\xb})$.

La compactification de $\Xb(\xb)$ ainsi construite poss\`ede de bonnes
propri\'et\'es \cite[lemme 9.11]{delu}~:

\begin{lem}\label{DCN DL}
La vari\'et\'e $\Xb(\overline{\xb})$ est projective
lisse et le compl\'ementaire de
$\Xb(\xb)$ est un diviseur \`a croisements normaux
$\bigcup_\vb\Xb(\overline{\vb})$ o\`u $\vb$ d\'ecrit l'ensemble des
suites telles que $\vb \le \xb$ et $l(\xb)-l(\vb)=1$.
\end{lem}

Nous noterons
$$\FC^{\yb;\xb} :
D^b(\Lambda\Tb^{\yb F}/N_\yb(Y_{\yb,\xb})) \to D^b_\Lambda(\Xb(\yb;\xb))
\index{fyx@$\FC^{\yb;\xb}$}$$
le foncteur
$\FC_{\Tb^{\yb F}/N_\yb(Y_{\yb;\xb})}^{\Yb(\yb;\xb)}$ d\'efini au \S\ref{sub foncteur}.
Nous fixons aussi dans cette section un \car lin\'eaire
$\theta:\Tb^{\wb F}\to \Lambda^\times$.
Si $\wb_\theta \le \xb \le
\yb \le \wb$, il r\'esulte de \S \ref{sousgroupestores}
que $\theta$ d\'efinit
un \car  lin\'eaire $\theta_{\yb;\xb}$ du groupe $\Tb^{\yb F}/N_\yb(Y_{\yb,\xb})$
(pour la d\'efinition de $\wb_\theta$, voir \S\ref{sousgroupestores}). Nous noterons
$\FC_\theta^{\yb;\xb}$\index{fthetayx@$\FC_\theta^{\yb;\xb}$} le faisceau localement constant
$\FC^{\yb;\xb}(\Lambda_{\theta_{\yb;\xb}})$
sur la vari\'et\'e $\Xb(\yb;\xb)$. Compte tenu du
corollaire \ref{iso recollement} et de la proposition \ref{Sw} (3),
$\FC_\theta^{\yb;\xb}$ est la
restriction de $\FC_\theta^{\wb;\wb_\theta}$ \`a $\Xb(\yb;\xb)$. En d'autres
termes,
\equat\label{r}
\FC_\theta^{\yb;\xb} \simeq (j_{\yb;\xb}^{\wb;\wb_\theta})^*\FC_\theta^{\wb;\wb_\theta}.
\endequat
Comme expliqu\'e dans l'introduction de cette section, nous allons calculer
ici les faisceaux $\RRR^i(j_\wb^{\overline{\wb}})_* \FC_\theta^\wb$. La proposition 
suivante combine un r\'esultat fondamental de Deligne et Lusztig
\cite[lemme 9.13]{delu}
et une cons\'equence imm\'ediate de la construction
des vari\'et\'es $\Yb(\wb;\vb)$.

\begin{prop}\label{913}
Soient $\vb$, $\vb'$ et $\xb$ des suites d'\'el\'ements de $\bar{S}$
telles que $\wb_\theta \le \vb\le \xb$
et $\vb'<\vb$.

\tete{1} Si $\vb$ est le plus petit \'el\'ement tel que
$\vb'\le\vb$ et $\wb_\theta\le\vb$, alors le faisceau
$\FC_\theta^\xb$ se ramifie le long de chaque composante du diviseur
$\Xb(\xb;\vb')-\Xb(\xb;\vb)$ et
le morphisme
canonique $(j_{\xb;\vb}^{\xb;\vb'})_! \FC_\theta^{\xb;\vb}
\isom \RRR(j_{\xb;\vb}^{\xb;\vb'})_* \FC_\theta^{\xb;\vb}$
est un isomorphisme.

\tete{2} Si $\wb_\theta \le \vb'$, alors le faisceau
$\FC_\theta^\xb$ ne se ramifie pas le long du diviseur
$\Xb(\xb;\vb')-\Xb(\xb;\vb)$.
\end{prop}

\begin{proof}
On a
$$\Xb(\xb;\vb')-\Xb(\xb;\vb)=\bigcup_{\xb'}\Xb(\xb';\vb')$$
o\`u $\vb'\le\xb'\le\xb$, $\vb\not\le\xb'$ et $l(\xb')=l(\xb)-1$.

Sous l'hypoth\`ese de (1), on a $\wb_\theta\not\le \xb'$.
D'apr\`es
\cite[lemme 9.13]{delu} (la preuve, faite pour des $\bar{\QM}_\ell$-faisceaux
reste valable pour des $\Lambda$-faisceaux),
 le faisceau $\FC_\theta^\xb$ se ramifie le long du
diviseur $\Xb(\overline{\xb'})$ de $\bar{\Xb}(\xb)$, i.e., la restriction de
$(j_\xb^{\overline{\xb}})_* \FC_\theta^\xb$ \`a $\Xb(\overline{\xb'})$ est nulle.
On en d\'eduit que 
$\FC_\theta^{\xb;\vb}$ se ramifie le long du diviseur $\Xb(\xb';\vb')$. Par
cons\'equent, le morphisme
canonique $(j_{\xb;\vb}^{\xb;\vb'})_! \FC_\theta^{\xb;\vb}
\to (j_{\xb;\vb}^{\xb;\vb'})_* \FC_\theta^{\xb;\vb}$
est un isomorphisme. L'annulation des
$\RRR^i(j_{\xb;\vb}^{\xb;\vb'})_* \FC_\theta^{\xb;\vb}$ pour $i>0$ r\'esulte
de \cite[exemple 1.19.1 p.180]{SGA4.5}
(on utilise ici le lemme \ref{DCN DL} et la propri\'et\'e de
$\FC_\theta^{\xb;\vb}$ d'\^etre mod\'er\'ement ramifi\'e), d'o\`u (1).

(2) d\'ecoule directement de \ref{r}. 
%
%
\end{proof}

D\'efinissons maintenant quelques
sous-cat\'egories de $D^b_\Lambda(\Xb(\overline{\wb}))$. Prenons
$\wb_\theta \le \yb$. On pose

\begin{align*}
A_{\le \yb}^1(\theta)&=\{ (j_{\xb}^{{\overline{\wb}}})_! \FC_\theta^\xb~|~\wb_\theta \le \xb
\le \yb\},\\
A_{<\yb}^1(\theta)&=\{ (j_{\xb}^{{\overline{\wb}}})_! \FC_\theta^\xb~|~\wb_\theta \le \xb<\yb\},\\
A_{\le \yb}^2(\theta)&=\{ (j_{\xb;\xb'}^{{\overline{\wb}}})_! \FC_\theta^{\xb;\xb'}~|~
\wb_\theta \le \xb' \le \xb \le \yb\},\\
A_{<\yb}^2(\theta)&=\{ (j_{\xb;\xb'}^{{\overline{\wb}}})_! \FC_\theta^{\xb;\xb'}~|~
\wb_\theta \le \xb' \le \xb <\yb\},\\
A_{\le \yb}^3(\theta)&=\{ (j_{\xb;\wb_\theta}^{{\overline{\wb}}})_! \FC_\theta^{\xb;\wb_\theta}~|~
\wb_\theta \le \xb \le \yb\},\\
A_{<\yb}^3(\theta)&=\{ (j_{\xb;\wb_\theta}^{{\overline{\wb}}})_! \FC_\theta^{\xb;\wb_\theta}~|~
\wb_\theta \le \xb < \yb\},\\
A_{\le \yb}^4(\theta)&=\{ \RRR(j_{\xb}^{{\overline{\wb}}})_* \FC_\theta^\xb~|~\wb_\theta \le \xb
\le \yb\},\\
A_{<\yb}^4(\theta)&=\{ \RRR(j_{\xb}^{{\overline{\wb}}})_* \FC_\theta^\xb~|~\wb_\theta \le \xb
< \yb\},\\
A_{\le \yb}^5(\theta)&=\{ \RRR(j_{\xb;\xb'}^{{\overline{\wb}}})_* \FC_\theta^{\xb;\xb'}~|~
\wb_\theta \le \xb' \le \xb \le \yb\},\\
 A_{<\yb}^5(\theta)&=\{ \RRR(j_{\xb;\xb'}^{{\overline{\wb}}})_* \FC_\theta^{\xb;\xb'}~|~
\wb_\theta \le \xb' \le \xb <\yb\},\\
A_{\le \yb}^6(\theta)&=\{ \RRR(j_{\xb;\wb_\theta}^{{\overline{\wb}}})_* \FC_\theta^{\xb;\wb_\theta}~|~
\wb_\theta \le \xb \le \yb\},\\
\textrm{et}\ \ \ \ 
A_{<\yb}^6(\theta)&=\{ \RRR(j_{\xb;\wb_\theta}^{{\overline{\wb}}})_* \FC_\theta^{\xb;\wb_\theta}~|~
\wb_\theta \le \xb < \yb\}.
\end{align*}
Si $(i,\ast) \in \{1,2,3,4,5,6\} \times \{\le \yb,<\yb\}$,
nous noterons $\AC_\ast^i(\theta)$ la sous-cat\'egorie de
$D^b_\Lambda(\Xb(\overline{\wb}))$ engendr\'ee par $A_\ast^i(\theta)$.

\begin{prop}\label{egalite derivee}
Pour $\ast \in \{\le \yb,< \yb\}$, on a 
$$\AC_\ast^1(\theta)=\AC_\ast^2(\theta)=
\AC_\ast^3(\theta)=\AC_\ast^4(\theta)=\AC_\ast^5(\theta)=\AC_\ast^6(\theta).$$
\end{prop}

\begin{proof}
Notons pour commencer que le morphisme canonique
$$(j_{\xb;\wb_\theta}^{{\overline{\wb}}})_! \FC_\theta^{\xb;\wb_\theta}\isom
\RRR(j_{\xb;\wb_\theta}^{{\overline{\wb}}})_* \FC_\theta^{\xb;\wb_\theta}$$
est un isomorphisme (proposition \ref{913} (1)),
donc $\AC_\ast^3(\theta)=\AC_\ast^6(\theta)$.

\medskip
Notons que $\AC_\ast^1(\theta)\incl \AC_\ast^2(\theta)$ et
$\AC_\ast^3(\theta)\incl \AC_\ast^2(\theta)$.

Prenons $\wb_\theta\le \vb'\le\vb\le\xb\le\wb$ avec
$l(\vb)=l(\vb')+1$. Soit $\xb'\le \xb$ tel que
$\vb'\le\xb'$, $\vb\not\le\xb'$ et $l(\xb')=l(\xb)-1$.
On a $\Xb(\xb;\vb')-\Xb(\xb;\vb)=\Xb(\xb';\vb')$, d'o\`u on d\'eduit
une suite exacte
$$0\to (j_{\xb;\vb}^{\overline{\wb}})_!\FC_\theta^{\xb;\vb}\to
(j_{\xb;\vb'}^{\overline{\wb}})_! \FC_\theta^{\xb;\vb'}\to
 (j_{\xb';\vb'}^{\overline{\wb}})_!\FC_\theta^{\xb';\vb'}\to 0.\leqno{(1)}$$

Par cons\'equent, 
$(j_{\xb;\vb'}^{\overline{\wb}})_! \FC_\theta^{\xb;\vb'}$ est dans la
sous-cat\'egorie engendr\'ee par les
$(j_{\xb'';\vb''}^{\overline{\wb}})_!\FC_\theta^{\xb'';\vb''}$, o\`u
$l(\xb'')-l(\vb'')<l(\xb)-l(\vb')$ et
$\wb_\theta\le \vb''\le\xb''\le\xb$. Par r\'ecurrence, on en
d\'eduit que $\AC_\ast^2(\theta)\incl\AC^1_\ast(\theta)$.

On d\'eduit aussi de la suite exacte (1) que
$(j_{\xb;\vb}^{\overline{\wb}})_!\FC_\theta^{\xb;\vb}$ est dans la
sous-cat\'egorie engendr\'ee par les
$(j_{\xb'';\vb''}^{\overline{\wb}})_!\FC_\theta^{\xb'';\vb''}$, o\`u
$\vb''<\vb$ et $\xb''\le\xb$. Par r\'ecurrence, on en
d\'eduit que $\AC_\ast^2(\theta)\incl\AC^3_\ast(\theta)$.

On a donc d\'emontr\'e que
$\AC^1_\ast(\theta)=\AC_\ast^2(\theta)=\AC^3_\ast(\theta)$.

\medskip
Le th\'eor\`eme de puret\'e cohomologique 
\cite[chapitre VI, th\'eor\`eme 5.1]{milne} fournit un triangle distingu\'e
$$
\FC_\theta^{\xb;\vb'} \to \RRR(j_{\xb;\vb}^{\xb;\vb'})_* \FC_\theta^{\xb;\vb}
\to (j_{\xb';\vb'}^{\xb;\vb'})_!\FC_\theta^{\xb';\vb'}[-1] \rightsquigarrow$$
d'o\`u
$$
\RRR(j_{\xb;\vb'}^{\overline{\wb}})_*\FC_\theta^{\xb;\vb'} \to
\RRR(j_{\xb;\vb}^{\overline{\wb}})_* \FC_\theta^{\xb;\vb} \to 
\RRR(j_{\xb';\vb'}^{\overline{\wb}})_*\FC_\theta^{\xb';\vb'}[-1] \rightsquigarrow.$$
On prouve alors, comme pr\'ec\'edemment, que
$\AC^4_\ast(\theta)=\AC_\ast^5(\theta)=\AC^6_\ast(\theta)$.
\end{proof}

\begin{rem}
Notons que la dualit\'e de Verdier \'echange les cat\'egories 
$\AC_*^i(\theta)$ et $\AC_*^i(\theta^{-1})$.
\end{rem}

\begin{coro}\label{image directe 2}
Le c\^one du \mor canonique de complexes
$(j_\wb^{\overline{\wb}})_! \FC_\theta^\wb\longto \RRR(j_\wb^{\overline{\wb}})_*
\FC_\theta^\wb$ est dans la sous-cat\'egorie $\AC_{<\wb}^i(\theta)$
pour tout $i\in \{1,2,3,4,5,6\}$.
\end{coro}

\begin{proof}
D'apr\`es la proposition \ref{egalite derivee}, ce c\^one appartient \`a la
cat\'egorie $\AC_{\le \wb}^1(\theta)$. D'autre part, son support
est contenu dans le compl\'ementaire de $\Xb(\wb)$ dans
$\Xb(\overline{\wb})$, d'o\`u le r\'esultat.
\end{proof}

\begin{theo}\label{theo mono}
Soit $i$ un entier naturel. Alors
$$(j_\vb^{\overline{\wb}})^*
\bigl(\RRR^i (j_\wb^{\overline{\wb}})_* \FC_\theta^\wb \bigr) \simeq
\left\{\begin{array}{ll}
(\FC_\theta^\vb)^{\oplus {l(\wb)-l(\vb) \choose i}}
& {\mathit{si}~} \vb \in \IC(\wb,\theta), \\
&\\
0 & {\mathit{sinon}}.
\end{array}\right.$$
\end{theo}

\begin{proof}
Factorisons l'immersion ouverte $j_\wb^{\overline{\wb}}$ de la mani\`ere suivante
$$\diagram
\Xb(\wb) \rrto^{j_1} && \Xb(\wb;\wb_\theta) \rrto^{j_2} && \overline{\Xb}(\wb).
\enddiagram$$
D'apr\`es la proposition \ref{egalite derivee}, le support du complexe
$\RRR (j_\wb^{\overline{\wb}})_* \FC_\theta^\wb$ est contenu dans $\Xb(\wb;\wb_\theta)$, ce qui
montre que
$$\RRR (j_\wb^{\overline{\wb}})_* \FC_\theta^\wb\simeq (j_2)_!~\RRR(j_1)_* \FC_\theta^\wb.\leqno{(1)}$$
Mais, $\FC_\theta^\wb\isom j^*_1 \FC_\theta^{\wb;\wb_\theta}$,
par cons\'equent, d'apr\`es le lemme \ref{*},
$$\RRR(j_1)_* \FC_\theta^\wb\isom
(\RRR(j_1)_* \Lambda_{\Xb(\wb)}) \otimes_\Lambda \FC_\theta^{\wb;\wb_\theta}.
\leqno{(2)}$$
Le th\'eor\`eme \ref{theo mono} d\'ecoule alors imm\'ediatement de (1) et (2),
en appliquant les lemmes \ref{DCN DL} et \ref{DCN}.
\end{proof}

\section{Induction de Deligne-Lusztig et s\'eries rationnelles\label{sec fdlesr}}~

Nous allons appliquer les r\'esultats de la section pr\'ec\'edente aux
foncteurs de Lusztig.
Pour cela, nous ferons l'hypoth\`ese suivante.

\medskip

\subsection{Orthogonalit\'e}
Soient $(\wb,\theta)$ et $(\wb',\theta')$ deux \eles de
$\nablab(\Tb,W,F)$.

Le r\'esultat suivant g\'en\'eralise \cite[th\'eor\`eme 6.2]{delu} qui
consid\`ere le cas $\Lambda=K$ et la conjugaison g\'eom\'etrique.

\begin{theo}\label{00}
Si $(\wb,\theta)$ et $(\wb',\theta')$ sont dans deux s\'eries
rationnelles diff\'erentes, alors
$$(\RC_{\dot{\wb}} \Lambda\Tb^{\wb F} e_{\theta^{-1}})^\opp \otimes_{\Lambda\Gb^F}^\LLL
\RC_{\dot{\wb}'} \Lambda\Tb^{\wb'F}e_{\theta'}=0.$$
\end{theo}

\begin{proof}
Par la formule des coefficients universels, 
il suffit de d\'emontrer le th\'eor\`eme lorsque $\Lambda$
est un corps, ce que nous supposerons dans la suite.
Commen\c{c}ons par d\'emontrer la proposition lorsque le centre de $\Gb$
est connexe (lorsque $\Lambda=K$, c'est
\cite[th\'eor\`eme 6.2]{delu}). Dans ce cas,
$(\wb,\theta)$ et $(\wb',\theta')$ ne sont pas
g\'eom\'etriquement conjugu\'es (cf remarque \ref{remarque conju}).

D'apr\`es le lemme \ref{lemme foncteur} et la proposition \ref{goodguy},
on a un isomorphisme dans
$D^b(\Lambda(\Tb^{\wb F}\times \Tb^{\wb'F}))$
$$
\RRR\Gamma_c((\Yb(\dot{\wb})\times\Yb'(\dot{\wb}'))/\Delta(\Gb^F))
\isom \RRR\Gamma_c(\Yb(\dot{\wb}))^\opp\otimes_{\Lambda\Gb^F}^\LLL
\RRR\Gamma_c(\Yb(\dot{\wb}')).$$

La preuve est alors similaire \`a celle de \cite[th\'eor\`eme 6.2]{delu}~:
on se ram\`ene \`a prouver l'assertion suivante.
 
Soit $\Xb$ une vari\'et\'e munie d'une action d'un groupe alg\'ebrique
$\Hb$ contenant $\Tb^{\wb F}\times \Tb^{\wb'F}$. On suppose le caract\`ere
$\theta^{-1}\times\theta'$ non trivial sur
$Q=(\Tb^{\wb F}\times \Tb^{\wb'F})\cap \Hb^\circ$.
Alors, $(e_{\theta^{-1}}\otimes_\Lambda e_{\theta'})\RRR\Gamma_c(\Xb)=0$.
 
\smallskip
D\'emontrons donc cette assertion. D'apr\`es \cite[proposition 6.4]{delu},
le groupe connexe $\Hb^\circ$ agit trivialement sur $H^*_c(\Xb)$.
Puisque $\Res_Q^{\Tb^{\wb F} \times \Tb^{\wb' F}}(\theta^{-1}\times\theta')$ 
n'est pas dans le bloc principal de $\Lambda Q$, on a
$$\RRR\Hom^\bullet_{\Lambda Q}(\Lambda\Tb^{\wb F}e_{\theta^{-1}}\otimes_\Lambda
\Lambda\Tb^{\wb'F} e_{\theta'},\RRR\Gamma_c(\Xb))=0,$$ donc
$$\RRR\Hom^\bullet_{\Lambda(\Tb^{\wb F}\times \Tb^{\wb'F})}(
\Ind_Q^{\Tb^{\wb F}\times \Tb^{\wb'F}}\left(
 \Lambda\Tb^{\wb F}e_{\theta^{-1}}\otimes_\Lambda
\Lambda\Tb^{\wb'F} e_{\theta'}\right),\RRR\Gamma_c(\Xb))=0$$
 et a fortiori
$$\RRR\Hom^\bullet_{\Lambda(\Tb^{\wb F}\times \Tb^{\wb'F})}
(\Lambda\Tb^{\wb F}e_{\theta^{-1}}\otimes_\Lambda \Lambda\Tb^{\wb'F} e_{\theta'},
\RRR\Gamma_c(\Xb))=0.$$

\medskip

Revenons maintenant \`a $\Gb$ quelconque. Soit $A$ (\resp $A'$)
l'ensemble des \cars lin\'eaires de $\tilde{\Tb}^{\wb F}$ (\resp $\tilde{\Tb}^{\wb'F}$)
\'etendant $\theta$ (\resp $\theta'$).

Le complexe
$$C=\bigl( \RC_{\dot{\wb}}(\Lambda\Tb^{\wb F} e_{\theta^{-1}})\bigr)^\opp
\otimes_{\Lambda\Gb^F}^\LLL
\RC_{\dot{\wb}'} (\Lambda\Tb^{\wb'F} e_{\theta'})$$
est un facteur direct du complexe
$$\tilde{C}=\bigl(\Ind_{\Gb^F}^{\tilde{\Gb}^F}
\RC_{\dot{\wb}}(\Lambda\Tb^{\wb F} e_{\theta^{-1}})\bigr)^\opp
\otimes_{\Lambda\tilde{\Gb}^F}^\LLL
\bigl( \Ind_{\Gb^F}^{\tilde{\Gb}^F} \RC_{\dot{\wb}'} (\Lambda\Tb^{\wb'F} e_{\theta'})
\bigr).$$
De plus, d'apr\`es (\ref{ind et lusztig}), on a
$$\Ind_{\Gb^F}^{\tilde{\Gb}^F} \RC_{\dot{\wb}} (\Lambda\Tb^{\wb F} e_{\theta^{-1}})
\simeq \tilde{\RC}_{\dot{\wb}} (\Lambda\tilde{\Tb}^{\wb F} e_{\theta^{-1}}) $$
$$\Ind_{\Gb^F}^{\tilde{\Gb}^F} \RC_{\dot{\wb}'} (\Lambda\Tb^{\wb'F} e_{\theta'})
\simeq \tilde{\RC}_{\dot{\wb}} (\Lambda\tilde{\Tb}^{\wb'F} e_{\theta'}).
\leqno{\mathrm{et}}$$
Mais
$$\Lambda\tilde{\Tb}^{\wb F} e_{\theta^{-1}} = 
\bigoplus_{\tilde{\theta} \in A} \Lambda\tilde{\Tb}^{\wb F} e_{\tilde{\theta}^{-1}}$$
$$\Lambda\tilde{\Tb}^{\wb'F} e_{\theta'} = 
\bigoplus_{\tilde{\theta}' \in A'} \Lambda\tilde{\Tb}^{\wb'F}
e_{\tilde{\theta}'}.\leqno{\mathrm{et}}$$
Par cons\'equent,
$$\tilde{C}\simeq\bigoplus_{\substack{\tilde{\theta} \in A\\
\tilde{\theta}' \in A'}}
\bigl( \tilde{\RC}_{\dot{\wb}}(\Lambda\tilde{\Tb}^{\wb F} 
e_{\tilde{\theta}^{-1}})\bigr)^\opp
\otimes_{\Lambda\tilde{\Gb}^F}^\LLL
\tilde{\RC}_{\dot{\wb}'} (\Lambda\tilde{\Tb}^{\wb'F} e_{\tilde{\theta}'}).$$
Puisque $(\wb,\theta)$ et $(\wb',\theta')$ ne sont pas dans la m\^eme
s\'erie rationnelle, les paires $(\wb,\tilde{\theta})$ et $(\wb',\tilde{\theta}')$ ne
sont pas g\'eom\'etriquement conjugu\'ees dans $\tilde{\Gb}$ pour tous
$\tilde{\theta} \in A$ et $\tilde{\theta}' \in A'$. On a alors
$\tilde{C}=0$ d'apr\`es la premi\`ere partie de la preuve,
ce qui implique que $C=0$.
\end{proof}

\begin{coro}\label{0000}
Si $(\wb,\theta)$ et $(\wb',\theta')$ appartiennent \`a deux s\'eries
rationnelles diff\'erentes de
$\nablab(\Tb,W,F)$, alors
$$(\RC_{\dot{\wb}} \Lambda_{\theta^{-1}})^\opp \otimes_{\Lambda\Gb^F}^\LLL
\RC_{\dot{\wb}'} \Lambda_{\theta'}=0.$$
\end{coro}

\subsection{Profondeur\label{sub f}}
{\it Dans toute cette section \S \ref{sub f},
nous supposerons que $\Lambda$ est un corps.}

\smallskip

En appliquant le foncteur
$\RRR\Gamma(\Xb(\overline{\wb}),-)$ \`a l'\'enonc\'e du corollaire
\ref{image directe 2}, on obtient~:

\begin{theo}\label{theo cone}
Soit $\wb\in\Sigma(\bar{S})$ et $\theta$ un caract\`ere lin\'eaire de
$\Tb^{\wb F}$. Alors, le c\^one du \mor
canonique
$\RC_{\dot{\wb}} \Lambda_\theta \to \SC_{\dot{\wb}} \Lambda_\theta$
est dans la sous-cat\'egorie triangul\'ee de $D^b(\Lambda\Gb^F)$
engendr\'ee par les complexes $\RC_{\dot{\xb}} \Lambda_{\theta_\xb}$
o\`u $\xb \in \IC(\wb,\theta)$ et $\xb < \wb$.
\end{theo}

\begin{coro}\label{coro cone}
Soit $\wb\in\Sigma(\bar{S})$ et $\theta$ un caract\`ere lin\'eaire de
$\Tb^{\wb F}$. Alors, le c\^one du \mor
canonique
$\RC_{\dot{\wb}} \Lambda_\theta \to \SC_{\dot{\wb}} \Lambda_\theta$
est dans la sous-cat\'egorie triangul\'ee de $D^b(\Lambda\Gb^F)$
engendr\'ee par les complexes $\RC_{\dot{\xb}} \Lambda_{\theta'}$
o\`u $(\xb,\theta') \equivb{W} (\wb,\theta)$ et $\xb < \wb$.
\end{coro}

\begin{proof}
Le corollaire d\'ecoule du th\'eor\`eme \ref{theo cone} et du corollaire
\ref{iso quotient}.
\end{proof}

\begin{coro}\label{00000}
Si $(\wb,\theta)$ et $(\wb',\theta')$ appartiennent \`a deux s\'eries
rationnelles diff\'erentes de
$\nablab(\Tb,W,F)$, alors
$$\RRR\Hom^\bullet_{\Lambda\Gb^F}\bigl(\SC_{\dot{\wb}}\Lambda\Tb^{\wb F} e_\theta,
\SC_{\dot{\wb}'}\Lambda\Tb^{\wb'F} e_{\theta'}\bigr)=0.$$
\end{coro}

\begin{proof}
La proposition \ref{remark tresse} permet de
se ramener au cas o\`u $\wb,\wb'\in\Sigma(\bar{S})$.

On a un isomorphisme
$$\RRR\Hom^\bullet_{\Lambda\Gb^F}\bigl(\SC_{\dot{\wb}}\Lambda\Tb^{\wb F} e_\theta,
\SC_{\dot{\wb}'}\Lambda\Tb^{\wb'F} e_{\theta'}\bigr)\simeq
(\SC_{\dot{\wb}} \Lambda\Tb^{\wb F} e_{\theta})^* \otimes_{\Lambda\Gb^F}^\LLL
\SC_{\dot{\wb}'} \Lambda\Tb^{\wb'F}e_{\theta'}.$$
D'apr\`es (\ref{iso dual}), on a
$$(\SC_{\dot{\wb}} \Lambda\Tb^{\wb F} e_{\theta})^*\simeq
(\RC_{\dot{\wb}} \Lambda\Tb^{\wb F} e_{\theta^{-1}})^\opp[2 l(\wb)].$$
Le corollaire \ref{00000} r\'esulte alors du th\'eor\`eme \ref{00}
et du fait que $\SC_{\dot{\wb}'} \Lambda\Tb^{\wb' F} e_{\theta'}$ est dans
la cat\'egorie engendr\'ee par les $\RC_{\dot{\xb}} \Lambda\Tb^{\xb F}
e_{\theta_1}$, o\`u $(\xb,\theta_1)$ appartient \`a la m\^eme s\'erie
rationnelle que $(\wb',\theta')$ (corollaire \ref{coro cone}).
\end{proof}

On a un accouplement parfait
$$\fonction{<,>}{\K0(\Lambda\Gb^F) \times
\K0(\Lambda\Gb^F\proj)}{\ZM}{([M],[P])}{\dim(\Hom_{\Lambda\Gb^F}(P,M)).}$$

D'apr\`es \cite[proposition 7.5]{delu}, on a
\equat
|W|\cdot [\Lambda\Gb^F]=\sum_{w \in W}
(-1)^{l(w)}\frac{|\Gb^F|_{p'}}{|\Tb^{w F}|}[\RC_{(\dot{w})} \Lambda\Tb^{w F}].
\endequat
En appliquant (\ref{iso dual}), on obtient
\equat
|W| \cdot [\Lambda\Gb^F]=\sum_{w \in W}
(-1)^{l(w)} \frac{|\Gb^F|_{p'}}{|\Tb^{w F}|}[\SC_{(\dot{w})} \Lambda\Tb^{w F}].
\endequat
Par cons\'equent, si $M$ est un $\Lambda\Gb^F$-module simple, il existe
$((\dot{w}),\theta) \in \nablab(\Tb,W,F)$ tel que
$$<[M],[\SC_{(\dot{w})} \Lambda\Tb^{w F}e_\theta] > \not= 0.$$

On a donc~:

\begin{lem}\label{truc}
Si $M$ est un $\Lambda\Gb^F$-module simple, il existe
$(\wb,\theta) \in \nablab(\Tb,W,F)$ tel que
$$\RRR\Hom^\bullet_{\Lambda\Gb^F}(\SC_{\dot{\wb}} \Lambda\Tb^{\wb F}e_\theta,M) \not=0.$$
\end{lem}

\begin{rem}
Il serait souhaitable de trouver une preuve plus directe de ce lemme.
\end{rem}

\begin{prop}\label{coro crucial}
Si $M$ est un $\Lambda\Gb^F$-module simple et si
$\wb$ est une suite d'\eles de $W$ avec $l(\wb)$ minimal telle que
$$\RRR\Hom_{\Lambda\Gb^F}^\bullet(\SC_{\dot{\wb}} \Lambda\Tb^{\wb F},M) \not=0,$$
alors la cohomologie du complexe
$\RRR\Hom^\bullet_{\Lambda\Gb^F}(\SC_{\dot{\wb}} \Lambda\Tb^{\wb F},M)$ est
concentr\'ee en degr\'e $-l(\wb)$.
\end{prop}

\begin{proof}
La proposition \ref{remark tresse} permet de 
se ramener au cas o\`u $\wb\in\Sigma(\bar{S})$.
D'apr\`es (\ref{iso dual}), on a
$$\RRR\Hom_{\Lambda\Gb^F}(\SC_{\dot{\wb}} \Lambda\Tb^{\wb F},M)\simeq (\RC_{\dot{\wb}}
\Lambda\Tb^{\wb F})^\opp \otimes_{\Lambda\Gb^F}^\LLL M [2l(\wb)].$$
Cela montre que, si $\vb < \wb$, alors
$$(\RC_{\dot{\vb}} \Lambda\Tb^{\vb F})^\opp \otimes_{\Lambda\Gb^F}^\LLL M = 0.$$
Donc, pour tout \car lin\'eaire $\theta : \Tb^{\vb F} \to \Lambda^\times$, on a
$$(\RC_{\dot{\vb}} \Lambda_\theta)^\opp \otimes_{\Lambda\Gb^F}^\LLL M = 0.$$
Mais, d'apr\`es le th\'eor\`eme \ref{theo cone}, le c\^one du \mor
canonique $\RC_{\dot{\wb}} \Lambda\Tb^{\wb F} \to \SC_{\dot{\wb}} \Lambda\Tb^{\wb F}$
est dans la sous-cat\'egorie engendr\'ee par les $\RC_{\dot{\vb}} \Lambda_\theta$
o\`u $\vb < \wb$ et $\theta$ est un \car lin\'eaire de $\Tb^{\vb
F}$.
Par suite, le morphisme
$$(\RC_{\dot{\wb}} \Lambda\Tb^{\wb F})^\opp \otimes_{\Lambda\Gb^F}^\LLL M \longto
(\SC_{\dot{\wb}} \Lambda\Tb^{\wb F})^\opp \otimes_{\Lambda\Gb^F}^\LLL M$$
est un isomorphisme. La preuve de la proposition \ref{coro crucial} est
alors compl\'et\'ee par le lemme \ref{quasi lemme} et la proposition
\ref{goodguy}.
\end{proof}

\begin{rem}
\label{cabanes}
Comme nous l'a fait remarquer Marc Cabanes, on peut en fait d\'emontrer
directement que 
$\RRR\Hom^\bullet_{\Lambda\Gb^F}(\SC_{\dot{\wb}} \Lambda\Tb^{\wb F},
(M^*)^\opp)\not=0$.
En effet,
$$\RRR\Hom^\bullet_{\Lambda\Gb^F}(\SC_{\dot{\wb}} \Lambda\Tb^{\wb F},(M^*)^\opp)
\simeq
\RRR\Hom^\bullet_{\Lambda\Gb^F}((\SC_{\dot{\wb}} \Lambda\Tb^{\wb F})^*)^\opp,
M)$$
car $\SC_{\dot{\wb}} \Lambda\Tb^{\wb F}$ est parfait.
Maintenant, la dualit\'e de Poincar\'e fournit
$$\RRR\Hom^\bullet_{\Lambda\Gb^F}((\SC_{\dot{\wb}} \Lambda\Tb^{\wb F})^*)^\opp,
M)\simeq
\RRR\Hom^\bullet_{\Lambda\Gb^F}(\RC_{\dot{\wb}} \Lambda\Tb^{\wb F}),M)$$
d'o\`u le r\'esultat annonc\'e.
Par cons\'equent, l'application du lemme \ref{quasi lemme} se fait
dans un cas o\`u $M$ et $(M^*)^\opp$ v\'erifient les hypoth\`eses.
Le r\'esultat d'Haastert (cf remarque \ref{haastert}) suffit alors.
\end{rem}

La proposition suivante fournit plusieurs caract\'erisations 
d'un \'el\'ement $\wb$ de longueur minimale~:

\begin{prop}\label{equivalence}
Soit $M$ un $\Lambda\Gb^F$-module simple et soit $(\wb,\theta) \in \nablab(\Tb,W,F)$.
Les propri\'et\'es suivantes sont \'equivalentes~:

\tete{a} $(\wb,\theta)$ est de longueur minimale tel que
$<[M],[\RC_{\dot{\wb}} \Lambda\Tb^{\wb F}e_\theta] > \not= 0$~;

\tete{b} $(\wb,\theta)$ est de longueur minimale tel que
$\RRR\Hom^\bullet_{\Lambda\Gb^F}(\RC_{\dot{\wb}} \Lambda\Tb^{\wb F}e_\theta,M) \not= 0$~;

\tete{c} $(\wb,\theta)$ est de longueur minimale tel que
$\RRR\Hom^\bullet_{\Lambda\Gb^F}(M,\RC_{\dot{\wb}} \Lambda\Tb^{\wb F}e_\theta) \not= 0$~;

\tete{d} $(\wb,\theta)$ est de longueur minimale tel que
$\RRR\Hom^\bullet_{\Lambda\Gb^F}(\SC_{\dot{\wb}} \Lambda\Tb^{\wb F}e_\theta,M) \not= 0$~;

\tete{e} $(\wb,\theta)$ est de longueur minimale tel que
$\RRR\Hom^\bullet_{\Lambda\Gb^F}(M,\SC_{\dot{\wb}} \Lambda\Tb^{\wb F}e_\theta) \not= 0$~;

\tete{f} Un \'enonc\'e parmi $(a)$, $(b)$, $(c)$, $(d)$ ou $(e)$ en
rempla\c{c}ant $e_\theta$ par $e_{\theta^{-1}}$ et $M$ par $(M^*)^\opp$.
\end{prop}

\begin{proof}
D'apr\`es le corollaire \ref{coro cone}, le c\^one du morphisme canonique
$\RC_{\dot{\wb}} \Lambda\Tb^{\wb F}e_\theta\to
\SC_{\dot{\wb}} \Lambda\Tb^{\wb F}e_\theta$
est dans la sous-cat\'egorie engendr\'ee par les
$\RC_{\dot{\wb}} \Lambda\Tb^{\wb' F}$ avec $\wb'<\wb$. 
Ceci d\'emontre l'\'equivalence de (b) et (d) et de (c) et (e).

\medskip
Soit $C$ un complexe born\'e de $\Lambda \Gb^F$-modules projectifs, sans
facteur direct non nul homotope \`a $0$. Alors,
$\RRR\Hom^\bullet(C,M)\not=0$ si et seulement si 
$\RRR\Hom^\bullet(M,C)\not=0$, si et seulement si 
une enveloppe projective de $M$ est facteur direct d'une des composantes
de $C$. Cette remarque d\'emontre l'\'equivalence entre (b) et (c) et
entre (d) et (e).

\medskip
L'\'equivalence entre (a) et (d) r\'esulte de la proposition 
\ref{coro crucial}. L'\'equivalence avec les \'enonc\'es de (f)
r\'esulte de (\ref{iso dual}).
\end{proof}

\begin{rem}
\label{minimal}
Il n'est pas difficile de voir que si $\wb$ v\'erifie les propri\'et\'es
de la proposition, alors
$\mathbf{\sigma}(\wb)$ est dans $\sigma(W)$ et son image dans
$W$ est de longueur minimale dans sa classe de $F$-conjugaison (car
$[\RC_{\dot{\wb}} \Lambda\Tb^{\wb F}e_\theta]$ ne d\'epend
que de la classe de $F$-conjugaison de l'image de $\wb$ dans $W$).
Il est sans doute vrai que l'image de $\wb$ dans $W$ est dans une $F$-classe
de conjugaison de $W$ ne d\'ependant que de $M$.
C'est un r\'esultat connu lorsque $\Lambda=K$ \cite{dimirou}.
La preuve, qui requiert l'\'etude des valeurs propres de l'endomorphisme de
Frobenius sur le cohomologie, ne semble pas s'\'etendre \`a $\Lambda=k$.
\end{rem}

\section{Engendrement de $\Lambda\Gb^F\parf$}
\label{sec engendrement}

\subsection{Engendrement de cat\'egories de complexes parfaits.}
Soit $A$ une alg\`ebre libre de type fini
sur un anneau local complet $\OC$,
d'id\'eal maximal $\mG$ et de corps r\'esiduel ${\bar{\OC}}$.

Soit $E$ un ensemble fini de complexes parfaits de $A$, muni d'une
relation de pr\'eordre partiel.

\begin{lem} \label{lemmeengendrement} 
Supposons que pour tout $A$-module simple $M$, il existe $C \in E$ tel que
$\RRR\Hom^\bullet_{D^b(A\otimes_\OC {\bar{\OC}})}(C\otimes_\OC^\LLL {\bar{\OC}},M) \not= 0$
et que si $C$ est minimal
avec cette propri\'et\'e, alors la cohomologie du complexe
$\RRR\Hom^\bullet_{D^b(A\otimes_\OC {\bar{\OC}})}(C\otimes_\OC^\LLL {\bar{\OC}},M)$
est concentr\'ee en degr\'e $0$.

Alors, $A\parf$ est engendr\'ee par $E$.
\end{lem}

\begin{proof}
Nous  noterons $\bar{C}=C\otimes_\OC^\LLL \bar{\OC}$ pour $C$ un complexe de $A$-modules.
Soit $C \in E$. On note
$\SC(\le C)$ (\resp $\SC(< C)$) l'ensemble des classes d'\iso
de $A$-modules simples $M$ tels qu'il existe $C' \in E$ tel que $C' \le C$ 
(\resp $C' < C$) et $\RRR\Hom^\bullet(\bar{C}',M)\not=0$.

Quitte \`a remplacer les \eles de $E$ par des complexes quasi-isomorphes,
on peut supposer, et on le fera,
que ce sont des complexes born\'es de modules projectifs
sans facteurs directs non nuls homotopes \`a $0$. Alors, pour $C\in E$
et pour $M$ simple, on a $\RRR\Hom^\bullet(\bar{C},M[i])\not=0$ si et seulement
si une enveloppe projective $P_M$ de $M$ est facteur direct de $C^{-i}$.
Lorsque ces conditions
\'equivalentes sont v\'erifi\'ees, on a $M\in\SC(\le C)$. Lorsqu'en 
plus $i\not=0$, alors $M\in\SC(<C)$.

Soit $\langle E \rangle$ la sous-cat\'egorie de $K^b(A\proj)$ engendr\'ee par
$E$.

Consid\'erons l'assertion suivante~:

\medskip

{\it 
$(\bigstar_C)$ Pour $M \in \SC(\le C)$, alors
$P_M$ est dans $\langle E \rangle$. }

\medskip

Notons que le lemme \ref{lemmeengendrement} d\'ecoule imm\'ediatement
de $(\bigstar_C)$ pour tout $C$.
Nous allons montrer $(\bigstar_C)$ par r\'ecurrence sur $C$.

\bigskip

Soit $C$ un complexe minimal. Alors, $\RRR\Hom^\bullet(\bar{C},M[i])=0$ pour
tout module simple $M$ et tout entier non nul $i$. Par cons\'equent,
$C$ n'a qu'un terme non nul, $C^0$, qui est un $A$-module projectif de type
fini. Pour $M$ un $A$-module simple, alors $\Hom(C^0,M)\not=0$ si et seulement
si $P_M$ est facteur direct de $C^0$. On en d\'eduit imm\'ediatement
$(\bigstar_C)$.

\medskip

Prenons maintenant $C\in E$ tel que $(\bigstar_{C'})$ est vraie
pour tout $C'<C$.
Les facteurs directs ind\'ecomposables de $C^i$ sont de la forme
$P_M$ avec $M\in\SC(<C)$ pour $i\not=0$, donc il sont dans $\langle E\rangle$.
On en d\'eduit donc que $C^0$ est dans $\langle E\rangle$. Par cons\'equent,
$P_M$ est dans $\langle E\rangle$ pour $M\in\SC(C)$, ce qui montre
$(\bigstar_C)$.
\end{proof}

On se donne maintenant une relation d'\'equivalence $\sim$ sur $E$.

\begin{prop}\label{engendrement}
On garde les hypoth\`eses du lemme $\ref{lemmeengendrement}$ et on suppose
en plus que $\RRR\Hom^\bullet(C\otimes_\OC^\LLL {\bar{\OC}},C'\otimes_\OC^\LLL 
{\bar{\OC}})=0$
pour tous $C,C'\in E$ avec $C\not\sim C'$.

Soit $\langle X\rangle$ la sous-cat\'egorie de $A\parf$ engendr\'ee par
$X\in E/\sim$. Alors, 
on a une d\'ecomposition 
$$A\parf=\bigoplus_{X\in E/\sim} \langle X\rangle.$$

Par cons\'equent, pour $X\in E/\sim$, il existe un idempotent central $e_X$ de
$A$ tel que $\langle X\rangle=Ae_X\parf$ et on a une d\'ecomposition de 
l'unit\'e en somme d'idempotents orthogonaux
$$1=\sideset{}{^\perp}\sum_{X\in E/\sim} e_X.$$
\end{prop}

\begin{proof}
Pour $C$ et $C'$ deux complexes parfaits de $A$-modules, le morphisme
canonique $\RRR\Hom^\bullet(C,C')\otimes^\LLL_\OC \bar{\OC}\isom
\RRR\Hom^\bullet(C\otimes^\LLL_\OC \bar{\OC},C'\otimes^\LLL_\OC \bar{\OC})$
est un isomorphisme. Par cons\'equent, la nullit\'e de
$\RRR\Hom^\bullet(C,C')$ \'equivaut \`a celle de
$\RRR\Hom^\bullet(C\otimes^\LLL_\OC \bar{\OC},C'\otimes^\LLL_\OC \bar{\OC})$.

D'apr\`es le lemme \ref{lemmeengendrement}, la r\'eunion des
$\langle X\rangle$, o\`u $X\in E/\sim$, engendre $A\parf$. Puisque
ces sous-cat\'egories sont deux \`a deux orthogonales, la cat\'egorie
qu'elles engendrent est leur somme directe.
\end{proof}

\subsection{D\'ecomposition en s\'eries\label{sub engendrement}}

Nous allons d\'emontrer que la cat\'egorie
$\Lambda\Gb^F\parf$ est engendr\'ee par les
complexes $(\RC_{\dot{\wb}} \Lambda\Tb^{\wb F})_{\wb \in \Sigma(W)}$ et obtenir une
d\'ecomposition de cette cat\'egorie param\'etr\'ee par les s\'eries rationnelles.
Pour cela, il nous suffira d'appliquer la proposition
\ref{engendrement}, la v\'erification des hypoth\`eses r\'esultant du travail
effectu\'e dans le \S\ref{sec fdlesr}.

\medskip

Si $\mathbf{\XC}$ est une s\'erie rationnelle de $\nablab(\Tb,W,F)$,
nous noterons
$\mathbf{\XC}_{\min}$\index{xmin@$\mathbf{\XC}_{\min}$}
l'ensemble des couples $(w,\theta)$ o\`u $w\in W$ et
$\theta$ est un caract\`ere lin\'eaire de $\Tb^{wF}$
v\'erifiant la
propri\'et\'e suivante~:
il existe
un $\Lambda\Gb^F$-module simple $M$ tel que $((\dot{w}),\theta)$ est un \ele de
longueur minimale de $\nablab(\Tb,W,F)$ v\'erifiant
$$\RRR\Hom^\bullet_{\Lambda\Gb^F}(\SC_{\dot{w}} \Lambda\Tb^{wF}e_\theta,M) \not= 0.$$

Il r\'esulte de la remarque \ref{minimal}
que si $(w,\theta) \in \mathbf{\XC}_{\min}$,
alors $w$ est de longueur minimale dans sa classe de $F$-conjugaison.

\bigskip

\noindent{\bf Th\'eor\`eme A.} {\it 
Soit $\mathbf{\XC}$ une s\'erie rationnelle dans
$\nablab(\Tb,W,F)$. Alors, il existe un idempotent central
$e_\mathbf{\XC}=e_\mathbf{\XC}^{\Gb^F}$\index{ex@$e_\mathbf{\XC}$} de $\Lambda\Gb^F$ tel que
chacune des familles suivantes engendre la cat\'egorie
$\Lambda\Gb^Fe_\mathbf{\XC}\parf$~:
\begin{itemize}
\item $\{\RC_{\dot{\wb}} \Lambda\Tb^{\wb F} e_\theta\}_{(\wb,\theta) 
\in \mathbf{\XC}}$
\item $\{\SC_{\dot{\wb}} \Lambda\Tb^{\wb F} e_\theta\}_{(\wb,\theta) \in 
\mathbf{\XC}}$
\item $\{\RC_{(\dot{w})} \Lambda\Tb^{wF} e_\theta\}_{(w,\theta) \in 
\mathbf{\XC}_{\min}}$
\item $\{\SC_{(\dot{w})} \Lambda\Tb^{wF} e_\theta\}_{(w,\theta) \in 
\mathbf{\XC}_{\min}}$
\end{itemize}

On a une d\'ecomposition de $1$
en somme d'idempotents centraux deux \`a deux orthogonaux
$$1=\sideset{}{^\perp}\sum_{\mathbf{\XC} \in \nablab(\Tb,W,F)/\equivb{W}} e_\mathbf{\XC}.$$}

\begin{proof}
Soit $E$ l'ensemble des complexes
$\SC_{(\dot{w})} \Lambda\Tb^{wF} e_\theta[l(w)]$,
o\`u $(w,\theta)\in\mathbf{\XC}_{min}$ et $\mathbf{\XC}$ d\'ecrit les s\'eries
rationnelles.
On pose $(w,\theta)\le (w',\theta')$ si $w\le w'$.

D'apr\`es le lemme \ref{truc} et la proposition \ref{equivalence},
pour tout $\Lambda G^F$-module
simple $M$, il existe $C\in E$ tel que $\RRR\Hom^\bullet(C,M)\not=0$.
D'apr\`es la proposition \ref{coro crucial},
si $C$ est de longueur minimale avec
cette propri\'et\'e, alors $\RRR\Hom^\bullet(C,M)$ est concentr\'e
en degr\'e $0$.

On introduit la relation d'\'equivalence sur $E$ donn\'ee
par $(w,\theta)\sim (w',\theta')$ si $((w),\theta)$ et $((w'),\theta')$
sont rationnellement conjugu\'es.
Alors, $\RRR\Hom^\bullet(C,C')=0$ si $C\not\sim C'$
d'apr\`es le corollaire \ref{00000}.

Il r\'esulte alors de la proposition \ref{engendrement} que
$\{\SC_{(\dot{w})} \Lambda\Tb^{wF} e_\theta\}_{(w,\theta) \in \mathbf{\XC}_{\min}}$
engendre $\Lambda\Gb^Fe_\mathbf{\XC}\parf$. On montre de m\^eme que
$\{\SC_{\dot{\wb}} \Lambda\Tb^{\wb F} e_\theta\}_{(\wb,\theta) \in \mathbf{\XC}}$ engendre
$\Lambda\Gb^Fe_\mathbf{\XC}\parf$.

L'isomorphisme (\ref{iso dual}) (dualit\'e de Poincar\'e) montre
les assertions concernant les complexes
$\RC_{(\dot{w})} \Lambda\Tb^{w F} e_\theta$ et 
$\RC_{\wb} \Lambda\Tb^{\wb F} e_\theta$.
\end{proof}

\begin{rem}
\label{notation habituelle}
Consid\'erons tout d'abord le cas o\`u $\Lambda=K$.
La cat\'egorie $K\Gb^F\Mod$ \'etant semi-simple, la cat\'egorie
$K\Gb^F\parf$ est \'equivalente \`a la cat\'egorie des 
$K\Gb^F$-modules gradu\'es.

Pour $\mathbf{\XC} \in \nablab_K(\Tb,W,F)/\equivb{W}$, on note
$\EC(\Gb^F,\mathbf{\XC})$\index{eg@$\EC(\Gb^F,\mathbf{\XC})$}
 l'ensemble des caract\`eres irr\'eductibles
$\chi$ de $\Gb^F$ tels qu'il existe $(\wb,\theta) \in \mathbf{\XC}$ avec
$<\chi, {_K\RRR}_{\dot{\wb}} \theta > \not= 0$.
Il r\'esulte du th\'eor\`eme $A$ que
$$\Irr_K \Gb^F = \coprod_{\mathbf{\XC} \in \nablab_K(\Tb,W,F)/\equivb{W}} 
\EC(\Gb^F,\mathbf{\XC})$$
et que
$$e_\mathbf{\XC}=\sum_{\chi \in \EC(\Gb^F,\mathbf{\XC})} e_\chi,$$
o\`u $e_\chi$ d\'esigne l'unique idempotent primitif central de $K\Gb^F$
tel que $\chi(e_\chi)\not=0$. On retrouve ainsi la d\'ecomposition
de \cite{delu}.

\medskip 

Le morphisme canonique $R \to k$ induit une bijection
$\nablab_R(\Tb,W,F) \isom \nablab_k(\Tb,W,F),\ \mathbf{\XC}\mapsto
\overline{\mathbf{\XC}}$, compatible avec les
s\'eries rationnelles. Pour $\mathbf{\XC} \in \nablab_R(\Tb,W,F)/\equivb{W}$,
l'image de $e_\mathbf{\XC}\in R\Gb^F$ dans $k\Gb^F$ est
$e_{\overline{\mathbf{\XC}}}$.

Par ailleurs, nous avons une application surjective $\nablab_K(\Tb,W,F)
\to \nablab_k(\Tb,W,F)$ qui est obtenue par passage au quotient $R \to k$ (en effet,
un \car lin\'eaire $\Tb^{\wb F} \to K^\times$ est \`a valeurs dans $R^\times$).
Elle
induit une application surjective $\nablab_K(\Tb,W,F) \to \nablab_R(\Tb,W,F)$
en composant avec la bijection d\'ecrite pr\'ec\'edemment, d'o\`u
finalement
une application surjective entre s\'eries rationnelles
$\nablab_K(\Tb,W,F)/\equivb{W} \to \nablab_R(\Tb,W,F)/\equivb{W}$, $\mathbf{\XC}
\mapsto \mathbf{\XC}_{\ell'}$. Cette notation rappelle que cette application envoie 
un couple$(\wb,\theta)$ sur le couple $(\wb,\theta_{\ell'})$, o\`u $\theta_{\ell'}$ 
d\'esigne la $\ell'$-partie du \car lin\'eaire $\theta$.
Fixons une s\'erie rationnelle
$\mathbf{\XC}_0 \in \nablab_R(\Tb,W,F)/\equivb{W}$ et posons
$$\EC_\ell(\Gb^F,\mathbf{\XC}_0)=
\coprod_{\mathbf{\XC}_{\ell'}=\mathbf{\XC}_0} \EC(\Gb^F,\mathbf{\XC}).$$
On a alors
$$e_{\mathbf{\XC}_0}=\sum_{\mathbf{\XC}_{\ell'}=\mathbf{\XC}_0} e_\mathbf{\XC} =
\sum_{\chi \in \EC_\ell(\Gb^F,\mathbf{\XC}_0)} e_\chi\in R\Gb^F.$$
Nous retrouvons la compatibilit\'e des s\'eries avec les blocs \'etablie
par Brou\'e et Michel \cite[th\'eor\`eme 2.2]{BM}.
\end{rem}

\section{Foncteurs de Lusztig}

\subsection{D\'efinition}
Fixons une partie $I$\index{i@$I$}
de $\Delta$.
Nous notons $W_I$\index{wi@$W_I$}
le sous-groupe de $W$
engendr\'e par la famille $(s_\alpha)_{\alpha\in I}$, $\Pb_I$\index{pi@$\Pb_I$}
le \para $\Bb W_I \Bb$
de $\Gb$, $\Vb_I$\index{vi@$\Vb_I$} le radical unipotent de $\Pb_I$ et
$\Lb_I$\index{li@$\Lb_I$} le compl\'ement
de Levi de $\Vb_I$ dans $\Pb_I$ contenant $\Tb$. Le groupe de Weyl
de $\Lb_I$ relatif \`a $\Tb$ est alors \'egal \`a $W_I$. 
Nous notons aussi $\Bb_I$ le sous-groupe de 
Borel de $\Lb_I$ \'egal \`a $\Bb \cap \Lb_I$ et 
$\Ub_I$ son radical unipotent. Avec ces notations, on a
$\lexp{F}{{\boldsymbol{?}}_I}={\boldsymbol{?}}_{\phi(I)}$ si
$\boldsymbol{?}$ 
d\'esigne l'une des lettres $\Bb$, $\Lb$, $\Pb$, $\Ub$, $\Vb$ ou $W$.

Fixons maintenant un \ele $v$\index{v@$v$} de $W$ tel que
$v\phi(I)=I$. Alors, $\dot{v}F$ normalise $\Lb_I$ et on pose 
$$\Yb_{I,v}=\{g\Vb_I \in \Gb/\Vb_I~|~g^{-1}\hspace{0.1cm}\lexp{F}{g}
\in \Vb_I \dot{v} \hspace{0.1cm}
\lexp{F}{\Vb_I}\}.\index{yiv@$\Yb_{I,v}$}$$
Le groupe fini $\Gb^F$ agit par multiplication \`a gauche
sur la vari\'et\'e $\Yb_{I,v}$. De m\^eme, le groupe $\Lb_I^{\dot{v} F}$
agit par multiplication \`a droite. Munie de ces deux actions, $\Yb_{I,v}$
est une $\Gb^F$-vari\'et\'e-$\Lb_I^{\dot{v} F}$ de dimension
$l(v)$. Elle est
r\'eguli\`ere en tant que vari\'et\'e-$\Lb_I^{\dot{v} F}$. On
d\'efinit alors un foncteur entre cat\'egories d\'eriv\'ees
$$\fonction{\RC_{I,v}}{D^b(\Lambda\Lb_I^{\dot{v}
F})}{D^b(\Lambda\Gb^F)}{C}{\RRR\Gamma_c(\Yb_{I,v})
\otimes^\LLL_{\Lambda\Lb_I^{\dot{v} F}} C.}\index{riv@$\RC_{I,v}$}$$
Avec les notations du \S\ref{sub foncteur}, on a $\RC_{I,v}=
\RC_{\Lb_I^{\dot{v} F}}^{\Gb^F}(\Yb_{I,v})$.
Posons maintenant
$$\Xb_{I,v} =\{g\Pb_I \in \Gb/\Pb_I~|~ g^{-1}\hspace{0.1cm}\lexp{F}{g}
\in \Pb_I \dot{v} \hspace{0.1cm}
\lexp{F}{\Pb_I}\}.\index{xiv@$\Xb_{I,v}$}$$
Alors, d'apr\`es la proposition \ref{existence quotient} appliqu\'ee
\`a $\Pb=\Pb_I$, $\Hb=\Lb=\Lb_I$, $n=\dot{v}$ et $\Hb'=1$,
l'application naturelle
$$\fonction{\pi_{I,v}}{\Yb_{I,v}}{\Xb_{I,v}}{g\Vb_I}{g\Pb_I}$$
induit un \iso de $\Gb^F$-vari\'et\'es
        $$\Yb_{I,v} /\Lb_I^{\dot{v} F} \isom \Xb_{I,v}.$$

Pour finir, nous posons
$$\overline{\Xb}_{I,v} =\{g\Pb_I \in \Gb/\Pb_I~|~ g^{-1}\hspace{0.1cm}\lexp{F}{g}
\in \overline{\Pb_I \dot{v} \hspace{0.1cm}\lexp{F}{\Pb_I}}\}
\index{xivover@$\overline{\Xb}_{I,v}$}$$
et nous notons $j_{I,v}$\index{jiv@$j_{I,v}$} l'immersion ouverte canonique
$$j_{I,v} : \Xb_{I,v} \longto \overline{\Xb}_{I,v}.$$
C'est un morphisme de $\Gb^F$-vari\'et\'es, o\`u $\Gb^F$ agit sur 
$\overline{\Xb}_{I,v}$ par multiplication \`a gauche. De plus, 
$\overline{\Xb}_{I,v}$ est une vari\'et\'e projective. 
Nous noterons $\FC_{I,v} : D^b(\Lambda\Lb_I^{\dot{v} F}) \to
D^b_\Lambda(\Xb_{I,v})$\index{fiv@$\FC_{I,v}$}
le foncteur not\'e $\FC_{\Lb_I^{\dot{v} F}}^{\Yb_{I,v}}$ dans le
\S\ref{sub foncteur}.

\subsection{Transitivit\'e et s\'eries rationnelles}
On fixe une d\'ecomposition r\'eduite $v=s_1\cdots s_r$ ($s_i\in S$ et
$r=l(v)$) et on pose $\vb=(s_1,\ldots,s_r)$.

Soit $\wb=(w_1,\ldots,w_n)$ une suite finie d'\eles de $W_I$.
On a des \iso de $\Gb^F$-vari\'et\'es-$\Tb^{\wb\vb F}$
(cf \cite[lemme 3]{lufini} et (\ref{concatenationY}))
\refstepcounter{theo}

\begin{align*}
\label{fibre}
\Yb_{I,v} \times_{\Lb_I^{\dot{v} F}} \Yb^{\Lb_I,\dot{v} F}(\dot{\wb})
& \isom & \Yb^{\Gb,F}(\dot{w}_1,\ldots,\dot{w}_{n-1},\dot{w}_n\dot{v}) &
 \isom & \Yb^{\Gb,F}(\dot{\wb} \dot{\vb}) 
 \tag{$\boldsymbol{\arabic{section}.\arabic{theo}}$} \\
(g\Vb_I),(h_1\Ub_I,\ldots,h_n\Ub_I) & \mapsto & (gh_1\Ub,\ldots,gh_n\Ub)
\\
& & (g_1\Ub,\ldots,g_n\Ub) & \leftarrow & (g_1\Ub,\ldots,g_{n+r}\Ub).
\end{align*}
Par cons\'equent, d'apr\`es (\ref{composition}), on a
\equat\label{equat transit}
\RC_{I,v} \circ \RC_{\dot{\wb}}^{\Lb_I,\dot{v} F} \simeq
 \RC_{\dot{\wb} \dot{\vb}}^{\Gb,F}.
\endequat

\bigskip

Si $(\wb,\theta) \in \nablab(\Tb,W_I,\dot{v}F)$,
alors $(\wb \vb,\theta) \in \nablab(\Tb,W,F)$. De plus, si
$(\wb,\theta) \equiv_{W_I,\vb F} (\wb',\theta')$, alors
$(\wb \vb,\theta) \equiv_{W,F} (\wb' \vb,\theta')$.
Par suite, si $\mathbf{\XC}$ est une s\'erie rationnelle dans
$\nablab(\Tb,W_I,v F)$, il existe une unique
s\'erie rationnelle $\mathbf{\XC}^\Gb$ contenant $(\wb \vb,\theta)$ pour tout
$(\wb,\theta) \in \mathbf{\XC}$. Le th\'eor\`eme suivant (classique pour
$\Lambda=K$) montre que les foncteurs
de Lusztig pr\'eservent les s\'eries rationnelles.

\begin{theo}\label{series preservees}
Soit $\mathbf{\XC}$ une s\'erie rationnelle dans
$\nablab(\Tb,W_I,\dot{v} F)$
et soit $M$ dans\\
$D^b(\Lambda\Lb_I^{\dot{v} F}e_\mathbf{\XC}^{\Lb_I^{\dot{v} F}})$.
Alors, $\RC_{I,v}(M)$ est dans $D^b(\Lambda \Gb^F e_{\mathbf{\XC}^\Gb}^{\Gb^F})$.
\end{theo}

\begin{proof}
Par la formule des coefficients universels,
on peut supposer, et nous le ferons, que $\Lambda$ est un corps.

Il suffit en outre de consid\'erer le cas o\`u
$M=\Lambda\Lb_I^{\dot{v} F}e_\mathbf{\XC}^{\Lb_I^{\dot{v} F}}$. Mais, d'apr\`es le
th\'eor\`eme A,
$\Lambda\Lb_I^{\dot{v} F}e_\mathbf{\XC}^{\Lb_I^{\dot{v} F}}$ est dans la sous-cat\'egorie
engendr\'ee par les complexes 
$\RC_{\dot{\wb}}^{\Lb_I,\dot{v} F} \Lambda\Tb^{\wb\vb F}e_\theta$
o\`u $(\wb,\theta) \in \mathbf{\XC}$.
Cela montre qu'il suffit de consid\'erer le cas o\`u
$M=\RC_{\dot{\wb}}^{\Lb_I,\dot{v} F} \Lambda\Tb^{\wb\vb F} e_\theta$~:
alors,
$\RC_{I,v}(M)\simeq\RC_{\dot{\wb} \dot{\vb}}^{\Gb,F}
\Lambda\Tb^{\wb\vb F}e_\theta$ d'apr\`es
(\ref{equat transit}). La conclusion provient alors du th\'eor\`eme A.
\end{proof}

\subsection{Equivalence de Morita}
Commen\c{c}ons par une d\'efinition.
Le couple $(\wb,\theta) \in \nablab(\Tb,W_I,\dot{v}F)$ est dit {\it $(\Gb,I)$-r\'egulier}
(\resp {\it $(\Gb,I)$-super r\'egulier}) s'il v\'erifie la propri\'et\'e suivante~:

\medskip

\tete{R} {\it \hspace{1cm} Si $\alpha\in \Phi$ est
telle que $\theta \circ N_{\wb \vb}(\alpha^\vee)=1$, alors $\alpha\in \Phi_I$.}

\medskip

\noindent(\resp

\medskip

\tete{SR} {\it \hspace{1cm} Si $x \in W$ est tel que
$\lexp{x}{(\theta \circ N_{\wb \vb})}=\theta \circ N_{\wb \vb}$,
alors $x \in W_I$.})

\begin{lem}{\label{reg super reg}} Un couple $(\Gb,I)$-super r\'egulier est
$(\Gb,I)$-r\'egulier.
\end{lem}

\begin{proof}
Soient $(\wb,\theta)$ un couple $(\Gb,I)$-super r\'egulier
et $\alpha\in \Phi$ tel que $\theta \circ  N_{\wb \vb}(\alpha^\vee)=1$.
Pour tout $\lambda \in Y(\Tb)$, on a $(s_\alpha- 1)(\lambda) \in \ZM\alpha^\vee$, donc
$\theta \circ N_{\wb \vb}((s_\alpha-1)(\lambda)) = 1$, ce qui montre que
$\lexp{s_\alpha}{(\theta \circ N_{\wb \vb})}=\theta \circ N_{\wb \vb}$.
D'apr\`es (SR), on a alors $s_\alpha\in W_I$, donc $\alpha \in \Phi_I$.
\end{proof}

\begin{rem}
\label{reg serie}
Soient $(\wb,\theta)$ et $(\wb',\theta')$ deux couples
de $\nablab(\Tb,W_I,\dot{v}F)$ qui sont g\'eom\'e\-triquement conjugu\'es.
Alors, le couple
$(\wb,\theta)$ est $(\Gb,I)$-r\'egulier (\resp $(\Gb,I)$-super r\'egulier)
\ssi le couple $(\wb',\theta')$ l'est.
En effet, par d\'efinition de la conjugaison g\'eom\'etrique,
$\theta \circ N_{\wb \vb}$ et $\theta' \circ N_{\wb' \vb}$ sont 
conjugu\'es sous $W_I$.
\end{rem}

Une s\'erie rationnelle $\mathbf{\XC}$ dans $\nablab(\Tb,W_I,\dot{v} F)$
est dite {\it $(\Gb,I)$-r\'eguli\`ere} (\resp {\it $(\Gb,I)$-super r\'eguli\`ere})
si un de ses repr\'esentants l'est.
Compte tenu de la remarque \ref{reg serie}, $\mathbf{\XC}$ est $(\Gb,I)$-r\'eguli\`ere
(\resp $(\Gb,I)$-super r\'eguli\`ere) \ssi tous ses
repr\'esentants le sont.

La propri\'et\'e des couples r\'eguliers qui va nous \^etre utile par la suite
est contenue dans le lemme suivant.

\begin{lem}\label{wtheta}
Soit $(\wb,\theta) \in \nablab(\Tb,W_I,\dot{v} F)$ un couple
$(\Gb,I)$-r\'egulier avec $\wb\in\Sigma(I\cup\{1\})$.
Alors, $(\wb \vb)_\theta=\wb_\theta \vb$.
\end{lem}

\begin{proof}
Notons tout d'abord que $(\wb \vb)_\theta$ est calcul\'e
dans $(\Gb,F)$, alors que $\wb_\theta$ l'est dans $(\Lb,\dot{v}F)$.
Remarquons ensuite que $(\wb \vb)_\theta \le \wb_\theta \vb$.
\'Ecrivons $\wb=(t_1,\dots,t_n)$ et soit $i$ un entier naturel tel que
$1 \le i \le r$. Il s'agit de d\'emontrer que
$\theta(N_{\wb \vb}(\beta_{\wb \vb,n+i}^\vee)) \not= 1$. D'apr\`es la
propri\'et\'e (R), il suffit de d\'emontrer que $\beta_{\wb \vb,n+i}^\vee
\not\in \Phi_I$. Notons $\alpha$ la racine $\alpha_{\wb
\vb,n+i}=\alpha_{\vb,i}$.
Alors, $\beta_{\wb \vb,n+i}^\vee=t_1\dots t_n s_1 \dots s_{i-1}(\alpha^\vee)$.
Mais, $t_1 \dots t_n \in W_I$, donc on est ramen\'e \`a prouver l'assertion
suivante~:
$$s_1 \dots s_{i-1}(\alpha^\vee) \not\in \Phi_I.$$

Puisque $(s_1,\dots,s_r)$ est une d\'ecomposition r\'eduite de $v$, alors
$s_1 \dots s_{i-1}(\alpha^\vee)$ est une coracine positive dont l'image par
$v^{-1}$ est n\'egative. Or, par hypoth\`ese sur $v$, $v^{-1}(I)=\phi(I) \incl \Delta$.
Donc $v^{-1}(\Phi_I \cap \Phi^+) \incl \Phi^+$. Cela compl\`ete la preuve
du lemme \ref{wtheta}.
\end{proof}

Nous arrivons maintenant au r\'esultat principal de ce travail~:

\begin{theo}\label{i 0}
Soit $\mathbf{\XC}$  
une s\'erie rationnelle {\it $(\Gb,I)$-r\'eguli\`ere} dans
$\nablab(\Tb,W_I,v F)$.
Alors, le morphisme canonique
$$(j_{I,v})_!\FC_{I,v} \Lambda\Lb_I^{\dot{v} F}e_\mathbf{\XC}^{\Lb_I^{\dot{v} F}}
\isom \RRR(j_{I,v})_* \FC_{I,v} 
\Lambda\Lb_I^{\dot{v} F}e_\mathbf{\XC}^{\Lb_I^{\dot{v} F}}$$
est un isomorphisme.
\end{theo}

\begin{proof}
Notons ici
$i : \overline{\Xb}_{I,v}-\Xb_{I,v} \longto \overline{\Xb}_{I,v}$
l'immersion ferm\'ee canonique. Il s'agit de d\'emontrer que
$$i^* \RRR(j_{I,v})_* \FC_{I,v}
 \left(\Lambda\Lb_I^{\dot{v} F} e_\mathbf{\XC}^{\Lb_I^{\dot{v} F}}\right) = 0.
\leqno{(1)}$$
La formule des coefficients universels montre qu'il suffit de prouver
$(1)$ lorsque $\Lambda$ est un corps, ce que nous supposerons.
D'apr\`es le th\'eor\`eme A, le $\Lambda\Lb_I^{\dot{v} F}$-module
$\Lambda\Lb_I^{\dot{v} F}e_\mathbf{\XC}^{\Lb_I^{\dot{v} F}}$ est dans la sous-cat\'egorie
de $\Lambda\Lb_I^{\dot{v} F}\parf$ engendr\'ee par les complexes
$\SC_{\dot{\wb}}^{\Lb_I,\dot{v} F} \Lambda\Tb^{\wb\vb F} e_\theta$
pour $(\wb,\theta) \in \mathbf{\XC}$. Par suite, il suffit de d\'emontrer que
$$i^* \RRR(j_{I,v})_* \FC_{I,v} (\SC_{\dot{\wb}}^{\Lb_I,\dot{v} F}
\Lambda\Tb^{\wb\vb F} e_\theta) = 0\leqno{(2)}$$
pour tout $(\wb,\theta) \in \mathbf{\XC}$ avec $\wb\in\Sigma(I\cup\{1\})$.

Fixons un \ele $(\wb,\theta) \in \mathbf{\XC}$ avec
$\wb\in\Sigma(I\cup\{1\})$. On a un diagramme commutatif
$$\diagram
\Yb_{I,v}\times \Yb^{\Lb_I,\dot{v} F}(\dot{\wb}) \ar[rr]^a \ar[d]_{b} &&
 \Yb_{I,v}\times_{\Lb_I^{\dot{v} F}} \Yb^{\Lb_I,\dot{v} F}(\dot{\wb}) \ar[d]^{b'} \\
\Yb_{I,v}\times \Xb^{\Lb_I,\dot{v} F}(\wb) \ar[rr]^{a'} \ar[dr]_c &&
 \Yb_{I,v}\times_{\Lb_I^{\dot{v} F}} \Xb^{\Lb_I,\dot{v} F}(\wb) \ar[dl]^{c'} \\
& \Xb_{I,v}
\enddiagram$$
o\`u les fl\`eches horizontales sont les passages au quotient par
$\Lb_I^{\dot{v} F}$, les fl\`eches verticales les passages au quotient
par $\Tb^{\wb\vb F}$ et les fl\`eches diagonales proviennent du
morphisme structurel $ \Xb^{\Lb_I,\dot{v} F}(\wb)\to \Spec\FM$.

On a des isomorphismes dans $D^b_{\Lambda \Tb^{\wb\vb F}}(\Xb_{I,v})$
\begin{align*}
\RRR(c'b')_*\Lambda&\simeq \left((\RRR c'_*\RRR(b'a)_*\Lambda\right)
\otimes_{\Lambda\Lb_I^{\dot{v}
F}}^\LLL\Lambda \hspace{2cm} \textrm{(lemme \ref{lemme foncteur}})\\
&\simeq \left(\RRR c_* \RRR b_*\Lambda\right)\otimes_{\Lambda\Lb_I^{\dot{v}
F}}^\LLL\Lambda,
\end{align*}
ce qui montre que
$$\RRR(c'b')_*\Lambda\simeq \FC_{I,v}(\Lambda \Lb_I^{\dot{v} F})
\otimes_{\Lambda\Lb_I^{\dot{v} F}}^\LLL 
\RRR\Gamma(\Yb^{\Lb_I,\dot{v} F}(\wb)).\leqno{(3)}$$

Via l'isomorphisme (\ref{fibre}), le morphisme $c'$ devient
$$\fonction{\tau}{\Xb^{\Gb,F}
(\wb \vb)}{\Xb_{I,v}}{(g_1\Bb,\dots,g_n\Bb)}{g_1\Pb_I}$$
et apr\`es application de 
$-\otimes^\LLL_{\Lambda \Tb^{\wb\vb F}}\Lambda\Tb^{\wb\vb F}e_\theta$,
l'isomorphisme (3) devient
$$\RRR\tau_* \FC^{\wb \vb}(\Lambda\Tb^{\wb\vb F}e_\theta) \simeq
\FC_{I,v} (\SC_{\dot{\wb}}^{\Lb_I,\dot{v} F} \Lambda\Tb^{\wb\vb F} e_\theta).$$

Consid\'erons maintenant le diagramme commutatif suivant~:
$$\xymatrix{
\Xb^{\Gb,F}(\wb \vb) \rrto^{j_{\wb \vb}^{\overline{\wb \vb}}}
\ddto_{\tau} &&
\Xb^{\Gb,F}(\overline{\wb \vb}) \ddto^{\bar{\tau}}\\
&&\\
\Xb_{I,v} \rrto^{j_{I,v}} &&\overline{\Xb}_{I,v},
}$$
o\`u $\bar{\tau} : \Xb^{\Gb,F}(\overline{\wb \vb}) \to \overline{\Xb}_{I,v}$, 
$(g_1\Bb,\dots,g_n\Bb) \mapsto g_1 \Pb_I$ est un morphisme projectif. 
Pour d\'emontrer (2), il suffit, en utilisant le th\'eor\`eme de changement de
base propre, de d\'emontrer que la restriction du complexe
$\RRR(j_{\wb \vb}^{\overline{\wb \vb}})_*
\FC^{\wb \vb}(\Lambda\Tb^{\wb\vb F}e_\theta)$ \`a
$\Xb^{\Gb,F}(\overline{\wb \vb}) - \bar{\tau}^{-1}(\Xb_{I,v})$ est nulle.
Notons $i_1 : \Xb^{\Gb,F}(\overline{\wb \vb}) - \bar{\tau}^{-1}(\Xb_{I,v})
\injto \Xb^{\Gb,F}(\overline{\wb \vb})$
l'immersion ferm\'ee canonique. Il nous suffit de d\'emontrer que
$$i_1^*\RRR(j_{\wb \vb}^{\overline{\wb \vb}})_*
\FC^{\wb \vb}_\theta=0.\leqno{(4)}$$
On a
$$\Xb^{\Gb,F}(\overline{\wb \vb}) - \bar{\tau}(\Xb_{I,v})
=\coprod_{\substack{\wb' \le \wb\\ \vb' < \vb}}
\Xb^{\Gb,F}(\wb' \vb').$$
Mais le fait que $\mathbf{\XC}$ soit $(\Gb,I)$-r\'eguli\`ere implique
que $(\wb \vb)_\theta$ est \'egal \`a $\wb_\theta \vb$
(lemme \ref{wtheta}).
Le r\'esultat d\'ecoule alors du th\'eor\`eme \ref{theo mono}.
\end{proof}

\begin{rem}
Le lecteur int\'eress\'e obtiendra sans difficult\'e une version
sans hypoth\`ese de r\'egularit\'e
du th\'eor\`eme \ref{i 0} qui g\'en\'eralise le corollaire
\ref{image directe 2}.
\end{rem}

\begin{coro}
\label{concentration}
Soit $\mathbf{\XC}$ une s\'erie $(\Gb,I)$-r\'eguli\`ere dans
$\nablab(\Tb,W_I,\dot{v} F)$.
Alors,\\
$H^i_c(\Xb_{I,v},\FC_{I,v} \Lambda\Lb_I^{\dot{v} F}e_\mathbf{\XC}^{\Lb_I^{\dot{v} F}})$
est un $\Lambda$-module libre, nul pour $i\not=l(v)$.
\end{coro}

\begin{proof}
Par la formule des coefficients universels, il suffit de d\'emontrer
l'annulation dans le cas o\`u $\Lambda$ est un corps.

Par application du foncteur $\RRR(\pib_{\overline{\Xb}_{I,v}})_*$
\`a l'isomorphisme du th\'eor\`eme \ref{i 0}, on obtient que le
morphisme canonique
$$\RRR(\pib_{\Xb_{I,v}})_!\FC_{I,v} 
\Lambda\Lb_I^{\dot{v} F}e_\mathbf{\XC}^{\Lb_I^{\dot{v} F}}
\isom \RRR(\pib_{\Xb_{I,v}})_*\FC_{I,v} 
\Lambda\Lb_I^{\dot{v} F}e_\mathbf{\XC}^{\Lb_I^{\dot{v} F}}$$
est un isomorphisme et le lemme \ref{quasi lemme} (joint au lemme
\ref{lemme foncteur}) fournit la conclusion.
\end{proof}

\begin{rem}
Si la s\'erie $\XC$ est r\'eguli\`ere, alors la s\'erie
$\XC^{-1}=\{(\wb,\theta^{-1}) | (\wb,\theta)\in\XC\}$ est encore r\'eguli\`ere.
On en d\'eduit alors, comme dans la remarque \ref{cabanes}, que
l'utilisation du lemme \ref{quasi lemme} se fait dans le cas
particulier o\`u $M$ et $(M^*)^\opp$ v\'erifient les hypoth\`eses.
\end{rem}

\smallskip
Comme l'a montr\'e Brou\'e \cite[pp. 61--62]{B},
le th\'eor\`eme \ref{i 0}, qui est
de nature g\'eom\'etrique, admet pour cons\'equence le th\'eor\`eme
alg\'ebrique suivant (nous reprenons la preuve de \cite[th\'eor\`eme 3.3]{B}).
Lorsque $\Lambda=K$ et $q$ est plus grand que le nombre de Coxeter de
$\Gb$, ce r\'esultat est d\^u \`a Lusztig \cite[Proposition 6.6]{lugf}.
Pour le cas g\'en\'eral o\`u $\Lambda=K$, seule la concentration en un
degr\'e (non d\'etermin\'e) de la cohomologie \'etait connue.

\bigskip

\noindent{\bf Th\'eor\`eme B.} {\it Soit $\mathbf{\XC}$
une s\'erie {\it $(\Gb,I)$-super r\'eguli\`ere}
dans $\nablab(\Tb,W_I,\dot{v} F)$.
Alors, $\RC_{I,v}$ induit une \'equivalence de Morita entre
les $\Lambda$-alg\`ebres $\Lambda\Lb_I^{\dot{v} F}e_\mathbf{\XC}^{\Lb_I^{\dot{v} F}}$ et
$\Lambda\Gb^Fe_{\mathbf{\XC}^\Gb}^{\Gb^F}$.

Plus pr\'ecis\'ement,
le complexe $\RC_{I,v} \Lambda\Lb_I^{\dot{v} F}e_\mathbf{\XC}^{\Lb_I^{\dot{v} F}}$ de
$(\Lambda\Gb^F,\Lambda\Lb_I^{\dot{v} F}e_\mathbf{\XC}^{\Lb_I^{\dot{v} F}})$-bimodules
n'a de la cohomologie qu'en degr\'e $r=l(v)$ et son $r$-i\`eme bimodule de cohomologie
induit l'\'equivalence de Morita d\'ecrite.}

\bigskip

\begin{proof}
Pour montrer l'\'equivalence de Morita, il suffit de consid\'erer le
cas o\`u $\Lambda=R$.
D'apr\`es le lemme \ref{reg super reg} et le corollaire \ref{concentration},
le bimodule
$M=H^r_c(\Xb_{I,v},\FC_{I,v} R\Lb_I^{\dot{v} F}e_\mathbf{\XC}^{\Lb_I^{\dot{v} F}})$
est le seul groupe de cohomologie non nul
du complexe parfait de $R\Gb^F$-modules $\RC_{I,v} R\Lb_I^{\dot{v}
F}e_\mathbf{\XC}^{\Lb_I^{\dot{v} F}}$, donc il est projectif comme $R\Gb^F$-module.
De m\^eme, il est projectif comme
$R\Lb_I^{\dot{v} F}e_\mathbf{\XC}^{\Lb_I^{\dot{v} F}}$-module \`a droite.

Soit $M_K=K\otimes_R M$.
D'apr\`es \cite[th\'eor\`eme 8]{lufini}, la $(\Gb,I)$-super
r\'egularit\'e de $\XC$ montre que
$$M_K^*\otimes_{K\Gb^F}M_K\simeq K\Lb_I^{\dot{v}
F}e_\mathbf{\XC}^{\Lb_I^{\dot{v} F}}.$$
Puisque $K\Gb^Fe_{\mathbf{\XC}^\Gb}^{\Gb^F}$ est facteur direct du
$(K\Gb^Fe_{\mathbf{\XC}^\Gb}^{\Gb^F},K\Gb^Fe_{\mathbf{\XC}^\Gb}^{\Gb^F})$-bimodule
$M_K\otimes_{K\Lb_I^{\dot{v} F}}M_K^*$,
on en d\'eduit que $M_K$ induit une \'equivalence de Morita entre
$K\Lb_I^{\dot{v} F}e_\mathbf{\XC}^{\Lb_I^{\dot{v} F}}$ et $K\Gb^Fe_{\mathbf{\XC}^\Gb}^{\Gb^F}$.
L'\'equivalence de Morita sur $R$
r\'esulte maintenant de \cite[th\'eor\`eme 2.4]{B}.
\end{proof}

\begin{rem}
Le th\'eor\`eme pr\'ec\'edent montre la co\"{\i}ncidence des matrices
de d\'ecompo\-sition. Cette propri\'et\'e avait \'et\'e conjectur\'ee par
Hi\ss\ \cite[p.342]{Hi}.
\end{rem}

\begin{rem}
D'apr\`es \cite[proposition 12]{lufini}, on a $\dim
M_K\otimes_{K\Lb_I^{\dot{v} F}}V=n\dim V$, pour tout 
$K\Lb_I^{\dot{v} F}e_\mathbf{\XC}^{\Lb_I^{\dot{v} F}}$-module simple
$V$, o\`u $n=\frac{|\Gb^F|_{p'}}{|\Lb_I^{\dot{v}F}|_{p'}}$.
On d\'eduit alors du th\'eor\`eme B que
$\Lambda\Gb^Fe_{\mathbf{\XC}^\Gb}^{\Gb^F}$ 
est isomorphe \`a une alg\`ebre de matrices de dimension $n^2$ sur
$\Lambda\Lb_I^{\dot{v} F}e_\mathbf{\XC}^{\Lb_I^{\dot{v} F}}$.
\end{rem}

\section{Changement de point de vue}

Nous avons pris le parti, pour la commodit\'e des preuves, de travailler
syst\'ematiquement avec des sous-groupes de Levi standards, quitte \`a
modifier l'isog\'enie $F$. Le but de cette partie est de traduire
les r\'esultats principaux qui pr\'ec\`edent en ne s'int\'eressant
qu'aux sous-groupes de Levi et aux tores maximaux $F$-stables. Ceci
permet d'all\'eger les notations et certains lecteurs y retrouveront
un cadre plus familier. 

Nous nous bornerons \`a \'enoncer les analogues des th\'eor\`emes A, 
\ref{series preservees}, \ref{i 0} et B.

\subsection{Notations} 
Soit $\Pb$ un \para de $\Gb$ de radical unipotent $\Vb$, avec un
compl\'ement de Levi $F$-stable $\Lb$.
Nous posons
$$\Yb_\Vb^\Gb=\{g\Vb \in \Gb/\Vb~|~g^{-1}F(g) \in \Vb\cdot\lexp{F}{\Vb}\}
\index{yvg@$\Yb_\Vb^\Gb$}$$
$$\Xb_\Pb^\Gb=\{g\Pb \in \Gb/\Pb~|~g^{-1}F(g) \in \Pb\cdot\lexp{F}{\Pb}\}.
\leqno{\mathrm{et}}
\index{xpg@$\Xb_\Pb^\Gb$}$$
Alors,
$\Yb_\Vb^\Gb$ est une $\Gb^F$-vari\'et\'e-$\Lb^F$ (pour les actions par
multiplication
\`a gauche et \`a droite), $\Xb_\Pb^\Gb$ est une $\Gb^F$-vari\'et\'e et le
morphisme \'etale $\pi_\Vb : \Yb_\Vb^\Gb \to \Xb_\Pb^\Gb$ induit un \iso de
vari\'et\'es $\Yb_\Vb^\Gb/\Lb^F \simeq \Xb_\Pb^\Gb$ (cf (\ref{YversX})
et \S\ref{traduction}). 
Nous d\'efinissons alors les foncteurs
$$\RC_{\Lb \incl \Pb}^\Gb : D^b(\Lambda\Lb^F) \longto D^b(\Lambda\Gb^F)
\index{rlg@$\RC_{\Lb \incl \Pb}^\Gb$}$$
$$\SC_{\Lb \incl \Pb}^\Gb : D^b(\Lambda\Lb^F) \longto D^b(\Lambda\Gb^F)
\leqno{\mathrm{et}}
\index{slg@$\SC_{\Lb \incl \Pb}^\Gb$}$$
comme au \S \ref{sub foncteur}. 

Pour finir, posons
$$\overline{\Xb}_\Pb^\Gb =\{g\Pb \in \Gb/\Pb~|~g^{-1}F(g) \in
\overline{\Pb\cdot \lexp{F}{\Pb}}\}
\index{xg@$\overline{\Xb}_\Pb^\Gb$}$$
et notons $j_\Pb^\Gb : \Xb_\Pb^\Gb \injto
\overline{\Xb}_\Pb^\Gb$\index{jpg@$j_\Pb^\Gb$} l'immersion
ouverte canonique. 

\subsection{Rapport avec les groupes $\Lb_I^{\dot{v} F}$\label{traduction}}
Nous expliquons ici le lien entre les vari\'et\'es $\Yb_\Vb^\Gb$ et les 
vari\'et\'es $\Yb_{I,v}$ d\'efinies pr\'ec\'edemment. 
Soit $I$ l'unique partie de $\Delta$ telle que 
$\Pb$ soit conjugu\'e \`a $\Pb_I$. Il existe alors un \ele $x_1 \in \Gb$ 
tel que $\lexp{x_1}{(\Lb_I,\Pb_I)}=(\Lb,\Pb)$. Posons $v_1=x_1^{-1}F(x_1)$. 
Puisque $\Lb$ est $F$-stable, on a $\lexp{v_1F}{\Lb_I}=\Lb_I$. Par 
cons\'equent, il existe $a \in \Lb_I$ tel que $\lexp{av_1F}{(\Tb,\Bb_I)}=(\Tb,\Bb_I)$. 
Notons $v$ l'image de $av_1 \in N_\Gb(\Tb)$ dans $W$. Nous choisissons 
$a$ de sorte que $av_1=\dot{v}$.
Par le th\'eor\`eme de Lang, on trouve $b \in \Lb_I$ tel que 
$b^{-1}~\lexp{v_1F}{b}=a$. 

Posons maintenant $x=x_1b$. Alors $x^{-1}F(x)=b^{-1}v_1F(b)=av_1=\dot{v}$. 
De plus $\dot{v} F$ stabilise $\Bb_I$, donc $v\phi(I)=I$. 
D'autre part, $\lexp{x}{(\Lb_I,\Pb_I)}=(\Lb,\Pb)$ et l'application 
$$\fonctio{\Lb_I}{\Lb}{g}{xgx^{-1}}$$
induit un isomorphisme $\Lb_I^{\dot{v}F}\isom \Lb^F$, que nous noterons aussi 
par la lettre $x$. De plus, $\lexp{x}{\Vb_I}=\Vb$. 

Les applications 
$$\fonctio{\Yb_\Vb^\Gb}{\Yb_{I,v}}{g\Vb}{gx\Vb_I,}$$
$$\fonctio{\Xb_\Pb^\Gb}{\Xb_{I,v}}{g\Pb}{gx\Pb_I}$$
$$\fonctio{\overline{\Xb}_\Pb^\Gb}{\overline{\Xb}_{I,v}}{g\Pb}{gx\Pb_I}
\leqno{\mathrm{et}}$$
sont des isomorphismes, de $\Gb^F$-vari\'et\'es-$\Lb^F$ pour le premier, 
de $\Gb^F$-vari\'et\'es pour le second et le troisi\`eme. Ici, $\Yb_{I,v}$ est vue comme 
une vari\'et\'e-$\Lb^F$ via l'isomorphisme $x^{-1}$. 
De plus, le diagramme correspondant
$$\diagram
\Yb_\Vb^\Gb \rrto^{\DS{\sim}} \ddto_{\DS{\pi_\Vb}} && \Yb_{I,v} \ddto^{\DS{\pi_{I,v}}} \\
&& \\
\Xb_\Pb^\Gb \rrto^{\DS{\sim}} \ddto|<\ahook_{\DS{j_\Pb^\Gb}} && 
\Xb_{I,v} \ddto|<\ahook^{\DS{j_{I,v}}} \\
&& \\
\overline{\Xb}_\Pb^\Gb \rrto^{\DS{\sim}} && \overline{\Xb}_{I,v} 
\enddiagram$$
est commutatif.

Ces isomorphismes permettent de faire le lien entre tous les objets 
d\'efinis pr\'ec\'edemment. Par exemple, on a un diagramme commutatif
\equat\label{commutatif RLG}
\diagram
D^b(\Lambda\Lb^F) \drto_{\DS{\RC_{\Lb \incl \Pb}^\Gb}}
 \rrto^{\DS{x_*}}_\sim && 
D^b(\Lambda\Lb_I^{\dot{v} F}) \dlto^{\DS{\RC_{I,v}}} \\
&D^b(\Lambda\Gb^F).& \\
\enddiagram
\endequat
qui montre que l'on peut jongler entre les divers foncteurs de Lusztig.

\begin{rem}
Il est facile de montrer que, si $(I',v')$ est un autre couple v\'erifiant les 
propri\'et\'es suivantes~:

\tete{1} $I' \incl \Delta$, $v' \in W$,

\tete{2} $v'\phi(I')=I'$,

\tete{3} il existe $x' \in \Gb$ tel que $\lexp{x'}{(\Lb_{I'},\Pb_{I'})}=(\Lb,\Pb)$ 
et $x^{\prime -1}~\lexp{F}{x'}=\dot{v}'$, 

\noindent alors $(I,v)=(I',v')$. On dira que le couple $(I,v)$ est {\it associ\'e} 
au couple $(\Lb,\Pb)$. 
\end{rem}

\subsection{S\'eries rationnelles} 
Notons $\nablab(\Gb,F)$\index{nablagf@$\nablab(\Gb,F)$}
l'ensemble des triplets $(\Tb',\Bb',\theta')$,
o\`u $\Tb'$ est un \tor $F$-stable de $\Gb$, $\Bb'$ est un sous-groupe de
Borel de $\Gb$
contenant $\Tb'$ et $\theta' : \Tb^{\prime F} \to \Lambda^\times$ un \car lin\'eaire
d'ordre inversible dans $\Lambda$. Le choix de $\imath$ et $\jmath$ fournit
une partition de $\nablab(\Gb,F)$ en {\it s\'eries rationnelles}, {\em
i.e.}, selon 
les classes de $\Gb^{*F^*}$-conjugaison d'\eles
de $\Gb^{*F^*}_{\sem,\Lambda}$ (on note
$\Gb^{*F^*}_{\sem,\Lambda}$\index{gstar@$\Gb^{*F^*}_{\sem,\Lambda}$} l'ensemble des \eles
semi-simples de $\Gb^{*F^*}$ d'ordre inversible dans $\Lambda$).
Nous noterons
$\nablab(\Gb,F,(s)_{\Gb^{*F^*}})$\index{nablagfs@$\nablab(\Gb,F,(s)_{\Gb^{*F^*}})$}
la s\'erie rationnelle associ\'ee \`a $s\in \Gb^{*F^*}_{\sem,\Lambda}$. 

Si $(\Tb',\Bb',\theta') \in \nablab(\Lb,F)$, alors il existe 
un unique \'el\'ement $w \in W_I$ tel que $(\vide,wv)$ soit 
associ\'e \`a $(\Tb',\Bb'\Vb)$. Il existe un \ele $y \in \Gb$ tel que 
$\lexp{y}{(\Tb,\Bb)}=(\Tb',\Bb'\Vb)$ et $y^{-1}~\lexp{F}{y}=\dot{w}\dot{v}$ et, comme 
dans le paragraphe pr\'ec\'edent, $y$ induit un isomorphisme 
$\Tb^{\prime F} \isom \Tb^{wvF}$. On note $\theta$ le caract\`ere lin\'eaire 
de $\Tb^{wvF}$ correspondant \`a $\theta'$ via $y$. Alors,
$((w),\theta)$ est dans $\nablab(\Tb,W_I,\dot{v} F)$ et il est d\'efini de mani\`ere
unique par le triplet $(\Tb',\Bb',\theta')$. 


Nous avons donc construit une application
$$\psi_{\Lb \incl \Pb} : \nablab(\Lb,F) \longto
\nablab(\Tb,W_I,vF).$$  
Lorsque $\Lb=\Gb$ (auquel cas $I=\Delta$ et $v=1$), nous la noterons simplement $\psi$.   
Cette application conserve les s\'eries rationnelles. Mieux, si 
$\XC$ est une s\'erie rationnelle dans $\nablab(\Tb,W,F)$ et si 
$((\dot{w}),\theta) \in \XC$, alors il existe $(\Tb',\Bb',\theta')\in
\nablab(\Gb,F)$
tel que $\psi(\Tb',\Bb',\theta')=((\dot{w}),\theta)$. On a alors
$$\RC_{\Tb' \incl \Bb'}^{\Gb'} \Lambda\Tb^{\prime
F}e_{\theta'}\isom\RC_{\dot{w}}\Lambda\Tb^{wF}e_{\theta}.$$
Cette discussion montre que le th\'eor\`eme A devient~:

\bigskip

\noindent{\bf Th\'eor\`eme A'.} {\it
Pour $(s)_{\Gb^{*F^*}} \in \Gb_{\sem,\Lambda}^{*F^*}/\sim$, il existe
un idempotent central $e_s=e_s^{\Gb^F}=e_{\Lambda,s}^{\Gb^F}$ de $\Lambda\Gb^F$ 
tel que la famille
$(\RC_{\Tb' \incl \Bb'}^\Gb \Lambda\Tb^{\prime F} e_\theta)_{(\Tb',\Bb',\theta) \in
\nablab(\Gb,F,(s)_{\Gb^{*F^*}})}$ (\resp\\
$(\SC_{\Tb' \incl \Bb'}^\Gb \Lambda\Tb^{\prime F} e_\theta)_{(\Tb',\Bb',\theta) \in
\nablab(\Gb,F,(s)_{\Gb^{*F^*}})}$) engendre $\Lambda\Gb^Fe_s\parf$.

On a une d\'ecomposition de $1$ en somme d'idempotents centraux deux \`a
deux orthogonaux
$$1 = \sideset{}{^\perp}\sum_{(s)_{\Gb^{*F^*}}\in \Gb_{\sem,\Lambda}^{*F^*}/\sim}
e_s.$$}

\bigskip

\begin{rem}
\label{bofbof}
On a
$\Gb_{\sem,R}^{*F^*}=\Gb_{\sem,k}^{*F^*}$. Pour
$s \in \Gb_{\sem,R}^{*F^*}$, la
remarque \ref{notation habituelle} montre que
$$e_{R,s}=\sum_{\substack{(t)_{\Gb^{*F^*}} \in \Gb^{*F^*}_\sem/\sim\\
(t_{\ell'})_{\Gb^{*F^*}}=(s)_{\Gb^{*F^*}}}} e_{K,t}.$$
De plus, si $(t)_{\Gb^{*F^*}} \in \Gb^{*F^*}_\sem/\sim$, alors
$$e_{K,t}=\sum_{\chi \in \EC(\Gb^F,(t)_{\Gb^{*F^*}})} e_\chi,$$
o\`u $\EC(\Gb^F,(t)_{\Gb^{*F^*}})$ est la s\'erie de Lusztig rationnelle
associ\'ee \`a $(t)_{\Gb^{*F^*}}$.
\end{rem}

Le th\'eor\`eme \ref{series preservees} devient quant \`a lui~:

\begin{theo}\label{series preservees '}
Soit $s \in \Lb_{\sem,\Lambda}^{*F^*}$ et soit $M \in D^b(\Lambda\Lb^Fe_s^{\Lb^F})$.
Alors, $\RC_{\Lb \incl \Pb}^\Gb(M) \in D^b(\Lambda\Gb^Fe_s^{\Gb^F})$.
\end{theo}

\subsection{Equivalence de Morita}
Pour finir la traduction et obtenir des r\'esultats \'equivalents aux
th\'eor\`emes \ref{i 0} et B, nous devons voir ce que deviennent 
les d\'efinitions de s\'eries r\'eguli\`eres et super-r\'eguli\`eres 
\`a travers l'application $\psi$ d\'efinie au paragraphe pr\'ec\'edent.

\smallskip
Fixons un \'el\'ement semi-simple $s \in \Gb^{*F^*}_{\sem,\Lambda}$
et un triplet $(\Tb',\Bb',\theta') \in \nablab(\Gb,F,(s)_{\Gb^{*F^*}})$. 
Notons $\Tb^{\prime *}$ un tore maximal $F^*$-stable de $\Gb^*$ dual de $\Tb'$. 

Puisque $(\Tb',\Bb',\theta') \in \nablab(\Gb,F,(s)_{\Gb^{*F^*}})$,
on peut supposer que $s \in \Tb^{\prime *F^*}$. Alors il existe 
$\lambda \in Y(\Tb^{\prime *})=X(\Tb')$ tel que $s=N_{F^{*d}/F^*}(\lambda(\zeta))$. 
On a alors, pour tout $\mu \in Y(\Tb')=X(\Tb^{\prime *})$, 
\equat\label{s}
\theta'(N_{F^d/F}(\mu(\zeta)))=\kappa(\zeta^{\langle \lambda, N_{F^d/F}(\mu) \rangle})=
\kappa(\mu(s)).
\endequat

\smallskip
Fixons maintenant un sous-groupe de Levi $F^*$-stable $\Lb^*$ d'un sous-groupe 
parabolique de $\Gb^*$ dual de $\Lb$ et supposons que 
$s \in \Lb^{*F^*}_{\sem,\Lambda}$. La s\'erie 
$\nablab(\Lb,F,(s)_{\Lb^{*F^*}})$ est dite $(\Gb,\Lb)$-{\it r\'eguli\`ere}
(respectivement $(\Gb,\Lb)$-{\it super-r\'eguli\`ere}) si 
$\psi_{\Lb \incl \Pb}(\nablab(\Lb,F,(s)_{\Lb^{*F^*}}))$ est contenue dans 
une s\'erie $(\Gb,I)$-r\'eguli\`ere (respectivement $(\Gb,I)$-super-r\'eguli\`ere) 
de $\nablab(\Tb,W_I,\dot{v}F)$.

\begin{lem}\label{regul}
Avec les notations ci-dessus, $\nablab(\Lb,F,(s)_{\Lb^{*F^*}})$ est 
$(\Gb,\Lb)$-r\'eguli\`ere (respectivement $(\Gb,\Lb)$-super-r\'eguli\`ere) 
si et seulement si $C_{\Gb^*}^\circ(s) \incl \Lb^*$ (respectivement si
et seulement si $C_{\Gb^*}(s) \incl \Lb^*$).
\end{lem}

\begin{proof} 
L'\'egalit\'e \ref{s} montre qu'une coracine 
de $\alpha^\vee$ de $\Gb$ relative \`a $\Tb'$ v\'erifie  
$\theta'(N_{F^d/F}(\alpha^\vee(\zeta)))=1$ si et seulement si elle 
v\'erifie $\alpha^\vee(s)=1$ (ici, $Y(\Tb')$ a \'et\'e identifi\'e 
avec $X(\Tb^{\prime *})$). Donc $\nablab(\Lb,F,(s)_{\Lb^{*F^*}})$ est $(\Gb,\Lb)$-r\'eguli\`ere si et seulement si $C_{\Gb^*}^\circ(s) \incl \Lb^*$.

D'autre part, la m\^eme \'egalit\'e \ref{s} montre qu'un \'el\'ement 
$w \in N_\Gb(\Tb')/\Tb'$ stabilise le caract\`ere $\theta'$ si et
seulement si son correspondant $w^*$ dans $N_{\Gb^*}(\Tb^{\prime *})/\Tb^{\prime *}$ 
stabilise $s$. Donc $\nablab(\Lb,F,(s)_{\Lb^{*F^*}})$ est 
$(\Gb,\Lb)$-super-r\'eguli\`ere 
si et seulement si $C_{\Gb^*}(s) \incl \Lb^*$.
\end{proof}

Avant d'\'enoncer un \'equivalent du th\'eor\`eme \ref{i 0}, nous avons besoin
de quelques notations. Si $s \in \Lb_{\sem,\Lambda}^{*F^*}$, nous noterons
$\FC_s^\Lb$\index{fsl@$\FC_s^\Lb$}
 le syst\`eme local sur la vari\'et\'e $\Xb_\Pb^\Gb$ associ\'e
au $\Lambda\Lb^F$-module $\Lambda\Lb^Fe_s^{\Lb^F}$. Il r\'esulte 
alors imm\'ediatement du lemme \ref{regul} que le th\'eor\`eme 
\ref{i 0} est \'equivalent au th\'eor\`eme suivant~:

\begin{theo}\label{i 0 '}
Soit $s \in \Lb_{\sem,\Lambda}^{*F^*}$. Si $C_{\Gb^*}^\circ(s) \incl \Lb^*$, alors
le morphisme canonique de complexes
$(j_\Pb^\Gb)_! \FC_s^\Lb \isom \RRR(j_\Pb^\Gb)_* \FC_s^\Lb$
est un isomorphisme.
\end{theo}

De m\^eme, d'apr\`es le lemme \ref{regul}, le th\'eor\`eme 
B est \'equivalent au th\'eor\`eme
suivant, conjectur\'e par Brou\'e \cite[p.61]{B}~:

\bigskip

\noindent{\bf Th\'eor\`eme B'.} {\it  Soit $s \in \Lb_{\sem,\Lambda}^{*F^*}$ tel
que $C_{\Gb^*}(s) \incl \Lb^*$. Alors, le foncteur $\RC_{\Lb \incl \Pb}^\Gb$
induit une \'equivalence de Morita entre les $\Lambda$-alg\`ebres
$\Lambda\Lb^F e_s^{\Lb^F}$ et $\Lambda\Gb^Fe_s^{\Gb^F}$.

Plus pr\'ecis\'ement, le complexe $\RC_{\Lb \incl \Pb}^\Gb\Lambda\Lb^F e_s^{\Lb^F}$
n'a de la cohomologie qu'en degr\'e $r=\dim \Vb -\dim \Vb \cap \lexp{F}{\Vb}$
et son $r$-i\`eme bimodule de cohomologie induit l'\'equivalence de Morita
d\'ecrite.}

\subsection{D\'ecomposition de Jordan} 

Rappelons tout d'abord que si $s \in \Gb^{*F^*}_{\sem,\Lambda}$ est
central, alors il existe un caract\`ere lin\'eaire $\hat{s}$ 
de $\Gb^F$ tel que l'automorphisme d'alg\`ebre de $\Lambda \Gb^F$
donn\'e par $g\mapsto \hat{s}(g)g$ se restreint en un isomorphisme
$$
\Lambda\Gb^Fe_s^{\Gb^F}\isom \Lambda\Gb^Fe_1^{\Gb^F}.
\leqno{(\#)_\Gb}$$

Revenons au cas d'un \'el\'ement g\'en\'eral $s \in \Gb^{*F^*}_{\sem,\Lambda}$.
On note $\Lb^*(s)$ le sous-groupe de Levi minimal de $\Gb^*$ contenant
$C_{\Gb^*}(s)$ (il est $F^*$-stable). On note $\Lb(s)$ un sous-groupe
de Levi $F$-stable de $\Gb$ dual de $\Lb^*(s)$.
Le th\'eor\`eme B' fournit une \'equivalence de Morita entre
$\Lambda\Lb(s)^F e_s^{\Lb(s)^F}$ et $\Lambda\Gb^Fe_s^{\Gb^F}$.
Ceci ram\`ene l'\'etude des blocs,
\`a \'equivalence de Morita pr\`es, \`a
l'\'etude de blocs associ\'es \`a un \'el\'ement semi-simple quasi-isol\'e
({\it ie}, \`a un \'el\'ement $s\in \Gb^*_{\sem}$ tel que
$C_{\Gb^*}(s)$ n'est contenu
dans aucun sous-groupe de Levi strict).
Pour un groupe \`a centre connexe, un \'el\'ement semi-simple
quasi-isol\'e du dual est un \'el\'ement isol\'e ({\it ie}, un
\'el\'ement $s\in \Gb^*_{\sem}$ tel que $C_{\Gb^*}^\circ(s)$ n'est contenu
dans aucun sous-groupe de Levi strict).

\begin{theo}
\label{jordan}
Si $C_{\Gb^*}(s)$ est un sous-groupe de Levi de $\Gb^*$, alors
pour tout sous-groupe parabolique $\Pb(s)$ de $\Gb$ de compl\'ement de Levi
$\Lb(s)$, le foncteur
$\RC_{\Lb(s)\incl \Pb(s)}^{\Gb}(\hat{s}\otimes_\Lambda -)$
induit une \'equivalence de Morita entre 
$\Lambda \Lb(s)^F e_1^{\Lb(s)^F}$ et
$\Lambda\Gb^Fe_s^{\Gb^F}$.
\end{theo}

Deux raisons font que $C_{\Gb^*}(s)$ n'est pas n\'ecessairement un
sous-groupe de Levi de $\Gb^*$. D'abord, $C_{\Gb^*}(s)$ peut ne
pas \^etre connexe (cela ne peut \^etre arriver que si $Z(\Gb)$ n'est
pas connexe). Ensuite, le groupe $C_{\Gb^*}^\circ(s)$ peut aussi ne
pas \^etre un sous-groupe de Levi de $\Gb^*$ (cela ne peut arriver que
si l'ordre de $s$ est divisible par un nombre premier mauvais pour
$\Gb$).

Plus pr\'ecis\'ement, soit $\pi_1$ (respectivement $\pi_2$)
l'ensemble des nombres premiers qui divisent $|Z(\Gb)/Z(\Gb)^\circ|$
(respectivement qui sont mauvais pour $\Gb$). On pose
$\pi=\pi_1\cup\pi_2$.
Si $s$ est un $\pi'$-\'el\'ement, alors $C_{\Gb^*}(s)$ est un sous-groupe
de Levi de $\Gb^*$ et le th\'eor\`eme \ref{jordan} s'applique.

Pour un groupe de type $A$ \`a centre connexe (par exemple
un groupe lin\'eaire ou unitaire), alors $\pi=\emptyset$~: on a d\'emontr\'e
la d\'ecomposition de Jordan des blocs.

Pour un groupe de type $B$, $C$ ou $D$, alors $\pi=\{2\}$. On a en
particulier une d\'ecomposition de Jordan des blocs pour $\ell=2$ et
$\Lambda=R/\lG^n$.

\begin{rem}
On ne dispose pas de construction (g\'eom\'etrique) d'un
foncteur qui pourrait induire une \'equivalence fournissant une
d\'ecomposition de Jordan g\'en\'erale.
\end{rem}

\section{Appendice~: quelques lemmes g\'eom\'etriques}
\label{appendice}
Nous rassemblons ici plusieurs r\'esultats
g\'eom\'etriques utiles dans le reste du texte. Plus
pr\'ecis\'ement, nous donnons quelques propri\'et\'es des foncteurs entre
cat\'egories de $\Lambda$-faisceaux constructibles (images directes ou inverses,
\`a support propre ou non). La plupart de ces propri\'et\'es est bien connue
(formules d'adjonction, de K\"unneth, diviseurs lisses \`a croisements normaux).
Nous avons tout de m\^eme inclus une preuve lorsque nous n'avons pas trouv\'e de
r\'ef\'erence satisfaisante. Les lemmes \ref{quasi affine} et
\ref{quasi lemme}, qui semblent
\^etre nouveaux, sont une variation sur un th\`eme classique.

\subsection{Autour de la formule de K\"unneth}
Soit $f : \Xb_1 \to \Xb_2$ un \mor de vari\'et\'es alg\'ebriques, soit
$A$ une $\Lambda$-alg\`ebre finie
et soit $\XC_1$ (\resp $\XC_2$) un objet
de $D^b_{A^\opp}(\Xb_1)$ (\resp de $D^b_A(\Xb_2)$).
Le lemme suivant est un cas particulier de la formule
de K\"unneth (cf par exemple \cite[chapitre VI, remarque 8.14]{milne}).

\begin{lem}\label{!}
Le morphisme canonique
$(\RRR f_!\XC_1)\otimes_A^\LLL \XC_2
 \isom
\RRR f_! (\XC_1 \otimes_A^\LLL f^* \XC_2)$
 est un quasi-isomorphisme
bifonctoriel.
\end{lem}

Le lemme qui suit est lui aussi classique.

\begin{lem}\label{*}
Si $\XC_2$ est \`a cohomologie localement constante, alors le morphisme
canonique 
$(\RRR f_*\XC_1)\otimes_A^\LLL \XC_2
 \isom
\RRR f_* (\XC_1 \otimes_A^\LLL f^* \XC_2)$
est un quasi-isomorphisme bifonctoriel.
\end{lem}

\begin{proof}
Il suffit de d\'emontrer le lemme pour $\Lambda=R/\lG^n$
avec $n>0$. Il suffit aussi de traiter le cas o\`u $\XC_2$
est un faisceau.
Le faisceau $\XC_2$ est localement constant et
l'isomorphisme est de nature locale, donc on peut
supposer, en utilisant le th\'eor\`eme de changement de base lisse
\cite[chapitre VI, th\'eor\`eme VI.4.1]{milne} pour le morphisme $f$, qu'il est constant.
Dans ce cas, le lemme \ref{*} est clair.
\end{proof}

\subsection{Quasi-affinit\'e}
Soit $\Xb$ une vari\'et\'e et $\FC$ un faisceau constructible
sur $X$.
Soit $D$ le foncteur de dualit\'e de Verdier,
$D=\RRR{\HC}om^\bullet(-,\pi_{\Xb}^!\Lambda)$.

Consid\'erons le morphisme canonique de complexes
$\can_\FC : \RRR(\pi_\Xb)_! \FC \to \RRR(\pi_\Xb)_*\FC$\index{canf@$\can_\FC$}.
Le lemme suivant montre quel parti on peut tirer du fait que
ce morphisme est un isomorphisme.

\begin{lem}\label{quasi affine}
Supposons que $\Xb$ est quasi-affine, purement de dimension $d$,
que $\Lambda$ est un corps et que $D(\FC)$ a ses faisceaux de
cohomologie nuls en dehors du degr\'e $-2d$.

Si $\can_\FC$
est un isomorphisme, alors, $\can_{\HC^{-2d}D(\FC)}$ est un isomorphisme et
$H^i_c(\Xb,\FC)=0$ pour $i\not=d$.
\end{lem}

\begin{proof}
Soit $j:\Xb\to\Yb$ une immersion ouverte avec $\Yb$ affine.
On a un diagramme commutatif, o\`u les fl\`eches sont les morphismes
canoniques
$$\xymatrix{
\RRR(\pi_\Yb)_! j_!\FC\ar[rr]^{\can_{\FC}} \ar[dr] & &
  \RRR(\pi_\Yb)_* \RRR j_*\FC \\
& \RRR(\pi_\Yb)_* j_!\FC \ar[ur]
}$$
Puisque $\Yb$ est affine, on a 
$H^i(\Yb,j_!\FC)=0$ pour $i\not\in [0\ldots d]$, donc 
$H^i_c(\Xb,\FC)=0$ pour $i\not\in [0\ldots d]$.

On a (dualit\'e de Verdier=adjonction entre $\RRR(\pi_\Xb)_!$ et
$\pi_{\Xb}^!$)
$$H^i_c(\Xb,\GC)\simeq H^{2d-i}(\Xb,\FC)^*$$
o\`u $\GC=\HC^{-2d}D(\FC)$.
Par cons\'equent, la proposition r\'esultera de ce que
$\can_{\GC}$ est un isomorphisme.

\medskip
Nous allons d\'eduire ceci d'un r\'esultat plus g\'en\'eral.
Soit $f:\Xb_1\to\Xb_2$ un morphisme de vari\'et\'es. Soit
$C_i\in D^b_\Lambda(\Xb_i)$.

Nous allons maintenant rappeler une construction classique,
qui s'effectue en \'etudiant s\'epar\'ement le cas o\`u $f$ est
propre et le cas o\`u $f$ est une immersion ouverte
(cf \cite[exercice III.9 (iii)]{KaScha} pour le cas de la topologie classique).

On a un diagramme commutatif

$$\xymatrix{
\RRR f_!\RRR{\HC}om^\bullet(C_1,f^!C_2) \ar[rr]\ar[d] &&
 \RRR{\HC}om^\bullet(\RRR f_* C_1,C_2) \ar[d] \\
\RRR f_*\RRR{\HC}om^\bullet(C_1,f^!C_2) \ar[rr]_\simeq &&
 \RRR{\HC}om^\bullet(\RRR f_! C_1,C_2) 
}$$
o\`u les fl\`eches verticales proviennent des morphismes canoniques
$\RRR f_!\to \RRR f_*$, la seconde fl\`eche horizontale provient
de  l'adjonction entre $\RRR f_!$ et $f^!$.

\smallskip
On d\'eduit du diagramme pr\'ec\'edent un diagramme commutatif
$$\xymatrix{
\RRR(\pi_{\Xb})_! D(\FC) \ar[d]_{\can_{D(\FC)}}\ar[rr]^\simeq &&
 (\RRR(\pi_{\Xb})_*\FC)^*\ar[d]^{\can_\FC^*} \\
\RRR(\pi_{\Xb})_* D(\FC)\ar[rr]_\simeq &&
 (\RRR(\pi_{\Xb})_!\FC)^*
}$$

Par cons\'equent, $\can_{D(\FC)}$ est un isomorphisme, donc
$\can_{\GC}$ aussi.
\end{proof}

\begin{rem}
\label{haastert}
Dans \cite{Haastert}, Haastert d\'emontre ce r\'esultat en 
supposant que $\can_{D(\FC)}$ est un isomorphisme.
\end{rem}

\begin{rem}
Si, dans le lemme \ref{quasi affine}, on suppose que $\Xb$ est
affine et non pas seulement quasi-affine, alors la conclusion est
bien connue.
\end{rem}

Passons maintenant \`a l'application qui nous est utile.

Soit $G$ un groupe fini agissant sur $\Xb$.

\begin{lem}\label{quasi lemme}
Supposons la vari\'et\'e $\Xb$ quasi-affine, lisse, purement de dimension $d$,
et supposons que $\Lambda$ est un corps. Supposons de plus que les
stabilisateurs de points de $\Xb$ dans $G$ sont d'ordre inversible
dans $\Lambda$.

Soit $M$ un $\Lambda G$-module tel que le morphisme
canonique $\RRR\Gamma_c(\Xb)\otimes_{\Lambda G}^\LLL M \to
\RRR\Gamma(\Xb)\otimes_{\Lambda G}^\LLL M$ est un isomorphisme.
Alors, le complexe
$\RRR\Gamma_c(\Xb)\otimes_{\Lambda G}^\LLL M$ n'a de l'homologie
qu'en degr\'e $d$.
\end{lem}

\begin{proof}
Soit $\pi:\Xb\to\Xb/G$ le morphisme quotient et
$\FC=\pi_* \Lambda_\Xb \otimes_{\Lambda G} M_{\Xb/G}$.
D'apr\`es le lemme \ref{lemme foncteur}, le morphisme
$\can_\FC$ est un isomorphisme.

Un calcul simple utilisant les adjonctions entre $\pi_*$ et
$\pi^!$ et entre $\otimes$ et $\Hom$ fournit
$$D(\FC)\simeq D(\pi_*\Lambda_\Xb\otimes_{\Lambda G}^\LLL M_{\Xb/G})\simeq
\pi_* D(\Lambda_{\Xb})\otimes_{\Lambda G}^\LLL (M^*)^\opp_{\Xb/G}.$$
Puisque $\Xb$ est lisse, on a $D(\Lambda_\Xb)\simeq \Lambda_\Xb[2d]$.
L'hypoth\`ese sur les stabilisateurs de points de $\Xb$ permet d'appliquer le
lemme \ref{lemme foncteur} et donc de conclure que
$$D(\FC)\simeq \pi_*\Lambda_{\Xb}\otimes_{\Lambda G}(M^*)^\opp_{\Xb/G}[2d].$$

Puisque $\Xb/G$ est quasi-affine,
le lemme \ref{quasi affine} fournit le r\'esultat.
\end{proof}

\subsection{Diviseurs \`a croisements normaux}
Soit $\overline{\Xb}$ une vari\'et\'e lisse et soit $\Xb$ un ouvert de
$\overline{\Xb}$. Nous noterons $i : \Xb \injto \overline{\Xb}$ l'immersion
ouverte canonique. Nous supposons que le compl\'ementaire
de $\Xb$ dans $\overline{\Xb}$ est une r\'eunion finie
        $$\bigcup_{i \in I} ~\Db_i,$$
o\`u les $\Db_i$ sont des diviseurs lisses \`a croisements normaux.

\bigskip

Par ailleurs, si $J \incl I$, nous noterons
$\Xb_J=(\bigcap_{i \in J} \Db_i)-(\bigcup_{i \in I\setminus J} \Db_i)$.
Selon les conventions usuelles, $\Xb_\vide=\Xb$ et
$\Xb_I=\bigcap_{i \in I} \Db_i$. Notons $i_J : \Xb_J \injto \overline{\Xb}$
l'immersion localement ferm\'ee canonique.

Le r\'esultat suivant est classique (cf par exemple
\cite[(1.3.3.2) p. 255]{SGA4.5})~:

\begin{lem}\label{DCN}
Si $J \incl I$ et si $n \in \NM$, alors
$$i_J^* \RRR^n i_* \Lambda_\Xb \simeq
 \Lambda_{\Xb_J}^{\oplus {|J| \choose n}}.$$
\end{lem}

\printindex
\end{document}